# Coorbit Theory for Coefficients in Weighted Lebesgue Spaces, and its Application to the Wavelet Transform and the Short-Time Fourier Transform

Jan Zimmermann

# Contents

# Introduction

Coorbit theory was developed in the late 1980s by Feichtinger and Gröchenig [14; 15; 16; 21] and utilizes tools from functional and harmonic analysis. Starting with some unitary group representation, it describes the construction of a family of smoothness spaces that are invariant under said representation. These so called *Coorbit spaces* are characterized through decay conditions of the *voice transform* that is associated to the representation.

The voice transform generalizes many transforms that are used in science and engineering. Famous examples are the short-time Fourier transform and the wavelet transform, which are for instance used in time-frequency analysis [22] and signal processing [27]. Another example is the shearlet transform, a relative of the wavelet transform from image processing [5; 24].

Coorbit theory describes how smoothness spaces can be constructed that are tailored to the respective voice transform and therefore to the respective application, making coorbit theory a useful tool in a lot of contexts. In doing so it recovers some smoothness spaces that are already known and allows a new interpretation of these. The coorbit spaces that are associated to the wavelet transform can be identified with the homogeneous Besov spaces, while the coorbit spaces of the short-time Fourier transform are exactly the modulation spaces.

One of the main theorems of coorbit theory states that these coorbit spaces can be described equivalently in terms of *atomic decompositions* and *Banach frames*, i.e. they can be described by sequence spaces in a way that associates their elements to sequences of coefficients. Both concepts generalize the notion of frames in Hilbert spaces [8], which itself is a very successful tool in harmonic analysis. For many applications, such discrete characterisations are crucial, so coorbit theory provides spaces with the required properties.

The voice transform of a given vector (where the vector can be a signal, function, etc. depending on the context) is a continuous function that is defined on the underlying group. Its values can be understood as a continuous family of coefficients. Each coorbit space consists of all vectors for which the voice transform belongs to a certain function space that 'measures' the decay of the voice transform. The discrete characterisation via atomic decompositions and Banach frames is done similarly by considering certain sequence space norms of discrete families of coefficients, that are samplings of the voice transform.

In general, it is possible to construct coorbit spaces with respect to any *solid Banach*

*function space* that is defined on the underlying group (note that extensions to quasi-Banach spaces do exist [31]). We restrict ourself to the smaller class of *weighted Lebesgue spaces*, which are less abstract but still contain many interesting examples. However, by that choice we exclude the important class of *mixed norm Lebesgue spaces*, where different exponents for different parameters are combined into a single norm.

**Principle of Construction**

A unitary representation $\pi : \mathcal{G} \to \mathcal{U}(\mathcal{H})$ is a unitary group action of a locally compact group $\mathcal{G}$ on a Hilbert space $\mathcal{H}$ with some additional properties. The voice transform of $\pi$ is defined for $f, g \in \mathcal{H}$ as the function

$$\mathcal{V}_g f : \mathcal{G} \to \mathbb{C}, \quad \mathcal{V}_g f(x) = \langle f, \pi(x)g \rangle.$$

If $g$ satisfies certain decay conditions, it is possible to extend $\mathcal{V}_g$ from a transform on $\mathcal{H}$ to a transform on some suitable distribution space. In that case, $g$ is called an *analyzing vector*.

Now the coorbit space $\mathcal{C}o_m^p$ consists of all distributions $f$ for which $\mathcal{V}_g f$ is contained in the weighted Lebesgue space $L_m^p(\mathcal{G})$. It becomes a Banach space if we equip it with the norm $\|f\|_{\mathcal{C}o_m^p} = \|\mathcal{V}_g f\|_{L_m^p}$. The voice transform is then by definition an isometric isomorphism $\mathcal{V}_g : \mathcal{C}o_m^p \to L_m^p(\mathcal{G})$.

For the discrete description of $\mathcal{C}o_m^p$ we need a family of points $(x_i)_{i \in I} \subset \mathcal{G}$ that is spread uniformly in $\mathcal{G}$ in a geometric sense. Then, provided that $g$ satisfies some additional smoothness conditions, any $f \in \mathcal{C}o_m^p$ can be written as a series

$$f = \sum_{i \in I} c_i\, \pi(x_i) g,$$

where the sequence $(c_i)_{i \in I}$ lies in a weighted sequence space $\ell_{\widetilde{m}}^p(I)$ that is associated to $L_m^p(\mathcal{G})$. Such a series representation of $f$ is an *atomic decomposition*. The family $(\pi(x_i)g)_{i \in I}$ is also a *Banach frame*, which essentially means that

$$\left\| \langle f, \pi(x_i)g \rangle_{i \in I} \right\|_{\ell_{\widetilde{m}}^p}, \quad f \in \mathcal{C}o_m^p$$

defines an equivalent norm on $\mathcal{C}o_m^p$ (see chapter 5 for a precise definition). Both notions are closely related, and the proofs that $(\pi(x_i)g)_{i \in I}$ defines an atomic decomposition respectively a Banach frame are similar. We will only consider Banach frames in this exposition.

The above constructions are essentially independent from the chosen analyzing vector $g \in \mathcal{H} \backslash \{0\}$, as long as it satisfies the decay conditions. Of course, it is only possible to carry out these constructions if at least one such $g$ exists at all.

**Organization**

The second chapter is devoted to an introduction to the theories of integration and representations on locally compact groups. This includes the description of weighted Lebesgue spaces and their convolution relations, as well as a notion for Banach space valued integrals. In the last two sections of that chapter, we consider unitary representations and their voice transforms, with a special focus on the reconstruction and reproducing formulas that arise from the orthogonality relations of these voice transforms.

In the third chapter we carry out the construction of the coorbit spaces with respect to a unitary representation. For that, we extend the voice transform to a suitable space of distributions and define the coorbit spaces as explained above. We will also take a look at some important properties of the coorbit spaces.

We use the fourth chapter to examine two examples in more detail, namely the wavelet transform and the short-time Fourier transform. We first show that both indeed arise from unitary representations, and then apply the theory of chapter three to both cases. Additionally, we give sufficient conditions for a vector (i.e. a wavelet resp. Gabor window) to be an analyzing vector, that is, to be localized well enough to characterize the coorbit spaces.

In the fifth chapter we discretize the coorbit spaces through Banach frames. By requiring the used analyzing vector to satisfy a certain smoothness condition, we are able to approximate the voice transform effectively by discrete samples, so we can compare its continuous weighted Lebesgue norm with a discrete sequence norm.

The sixth chapter contains an application of the results regarding Banach frames to the wavelet transform and the short-time Fourier transform. We also examine under which conditions an analyzing vector satisfies the required smoothness condition to be suitable for Banach frame constructions in the respective context.

# CHAPTER 2

## Integration and Representations on Locally Compact Groups

Before we started with coorbit theory, we introduce the basic notions it is based upon. This covers essentially the theories of integration and unitary representations on locally compact groups.

We first give an introduction to the theory of integration on locally compact groups in section 2.1. After that, we investigate in section 2.2 special classes of weighted Lebesgue spaces and their convolution relations. Section 2.3 serves as a short overview on integrals with values in Banach spaces, as we will encounter those occasionally.

We are then able to describe unitary representations and how they give rise to transforms like the wavelet transform or the short-time Fourier transform in sections 2.4 and 2.5. In section 2.6, we take a closer look at how the orthogonality relations of these transforms can be expressed through Hilbert space valued integral formulas and convolution formulas.

## 2.1 Locally Compact Groups and the Haar Measure

We begin by gathering all the necessary definitions and properties of locally compact groups and their theory of integration. This includes Lebesgue spaces as well as a notion of convolution. Since the used measure is canonical and essentially unique, the mentioned terms are unique themselves.

For a more detailed introduction to the theory of locally compact groups and complete proofs of all statements, we refer the reader to [17, Ch. 2].

**Definition 2.1.1.** A *locally compact group* is a group $\mathcal{G}$ which is equipped with a topology such that the following is true:

  i. $\mathcal{G}$ is a topological group, i.e. the group operation and inversion is continuous.

  ii. The topology is locally compact, i.e. for every point $x \in \mathcal{G}$ and every neighbourhood $U \subset \mathcal{G}$ of $x$ there exists a compact neighbourhood $K \subset U$ of $x$. In other words, every point in $\mathcal{G}$ has a neighbourhood basis of compact sets.

  iii. The topology is Hausdorff, i.e. two distinct points in $\mathcal{G}$ are separated by disjoint neighbourhoods.

A special property of locally compact groups is the existence of a unique *Haar measure*. Such a measure is defined on the group's Borel $\sigma$-algebra and is invariant under the group action, that is, it satisfies $\mu_L(xM) = \mu_L(M)$ or $\mu_R(Mx) = \mu_R(M)$ for $x \in \mathcal{G}$ and measurable $M \subset \mathcal{G}$. If $\mathcal{G}$ is non abelian, the two equations might lead to different measures $\mu_L \neq \mu_R$, so $\mu_L$ and $\mu_R$ are more precisely called a *left Haar measure* respectively a *right Haar measure*. In the following, we will only use left Haar measures if not stated otherwise, and just write $\mu = \mu_L$.

**Theorem 2.1.2** (Existence and uniqueness of the Haar measure). *Let $\mathcal{G}$ be a locally compact group. Then there exists a Radon measure $\mu$ for which $\mu(xM) = \mu(M)$ is true for all $x \in \mathcal{G}$ and all Borel sets $M \subset \mathcal{G}$. Such a $\mu$ is called a* Haar measure.

*The Haar measure is unique up to a positive multiplicative constant. That is, if $\mu$ and $\nu$ are two non-trivial Haar measures on $\mathcal{G}$, then there exists some $c > 0$ such that $\mu = c\nu$.*

A simple example of a locally compact group is the additive space $(\mathbb{R}^d, +)$. For this group, the usual Lebesgue measure is a Haar measure. The latter can be seen as a generalization of the former.

Thanks to the Haar measure, there is a canonical notion of integration on locally compact groups. In the following, we will simply write $dx$ for the integration with respect to the Haar measure $d\mu(x)$ if there is no risk of confusion.

The invariance of a Haar measure $\mu$ on the locally compact group $\mathcal{G}$ leads to several computational rules for the substitutions $x \mapsto yx$, $x \mapsto xy$ and $x \mapsto x^{-1}$ regarding the integral $\int_{\mathcal{G}} dx$ and a constant $y \in \mathcal{G}$.

For the first substitution, we consider the indicator function

$$\chi_M(x) = \begin{cases} 1, & x \in M \\ 0, & x \notin M \end{cases}$$

of some measurable set $M \subset \mathcal{G}$. We have

$$\int_{\mathcal{G}} \chi_M(x)\, dx = \mu(M) = \mu(y^{-1}M) = \int_{\mathcal{G}} \chi_{y^{-1}M}(x)\, dx = \int_{\mathcal{G}} \chi_M(yx)\, dx.$$

Through linear combinations and a density argument, this calculation can be carried over to all integrable functions, hence we have $\int_{\mathcal{G}} f(x)\, dx = \int_{\mathcal{G}} f(yx)\, dx$ for all integrable $f : \mathcal{G} \to \mathbb{C}$. Therefore, the integral is invariant under the substitution $x \mapsto yx$, i.e. it is left invariant. This fact can be expressed heuristically by the rule $d(yx) = dx$.

The left invariance of the integral is not only necessary for a measure to be a Haar measure, but also sufficient. That is, if $\nu$ is any Radon measure on $\mathcal{G}$ that induces a left invariant integral in the above sense, then $\nu$ is a Haar measure, as for such $\nu$ we have

$$\nu(yM) = \int_{\mathcal{G}} \chi_{yM}(x)\, d\nu(x) = \int_{\mathcal{G}} \chi_M(y^{-1}x)\, d\nu(x) = \int_{\mathcal{G}} \chi_M(x)\, d\nu(x) = \nu(M)$$

for all measurable $M \subset \mathcal{G}$ and $y \in \mathcal{G}$.

The right sided substitution is more complicated. We would like to compare $\mu(M)$ with $\mu(My)$ as before, but these terms are not equal in general. However, the mapping $M \mapsto \mu(My)$ is again a non-trivial (left) Haar measure, which is a positive multiple of $\mu$ by the uniqueness of theorem 2.1.2. Thus, there exists $\Delta(y) \in (0, \infty)$ such that $\mu(My) = \Delta(y)\mu(M)$ is true for all measurable $M \subset \mathcal{G}$. It follows

$$\int_{\mathcal{G}} \chi_M(x)\, dx = \mu(M) = \Delta(y)\mu(My^{-1}) = \Delta(y) \int_{\mathcal{G}} \chi_M(xy)\, dx,$$

which again carries over to all integrable functions and leads to the heuristic rule $d(xy) = \Delta(y) \, dx$.

The map $\Delta : \mathcal{G} \to (0, \infty)$, $y \mapsto \Delta(y)$ is a continuous group homomorphism (with respect to the multiplicative group $(0, \infty)$) and is called the *Haar modulus* or *modular function* of $\mathcal{G}$. It is independent of the exact choice of Haar measure, so it describes an intrinsic property of $\mathcal{G}$. To some extent, $\Delta$ quantifies the difference between left and right Haar measures. Both coincide if and only if $\Delta \equiv 1$. In that case, the group is called *unimodular*. Abelian groups are obviously unimodular. The same is true for compact and discrete groups, as in both cases the only continuous group homomorphism $\mathcal{G} \to (0, \infty)$ is the trivial map $\Delta \equiv 1$.

Using the modular function, also the last substitution can be examined. For integrable $f : \mathcal{G} \to \mathbb{C}$ we have

$$\int_\mathcal{G} f(x^{-1}) \, dx = \int_\mathcal{G} f(x) \Delta(x^{-1}) \, dx,$$

which can be described by the heuristic rule $d(x^{-1}) = \Delta(x^{-1}) \, dx$. For a proof of this statement we refer to [17, Prop 2.31].

Note that the Haar modulus is uniquely determined by either of the relations $d(xy) = \Delta(y) \, dx$ and $d(x^{-1}) = \Delta(x^{-1}) \, dx$. This can be useful to determine the modulus of a given locally compact group $\mathcal{G}$.

*Example 2.1.3.*

i. We consider the group of multiplicative real numbers $\mathbb{R}^* = \mathbb{R}\backslash\{0\}$. It is abelian and therefore unimodular. Furthermore, if $f : \mathbb{R}^* \to \mathbb{C}$ is a measurable function and $a \in \mathbb{R}^*$, we have

$$\int_{\mathbb{R}^*} f(x) \frac{d\lambda(x)}{|x|} = \int_{\mathbb{R}^*} f(ax) \frac{d\lambda(ax)}{|ax|} = \int_{\mathbb{R}^*} f(ax) \frac{d\lambda(x)}{|x|},$$

where $\lambda$ is the usual Lebesgue measure.

The integral with respect to $\frac{d\lambda(x)}{|x|}$ is apparently invariant under multiplication with $a$, which is exactly the (left-sided) group action of $\mathbb{R}^*$ on itself. Therefore, $\frac{d\lambda(x)}{|x|}$ defines a measure that is invariant under the group operation, that is, it defines a Haar measure on $\mathbb{R}^*$.

ii. Similar to the first example, the multiplicative complex numbers $\mathbb{C}^* = \mathbb{C}\backslash\{0\}$ have $\frac{d\lambda(x,y)}{x^2+y^2}$ as a Haar measure.

iii. An example of a non-unimodular group is the *affine group* $\mathcal{A}\!f\!f$. It is given by the topological space $\mathbb{R} \times \mathbb{R}^*$ which is equipped with the group operation

$$(b_1, a_1)(b_2, a_2) = (b_1 + a_1 b_2, a_1 a_2).$$

This makes $\mathcal{A}\!f\!f$ a non-abelian locally compact group.

In order to find a Haar measure on $\mathcal{A}\!f\!f$, we consider the integral

$$\int_{\mathcal{A}\!f\!f} f(x,y) \frac{d\lambda(x,y)}{y^2}$$

with measurable $f : \mathcal{A}\!f\!f \to \mathbb{C}$ and the two-dimensional Lebesgue measure $d\lambda(x,y)$. We apply the substitution $(x,y) = (b,a)(u,v) = (b+au, av)$ for some constant $(b,a) \in$ $\mathcal{A}\textit{ff}$ and variables $(x, y), (u, v) \in \mathcal{A}\textit{ff}$. The Jacobian of the map $(u, v) \mapsto (b+au, av)$ is $a^2$. Therefore we have

$$\int_{\mathcal{G}} f(x,y) \frac{d\lambda(x,y)}{y^2} = \int_{\mathcal{A}\textit{ff}} f((b,a)(u,v)) |a|^2 \frac{d\lambda(u,v)}{(av)^2}$$
$$= \int_{\mathcal{A}\textit{ff}} f((b,a)(u,v)) \frac{d\lambda(u,v)}{v^2}.$$

We see that the integral with respect to $d\mu(x,y) = \frac{d\lambda(x,y)}{y^2}$ is invariant under the left group action, so $\mu$ defines a (left) Haar measure on $\mathcal{A}\textit{ff}$.

In general, a right Haar measure $\mu_R$ can be defined in terms of an already known left Haar measure $\mu_L$ by $d\mu_R(z) = d\mu_L(z^{-1})$. On the affine group, the inverse of $(x,y) \in \mathcal{A}\textit{ff}$ is given by $(-x/y, 1/y)$. Thus, we can define a right Haar measure on $\mathcal{A}\textit{ff}$ by

$$d\mu_R(x,y) = d\mu((x,y)^{-1}) = d\mu(-x/y, 1/y) = \frac{d\lambda(-x/y, 1/y)}{(1/y)^2}.$$

The Jacobian of $(x,y) \mapsto (-x/y, 1/y)$ is $1/y^3$, so we have

$$d\mu_R(x,y) = \frac{1}{|y|^3} \frac{d\lambda(x,y)}{(1/y)^2} = \frac{d\lambda(x,y)}{|y|}.$$

The modular function is now determined by $d\mu((x,y)^{-1}) = \Delta(x,y) d\mu(x,y)$. The above calculations show that $\Delta(x,y) = |y|$.

◁

As we now have a well-defined notion of integration on the locally compact group $\mathcal{G}$, we are able to define the Lebesgue spaces. As usual, for measurable $F : \mathcal{G} \to \mathbb{C}$ we define the norms

$$\|F\|_{L^p} = \left( \int_{\mathcal{G}} |F(x)|^p \, dx \right)^{1/p}$$

for $1 \leq p < \infty$ as well as

$$\|F\|_{L^\infty} = \operatorname*{ess\,sup}_{x \in \mathcal{G}} |F(x)|.$$

Then we set for $1 \leq p \leq \infty$

$$L^p(\mathcal{G}) = \{F : \mathcal{G} \to \mathbb{C} \text{ measurable } | \ \|F\|_{L^p} < \infty\},$$

where we identify functions that are equal pointwise almost everywhere. Equipped with their natural norm $\|\cdot\|_{L^p}$, the spaces $L^p(\mathcal{G})$ become Banach spaces. For $p = 2$ the inner product

$$\langle F, G \rangle_{L^2(\mathcal{G})} = \int_{\mathcal{G}} F(x) \overline{G(x)} \, dx$$

makes $L^2(\mathcal{G})$ a Hilbert space.

The spaces $L^p(\mathcal{G})$ are independent from the choice of Haar measure. As different Haar measures are just multiples of each other, so are the associated $L^p$-norms, meaning the associated $L^p$-spaces are equal as vector spaces and isometrically isomorphic as Banach spaces.

Since integration on locally compact groups is defined in terms of measure theory, we are able to make use of all the well-known facts from Lebesgue integration theory. However, there is a technical detail we need to keep in mind. Fubini's theorem is only valid for *σ-finite measure spaces*, that is, measure spaces which can be covered by countably many sets of finite measure. This is why we will always require the used groups to be σ-compact.

**Definition 2.1.4.** A locally compact group is called *σ-compact* if it can be covered by countably many compact subsets.

Compact sets in locally compact groups are measurable and have finite measure. Thus, σ-compact groups are also σ-finite.

We end this section by setting some notation. Let $F : \mathcal{G} \to \mathbb{C}$ be a function. We define two involutions of $F$ by

$$F^\vee(x) = F(x^{-1}), \quad F^\triangledown(x) = \overline{F(x^{-1})}, \quad x \in \mathcal{G}.$$

If $Y$ is some space of functions on $\mathcal{G}$, we write $Y^\vee = \{H^\vee \mid H \in Y\}$. We define the left and right translation operators by

$$L_y F(x) = F(y^{-1}x), \quad R_y F(x) = F(xy), \quad x \in \mathcal{G}$$

for $y \in \mathcal{G}$. We furthermore define for measurable functions $F, G : \mathcal{G} \to \mathbb{C}$ the convolution

$$(F * G)(x) = \int_{\mathcal{G}} F(y) G(y^{-1}x) \, dx, \quad x \in \mathcal{G},$$

provided the integral is well-defined. The convolution is associative but in general not commutative. In the next section, we will take a closer look at sufficient conditions for the existence of the convolution integral, and derive certain Young-like norm inequalities.

## 2.2 Weighted Lebesgue Spaces and their Convolution Relations

In this section, we introduce the *weighted $L^p$-spaces*. These extend the usual $L^p$-spaces by a weight function and allow us to describe the decay of functions $F : \mathcal{G} \to \mathbb{C}$ in terms of finiteness of some Banach space norm.

**Definition 2.2.1.** Let $\Omega$ be a measure space.

a) We call a positive, measurable and locally integrable function $v : \Omega \to (0, \infty)$ a *weight function* or simply a *weight*.

b) Given a weight function $v$ on $\Omega$ and an exponent $1 \leq p \leq \infty$, we define the weighted $L^p$-space
$$L^p_v(\Omega) = \{F : \Omega \to \mathbb{C} \mid Fv \in L^p(\Omega)\}$$
and equip it with the norm $\|F\|_{L^p_v} = \|Fv\|_{L^p}$.

As usual, two functions in $L^p_v(\Omega)$ are ought to be the same if they coincide almost everywhere.

The spaces defined in this way are Banach spaces. This is an immediate consequence of the fact that the map

$$L^p_v(\Omega) \to L^p(\Omega), \quad F \mapsto Fv$$

defines an isometric isomorphism between the weighted and unweighted $L^p$-space. We can also use this map to pull back the duality properties of $L^p$-spaces to their weighted counterparts.

**Lemma 2.2.2.** *Let $1 \leq p < \infty$ and $p^{-1} + q^{-1} = 1$. Then the dual space of $L_v^p(\Omega)$ can be identified with $L_{1/v}^q(\Omega)$ by defining for $G \in L_{1/v}^q(\Omega)$ the functional*

$$L_v^p(\Omega) \ni F \mapsto \int_\Omega F(x) G(x)\, dx \in \mathbb{C}.$$

*Proof.* First, every function $G \in L_{1/v}^q(\Omega)$ defines a functional through the above mapping, as by Hölder's inequality we have

$$\left| \int_\Omega F(x) G(x)\, dx \right| = \left| \int_\Omega F(x) v(x) \frac{G(x)}{v(x)}\, dx \right| \leq \|Fv\|_{L^p} \|G/v\|_{L^q} = \|F\|_{L_v^p} \|G\|_{L_{1/v}^q}$$

for all $F \in L_v^p(\Omega)$. Note that the norm of that functional is exactly $\|G\|_{L_{1/v}^q}$ since Hölder's inequality is sharp.

If $\alpha : L_v^p(\Omega) \to \mathbb{C}$ is an arbitrary functional on $L_v^p(\Omega)$, we can define its counterpart on $L^p(\mathcal{G})$ by

$$\widetilde{\alpha} : L^p(\Omega) \to \mathbb{C}, \quad \widetilde{\alpha}(\widetilde{F}) = \alpha(\widetilde{F}/v).$$

Now $\widetilde{\alpha}$ is represented by some $\widetilde{G} \in L^q(\Omega)$ in the sense that

$$\widetilde{\alpha}(\widetilde{F}) = \int_\Omega \widetilde{F}(x) \widetilde{G}(x)\, dx$$

holds for all $\widetilde{F} \in L^p(\Omega)$. By setting $G = \widetilde{G} v \in L_{1/v}^q(\Omega)$ we have

$$\alpha(F) = \widetilde{\alpha}(Fv) = \int_\Omega F(x) v(x) \widetilde{G}(x)\, dx = \int_\Omega F(x) G(x)\, dx$$

for all $F \in L_v^p(\Omega)$.

Thus, every functional on $L_v^p(\mathcal{G})$ can be represented by some $G \in L_{1/v}^q(\mathcal{G})$ which has the same norm as the functional. This completes the proof. $\square$

If $w, v$ are two weight functions such that $v \geq cw$ for some constant $c > 0$, we have the inequality

$$\|F\|_{L_v^p(\Omega)} = \left( \int_\Omega |F(x)|^p v(x)^p\, dx \right)^{1/p} \geq c \left( \int_\Omega |F(x)|^p w(x)^p\, dx \right)^{1/p} = c\|F\|_{L_w^p(\Omega)}$$

for all measurable $F : \Omega \to \mathbb{C}$. This implies the continuous embedding $L_v^p(\Omega) \hookrightarrow L_w^p(\Omega)$. For weights $v$ that are bounded away from 0 (i.e. $v \geq c > 0$), it particularly follows the continuous embedding $L_v^p(\Omega) \hookrightarrow L^p(\Omega)$.

In coorbit theory, we work with weighted $L^p$-spaces on locally compact groups with respect to their Haar measure, so from now on we assume all weights are defined on some locally compact group $\mathcal{G}$.

We will use the weighted $L^p$-spaces in two different ways:

- The spaces $L_m^p(\mathcal{G})$ for arbitrary $1 \leq p \leq \infty$ and *moderate* weights $m$ describe a wide range of decay conditions through their norms.

- The space $L_w^1(\mathcal{G})$ for a *submultiplicative* weight $w$ contains functions that are localised 'well enough'. These functions help us to handle the more general $L_m^p(\mathcal{G})$-spaces.

We will now describe moderate and submultiplicative weight functions before examining how they are related to each other. We begin with moderate weights.

**Definition 2.2.3.** Let $m : \mathcal{G} \to (0, \infty)$ be a weight function.

a) We call $m$ *left-moderate (with respect to $\alpha$)* if there exists some $\alpha : \mathcal{G} \to (0, \infty)$ such that $m(xy) \leq \alpha(x)m(y)$ for all $x, y \in \mathcal{G}$.

b) We call $m$ *right-moderate (with respect to $\beta$)* if there exists some $\beta : \mathcal{G} \to (0, \infty)$ such that $m(xy) \leq m(x)\beta(y)$ for all $x, y \in \mathcal{G}$.

c) We call $m$ *moderate* if it is left- and right-moderate with respect to some functions $\alpha$ and $\beta$.

A weight function $m$ is left-moderate with respect to some $\alpha$ if and only if

$$\alpha_0(y) = \sup_{x \in \mathcal{G}} \frac{m(yx)}{m(x)}$$

is finite for all $y \in \mathcal{G}$. In that case, $m$ is left-moderate with respect to $\alpha_0$ and $\alpha_0$ is the smallest function $\alpha$ with that property. Similarly, $m$ is right-moderate with respect to some $\beta$ if and only if the function

$$\beta_0(y) = \sup_{x \in \mathcal{G}} \frac{m(xy)}{m(x)}$$

is finite for all $y \in \mathcal{G}$, where $\beta_0$ is the smallest possible choice for such $\beta$.

We need the following technical properties of $\alpha_0$ and $\beta_0$. For the proof we refer to [10, Prop. 1.16].

**Lemma 2.2.4.** *If $m$ is moderate, the functions $\alpha_0$ and $\beta_0$ are measurable and locally bounded, and in particular locally integrable.*

The functions $\alpha_0$ and $\beta_0$ are therefore weight functions in our terms.

The weighted $L^p$-spaces with respect to (left-/right-) moderate weights now have certain properties which are important to us.

**Proposition 2.2.5.** *Let $m$ be a left-moderate weight with respect to the locally integrable function $\alpha : \mathcal{G} \to (0, \infty)$, and let $1 \leq p \leq \infty$.*

(i) *The space $L^p_m(\mathcal{G})$ is left-invariant, i.e. invariant under left translations, and we have the inequality*

$$\left\| L_y F \right\|_{L^p_m} \leq \alpha(y) \|F\|_{L^p_m} \tag{2.1}$$

*for all $y \in \mathcal{G}$ and $F \in L^p_m(\mathcal{G})$.*

(ii) *The space $L^p_m(\mathcal{G})$ is invariant under left-sided convolutions with functions from $L^1_\alpha(\mathcal{G})$ and we have the Young inequality*

$$\|G * F\|_{L^p_m} \leq \|G\|_{L^1_\alpha} \|F\|_{L^p_m} \tag{2.2}$$

*for all $G \in L^1_\alpha(\mathcal{G})$ and $F \in L^p_m(\mathcal{G})$.*

A similar statement is true for right-moderate weights.

**Proposition 2.2.6.** *Let $m$ be a right-moderate weight with respect to the locally integrable function $\beta : \mathcal{G} \to (0, \infty)$, and let $1 \leq p \leq \infty$.*

(i) *The space $L^p_m(\mathcal{G})$ is right-invariant, i.e. invariant under right translations, and we have the inequality*
$$\left\| R_y F \right\|_{L^p_m} \leq \beta(y^{-1}) \Delta(y)^{-1/p} \|F\|_{L^p_m} \tag{2.3}$$
*for all $y \in \mathcal{G}$ and $F \in L^p_m(\mathcal{G})$, where we use the convention $1/\infty = 0$.*

(ii) *The space $L^p_m(\mathcal{G})$ is invariant under right-sided convolutions with functions from $L^1_{\beta^\vee \Delta^{-1/p}}(\mathcal{G})^\vee$ and we have the Young inequality*
$$\left\| F * G^\vee \right\|_{L^p_m} \leq \|F\|_{L^p_m} \|G\|_{L^1_{\beta^\vee \Delta^{-1/p}}} \tag{2.4}$$
*for all $G \in L^1_{\beta^\vee \Delta^{-1/p}}(\mathcal{G})$ and $F \in L^p_m(\mathcal{G})$.*

Because we are interested in translation- and convolution-invariant function spaces, moderate weight functions are our weights of choice. We will only prove the statements regarding right-moderate weights, as the properties of left-moderate weights follow in the same way.

*Proof of proposition 2.2.6.* Suppose $m$ is right-moderate with respect to $\beta$. Then, using the substitution rule $d(xy^{-1}) = \Delta(y)^{-1} dx$ for constant $y \in \mathcal{G}$, we get for $p \in [1, \infty)$ and $F \in L^p_m(\mathcal{G})$

$$\begin{aligned}
\left\| R_y F \right\|_{L^p_m} = \left\| F(\cdot\, y) \right\|_{L^p_m} &= \left( \int_{\mathcal{G}} |F(xy)|^p m(x)^p \, dx \right)^{1/p} \\
&= \left( \int_{\mathcal{G}} |F(x)|^p m(xy^{-1})^p \Delta(y)^{-1} \, dx \right)^{1/p} \\
&\leq \left( \int_{\mathcal{G}} |F(x)|^p m(x)^p \beta(y^{-1})^p \Delta(y)^{-1} \, dx \right)^{1/p} \\
&= \Delta(y)^{-1/p} \beta(y^{-1}) \|F\|_{L^p_m}.
\end{aligned}$$

Similarly, we have for $p = \infty$
$$\left\| R_y F \right\|_{L^\infty_m} = \operatorname*{ess\,sup}_{x \in \mathcal{G}} |F(xy)| m(x) \leq \operatorname*{ess\,sup}_{x \in \mathcal{G}} |F(x)| m(x) \beta(y^{-1}) = \beta(y^{-1}) \|F\|_{L^\infty_m}.$$

By convention $\Delta(y)^{-1/\infty} = \Delta(y)^0 = 1$ holds, thus we have proven the first part of the proposition.

We show the Young inequality by duality. Suppose $G \in L^1_{\beta^\vee \Delta^{-1/p}}(\mathcal{G})$, $F \in L^p_m(\mathcal{G})$ and $H \in L^q_{1/m}(\mathcal{G})$ with $p^{-1} + q^{-1} = 1$. Then

$$\begin{aligned}
\int_{\mathcal{G}} |(F * G^\vee)(x) H(x)| \, dx &\leq \int_{\mathcal{G}} \int_{\mathcal{G}} |F(y) G^\vee(y^{-1} x)| \, dy \, |H(x)| \, dx \\
&= \int_{\mathcal{G}} \int_{\mathcal{G}} |F(y)| |G(x^{-1} y)| |H(x)| \, dy \, dx.
\end{aligned}$$

By applying the substitution $y \mapsto xy$ we get

$$\begin{aligned}
\int_{\mathcal{G}} |(F * G^{\vee})(x) H(x)| \, dx &\leq \int_{\mathcal{G}} \int_{\mathcal{G}} |F(xy)||G(y)||H(x)| \, dy \, dx \\
&= \int_{\mathcal{G}} \int_{\mathcal{G}} |R_y F(x)||H(x)| \, dx \, |G(y)| \, dy \\
&\leq \int_{\mathcal{G}} \|R_y F\|_{L^p_m} \|H\|_{L^q_{1/m}} |G(y)| \, dy \\
&\leq \|F\|_{L^p_m} \|H\|_{L^q_{1/m}} \int_{\mathcal{G}} |G(y)| \beta(y^{-1}) \Delta(y)^{-1/p} \, dy \\
&= \|F\|_{L^p_m} \|H\|_{L^q_{1/m}} \|G\|_{L^1_{\beta^{\vee} \Delta^{-1/p}}}.
\end{aligned}$$

By duality, Young inequality (2.4) follows. This also implies that the convolution integral of $F * G^{\vee}$ converges absolutely almost everywhere, and consequently defines that a function in $L^p(\mathcal{G})$. □

We now turn our attention to submultiplicative weights.

**Definition 2.2.7.** A weight $w$ is called *submultiplicative* if $w(xy) \leq w(x) w(y)$ holds for all $x, y \in \mathcal{G}$.

Obviously, every submultiplicative weight is left- and right-moderate with respect to itself. Thus, we can apply the propositions 2.2.5 and 2.2.6 to submultiplicative weights.

**Proposition 2.2.8.** *Let $w$ be a submultiplicative weight. Then the space $L^1_w(\mathcal{G})$ is left- and right-invariant and the inequalities*

$$\|L_y G\|_{L^1_w} \leq w(y) \|G\|_{L^1_w} \quad \text{and} \quad \|R_y G\|_{L^1_w} \leq w(y^{-1}) \Delta(y)^{-1} \|G\|_{L^1_w}$$

*hold for all $y \in \mathcal{G}$ and $G \in L^1_w(\mathcal{G})$. Moreover, the space $L^1_w(\mathcal{G})$ is invariant under convolutions and we have the Young inequality*

$$\|G_1 * G_2\|_{L^1_w} \leq \|G_1\|_{L^1_w} \|G_2\|_{L^1_w} \tag{2.5}$$

*for all $G_1, G_2 \in L^1_w(\mathcal{G})$.*

In algebraic terms, $L^1_w(\mathcal{G})$ is a *convolution algebra*. If now $m$ is a left-moderate weight with respect to a *submultiplicative weight* $\alpha$, then Young inequality (2.2) means that $L^p_m(\mathcal{G})$ is a *left convolution module over* $L^1_\alpha(\mathcal{G})$. A short calculation shows that, for instance, $\alpha_0$ is submultiplicative. Since it is locally integrable by lemma 2.2.4, it is a submultiplicative weight function. Therefore, $L^p_m(\mathcal{G})$ is a left convolution module over $L^1_{\alpha_0}(\mathcal{G})$ if $m$ is left-moderate.

Similarly, if $m$ is right-moderate, then $L^p_m(\mathcal{G})$ is a right convolution module over $L^1_{\beta^{\vee} \Delta^{-1/p}}(\mathcal{G})$ whenever $\beta \geq \beta_0$ is a submultiplicative weight; this is certainly true for $\beta = \beta_0$.

In order to summarize the just mentioned module properties, we make the following definition.

**Definition 2.2.9.** Suppose $m$ is a moderate weight and $1 \leq p \leq \infty$. We call a submultiplicative weight $w$ a *p-control-weight of $m$* if it satisfies the inequality

$$w(y) \geq w_0(y) = \max\{\alpha_0(y), \alpha_0(y^{-1}), \beta_0(y^{-1}) \Delta(y)^{-1/p}, \beta_0(y) \Delta(y)^{-1/q}\} \tag{2.6}$$

for all $y \in \mathcal{G}$, where as usual $p^{-1} + q^{-1} = 1$ and $1/\infty = 0$.

***Remark* 2.2.10.** The terms $\Delta(y)^{-1/p}$ and $\Delta(y)^{-1/q}$ always lie between 1 and $\Delta(y)^{-1}$. Thus, if the submultiplicative weight $w$ satisfies

$$w(y) \geq \max\{\alpha_0(y), \alpha_0(y^{-1}), \beta_0(y^{-1}), \beta_0(y^{-1})\Delta(y)^{-1}, \beta_0(y), \beta_0(y)\Delta(y)^{-1})\},$$

it is a $p$-control-weight for all $p \in [1, \infty]$. We then call it a *simultaneous control-weight*.

For unimodular groups (where $\Delta \equiv 1$), any $p$-control-weight is a simultaneous control-weight. We just call it a *control-weight* then. ◁

A $p$-control-weight $w$ satisfies by definition the inequalities $w(y) \geq \alpha_0(y)$ as well as $w(y) \geq \beta_0(y^{-1})\Delta(y)^{-1/p}$. This allows us to rephrase propositions 2.2.5 and 2.2.6 in terms of $p$-control-weights.

**Theorem 2.2.11.** *Suppose $m$ is a moderate weight, $1 \leq p \leq \infty$, and $w$ is a $p$-control-weight of $m$. Then $L_m^p(\mathcal{G})$ is left- and right-invariant and the inequalities*

$$\left\|L_y F\right\|_{L_m^p} \leq w(y)\|F\|_{L_m^p} \quad \text{and} \quad \left\|R_y F\right\|_{L_m^p} \leq w(y)\|F\|_{L_m^p}$$

*hold for all $y \in \mathcal{G}$ and $F \in L_m^p(\mathcal{G})$. Moreover, $L_m^p(\mathcal{G})$ is a convolution bimodule over $L_w^1(\mathcal{G})$ in the sense that the Young inequalities*

$$\|G * F\|_{L_m^p} \leq \|G\|_{L_w^1}\|F\|_{L_m^p}, \tag{2.7}$$

$$\left\|F * G^{\vee}\right\|_{L_m^p} \leq \|F\|_{L_m^p}\|G\|_{L_w^1} \tag{2.8}$$

*hold for all $G \in L_w^1(\mathcal{G})$ and $F \in L_m^p(\mathcal{G})$.*

The other two lower bounds in (2.6) are important for compatibility with duality, as the next lemma shows.

**Lemma 2.2.12.**

(i) *If $m$ is left-moderate (right-moderate) with respect to $\alpha$ (resp. $\beta$), then $1/m$ is left-moderate (right-moderate) with respect to $\alpha^{\vee}$ (resp. $\beta^{\vee}$).*

(ii) *If $m$ is moderate and $w$ is a $p$-control-weight of $m$ for some $p \in [1, \infty]$, then $w$ also is a $q$-control-weight for $1/m$, where $p^{-1} + q^{-1} = 1$.*

*Proof.*

(i) Suppose $m$ is left-moderate with respect to $\alpha$. Then for $x, y \in \mathcal{G}$ the inequality $m(x) = m(y^{-1}yx) \leq \alpha(y^{-1})m(yx)$ holds. Thus we have

$$\frac{1}{m(yx)} \leq \alpha(y^{-1})\frac{1}{m(x)},$$

which proves the moderating inequality.

We have yet to show that $1/m$ actually is a weight function, i.e. that $m$ is locally integrable. So for any compact set $K \subset \mathcal{G}$, we need to show that the integral $\int_K 1/m(x)\,dx$ is finite. We have

$$\int_K \frac{1}{m(x)}\,dx = \int_K \frac{1}{m(xe)}\,dx \leq \int_K \alpha_0(x^{-1})\frac{1}{m(e)}\,dx,$$

where we applied the already proven inequality $m(xe)^{-1} \leq \alpha^\vee(x)m(e)^{-1}$ to $\alpha = \alpha_0$. Through the substitution $x \mapsto x^{-1}$ with $d(x^{-1}) = \Delta(x)^{-1}\,dx$ this is equal to

$$\frac{1}{m(e)} \int_{K^{-1}} \alpha_0(x)\Delta(x)^{-1}\,dx,$$

which is finite since $K^{-1}$ is compact, $\Delta^{-1}$ is continuous, and $\alpha_0$ is locally integrable (lemma 2.2.4).

In a similar way, the statement can be proven for right-moderate weights.

(ii) By (i), the weight $1/m$ is moderate with respect to the functions $\alpha_0^\vee$ and $\beta_0^\vee$. If we replace $\alpha_0$ by $\alpha_0^\vee$, $\beta_0$ by $\beta_0^\vee$ and $p$ by $q$ in the defining inequality (2.6) of control-weights, the overall maximum remains unchanged. Thus, $w$ is a $q$-control-weight of $1/m$.

$\square$

As $p$-control-weights are compatible with duality, theorem 2.2.11 is just as true if we replace the space $L_m^p(\mathcal{G})$ with $L_{1/m}^q(\mathcal{G})$. Therefore, not only $L_m^p(\mathcal{G})$ is a convolution bimodule over $L_w^1(\mathcal{G})$, but so is $L_{1/m}^q(\mathcal{G})$ as well, which means we also have the Young inequalities

$$\|G * H\|_{L_{1/m}^q} \leq \|G\|_{L_w^1} \|H\|_{L_{1/m}^q} \quad \text{and} \quad \|H * G^\vee\|_{L_{1/m}^q} \leq \|H\|_{L_{1/m}^q} \|G\|_{L_w^1}$$

for $G \in L_w^1(\mathcal{G})$ and $H \in L_{1/m}^q(\mathcal{G})$.

We can further use the duality of $L_m^p(\mathcal{G})$ and $L_{1/m}^q(\mathcal{G})$ to obtain a convolution relation between them.

**Lemma 2.2.13.** *Suppose $m$ is a moderate weight, $1 \leq p \leq \infty$, $p^{-1} + q^{-1} = 1$ and $w$ is a $p$-control-weight of $m$. Then we have the Young inequality*

$$\|F * H^\vee\|_{L_{1/w}^\infty} \leq \|F\|_{L_m^p} \|H\|_{L_{1/m}^q} \tag{2.9}$$

*for $F \in L_m^p(\mathcal{G})$ and $H \in L_{1/m}^q(\mathcal{G})$.*

*Proof.* We can write the convolution as

$$(F * H^\vee)(x) = \int_\mathcal{G} F(y) H(x^{-1}y)\,dy = \int_\mathcal{G} F(y) L_x H(y)\,dy$$

for all $x \in \mathcal{G}$. Since $w$ is a $q$-control-weight of $1/m$, the duality of $L_m^p(\mathcal{G})$ and $L_{1/m}^q(\mathcal{G})$ implies

$$|(F * H^\vee)(x)| \leq \|F\|_{L_m^p} \|L_x H\|_{L_{1/m}^q} \leq w(x) \|F\|_{L_m^p} \|H\|_{L_{1/m}^q}.$$

This proves the stated inequality. $\square$

Yet another important property of $p$-control-weights $w$ of some moderate $m$ is that they are bounded away from zero. If $x \in \mathcal{G}$ is arbitrary, the submultiplicativity of $\alpha_0$ implies

$$0 < \alpha_0(e) = \alpha_0(xx^{-1}) \leq \alpha_0(x)\alpha_0^\vee(x) \leq w(x)^2 \tag{2.10}$$

where $e \in \mathcal{G}$ is the neutral element. This yields the continuous embedding $L_w^1(\mathcal{G}) \hookrightarrow L^1(\mathcal{G})$.

*Remark 2.2.14.*

(i) It is worthwhile to ponder the technical requirements of the weight functions we use.

   If $m$ is a moderate weight, the functions $\alpha_0$ and $\beta_0$ are locally integrable. Since the map $x \mapsto x^{-1}$ is a homeomorphism on $\mathcal{G}$ and maps compact sets to compact sets, the functions $\alpha_0^\vee$ and $\beta_0^\vee$ are locally integrable as well. Thus, the function $w_0$ from (2.6) is locally integrable.

   Moreover, $w_0$ is submultiplicative itself as it is the maximum of submultiplicative functions. So $w_0$ is in fact a $p$-control-weight of $m$. This implies in particular that for any moderate weight $m$ and $1 \leq p \leq \infty$, $p$-control-weights of $m$ *do exist*.

   On the other hand, if we assume that $m$ is a measurable (not necessarily locally integrable) function which satisfies $m(xy) \leq w(x)m(y)$ for some measurable, submultiplicative (not necessarily locally integrable) $w$, then the local integrability of $w$ and $m$ already follows from [31, Thm. 2.2.22]. So the local integrability of $m$ is also necessary if we want control-weights to exist.

(ii) Our definition of $p$-control-weights (or even weights at all) is by no means standard in the literature. Many authors require their (control-) weights to satisfy additional conditions, like continuity or certain symmetries. These different conventions lead to equivalent theories, though.

   We opted to define weights and $p$-control-weights with respect to properties that are essential to coorbit theory, like local integrability (see point (i)) or inequality (2.6) [cf. 15, eq. (4.10)], and hope this underlines why each of these assumed properties is important.

◁

We will also need the following properties of weighted $L^p$-spaces. Both statements are part of [31, Lemma 2.4.6].

**Lemma 2.2.15.** *If $m$ is moderate, $1 \leq p < \infty$ and $F \in L_m^p(\mathcal{G})$, then the translation maps*

$$y \mapsto L_y F, \quad y \mapsto R_y F, \quad y \in \mathcal{G}$$

*are continuous.*

**Lemma 2.2.16.** *Let $v : \mathcal{G} \to \mathbb{C}$ be a weight that is bounded away from zero and locally bounded. Then the space of compactly supported continuous functions $C_c(\mathcal{G})$ is a dense subspace of $L_v^p(\mathcal{G})$ for $1 \leq p < \infty$.*

For further technical details about moderate and submultiplicative weights, we refer the reader to [10; 11] and section 2.2.4 of [31]. In [23], the author gives an overview of different types of weight functions on $\mathbb{R}^d$ and their usage in time-frequency analysis, including moderate and submultiplicative weights.

## 2.3 Integrals with Values in Banach Spaces

In the following theory, it is useful to have a notion of integrals with values in Hilbert and Banach spaces. Such integrals are of the form

$$\int_\mathcal{G} \varphi(y)\, dy,$$

where $\mathcal{G}$ is a locally compact group, $B$ is a Banach space and $\varphi : \mathcal{G} \to B$ is a function. There are essentially two approaches to achieve this, namely the *weak* and the *strong integral*. Weak integrals are defined in terms of the pointwise action of $\varphi(y)$ on the dual space $B'$, and therefore take values in the double dual $B''$. We will only use this integral for Hilbert spaces $\mathcal{H} = B = B''$, hence we will only describe it in this special case. The strong integral is constructed in the same way as the common Lebesgue integral.

We will only summarize the most important properties of these integrals and pass over any proofs. For more details, the reader is referred to [28, Sec. 2.3 and 3.3]. A thorough introduction to strong the integral can also be found in [26, Ch. VI].

**The Strong Integral**

The strong integral with values in a Banach space $B$ is constructed like the usual Lebesgue integral with complex values. Vice versa, the Lebesgue integral can be seen as the specialization $B = \mathbb{C}$. So we begin to define the integral for step maps and extend it to $\mu$-measurable functions afterwards.

A map $\tau : \mathcal{G} \to B$ is called *step map* if it can be written as the sum

$$\tau = \sum_{i \in I} \chi_{M_i} v_i,$$

where $M_i \subset \mathcal{G}$, $i \in I$ is a finite family of measurable sets with finite measure, and $v_i \in B$ for $i \in I$. We define the strong integral for such step maps as

$$\int_{\mathcal{G}} \tau(x)\, dx = \sum_{i \in I} \mu(M_i) v_i \in B.$$

Since $x \mapsto \|\tau(x)\|_B$ is a real-valued step map, the integral

$$\|\tau\|_1 = \int_{\mathcal{G}} \|\tau(x)\|\, dx \tag{2.11}$$

is finite. Now $\|\cdot\|_1$ defines a seminorm on the space of $B$-valued step maps.

A function $\varphi : \mathcal{G} \to B$ is called *$\mu$-measurable* if there exists a sequence of step maps that converges to $\varphi$ pointwise almost everywhere. If $\varphi$ is $\mu$-measurable, so is $x \mapsto \|\varphi(x)\|_B$ as a real-valued function.

We call a $\mu$-measurable function $\varphi$ *strongly integrable* if there exists a sequence of step maps $(\tau_n)_{n \in \mathbb{N}}$ that is Cauchy with respect to the seminorm $\|\cdot\|_1$ and that converges pointwise almost everywhere to $\varphi$. In that case, the sequence of integrals $\int_{\mathcal{G}} \tau_n(x)\, dx$ is a Cauchy sequence in $B$. We finally define the *strong integral of $\varphi$* by the limit

$$\int_{\mathcal{G}} \varphi(x)\, dx = \lim_{n \to \infty} \int_{\mathcal{G}} \tau_n(x)\, dx \in B.$$

Of course, the value of the strong integral is independent from the choice of approximating sequence $(\tau_n)_{n \in \mathbb{N}}$. Additionally, we have the following important statement known as *Bochner's theorem*.

**Proposition 2.3.1.** *[28, Prop. 2.16] Let $\varphi : \mathcal{G} \to B$ be a $\mu$-measurable function. Then $\varphi$ is strongly integrable if and only if $\|\varphi\|_B : \mathcal{G} \to \mathbb{C}$ is integrable.*

**Remark 2.3.2.** Continuous functions on $\mathcal{G}$ are always $\mu$-measurable if we assume the group to be $\sigma$-compact. According to [26, Ch. VI, M11], a function $\varphi : \mathcal{G} \to B$ is $\mu$-measurable if its domain is $\sigma$-finite, its image is separable and its preimage of any open set is measurable.

We always assume $\mathcal{G}$ to be $\sigma$-compact, so the domain of $\varphi$ is $\sigma$-finite. Since the image of locally compact spaces under continuous maps is locally compact, and locally compact spaces are separable, the second condition is also satisfied. Finally, the preimage of open sets under continuous function is open and therefore measurable, so the $\mu$-measurability of continuous functions follows. ◁

### The Weak Integral

Let $\mathcal{H}$ be a Hilbert space, where we assume the inner product $\langle\,\cdot\,,\,\cdot\,\rangle_\mathcal{H}$ to be linear in the first argument. We call $\varphi : \mathcal{G} \to \mathcal{H}$ *weakly $\mu$-measurable* if for every $f \in \mathcal{H}$ the function

$$\lambda_f : \mathcal{G} \to \mathbb{C}, \quad \lambda_f(x) = \langle \varphi(x), f \rangle_\mathcal{H}$$

is $\mu$-measurable with respect to the Haar measure $\mu$. Similarly, $\varphi$ is called *weakly integrable* if $\lambda_f$ is integrable for every $f \in \mathcal{H}$. If this is the case, we can define the antilinear map

$$\Lambda : \mathcal{H} \to \mathbb{C}, \quad f \mapsto \int_\mathcal{G} \lambda_f(x)\,dx = \int_\mathcal{G} \langle \varphi(x), f \rangle_\mathcal{H}\,dx.$$

Using the closed graph theorem it is possible to show that $\Lambda$ is a bounded antifunctional. Thus, there exists a Riesz representation $\Phi \in \mathcal{H}$ of $\Lambda$ such that

$$\langle \Phi, f \rangle_\mathcal{H} = \Lambda(f) = \int_\mathcal{G} \langle \varphi(x), f \rangle_\mathcal{H}\,dx$$

holds for all $f \in \mathcal{H}$. We call $\Phi$ the *weak integral of $\varphi$* and write $\Phi = \int_\mathcal{G} \varphi(x)\,dx$. Note that we have by definition

$$\left\langle \int_\mathcal{G} \varphi(x)\,dx, f \right\rangle_\mathcal{H} = \int_\mathcal{G} \langle \varphi(x), f \rangle_\mathcal{H}\,dx$$

for all $f \in \mathcal{H}$ and all weakly integrable functions $\varphi : \mathcal{G} \to \mathcal{H}$.

We need to following three important properties of weak and strong integrals. The first and second state that both integrals are compatible with bounded linear operators. The third says that the strong integral is indeed stronger than the weak one.

**Proposition 2.3.3.**

(i) *[28, Prop. 3.8] If $\varphi : \mathcal{G} \to \mathcal{H}_1$ is weakly integrable and $T : \mathcal{H}_1 \to \mathcal{H}_2$ is a bounded operator between Hilbert spaces, then $T \circ \varphi : \mathcal{G} \to \mathcal{H}_2$ is weakly integrable and we have*

$$\int_\mathcal{G} (T \circ \varphi)(x)\,dx = T\left( \int_\mathcal{G} \varphi(x)\,dx \right),$$

*where the integrals have to be understood in the weak sense.*

(ii) *[28, Prop. 2.18] If $\varphi : \mathcal{G} \to B_1$ is strongly integrable and $T : B_1 \to B_2$ is a bounded operator between Banach spaces, then $T \circ \varphi : \mathcal{G} \to B_2$ is strongly integrable and we have*

$$\int_\mathcal{G} (T \circ \varphi)(x)\,dx = T\left( \int_\mathcal{G} \varphi(x)\,dx \right),$$

*where the integrals have to be understood in the strong sense.*

(iii) If $\varphi : \mathcal{G} \to \mathcal{H}$ is *strongly integrable*, then it is *weakly integrable* as well and the strong and weak integral coincide.

## 2.4 Representations of Locally Compact Groups

In general, a representation is a linear group action on some vector space. When the group and vector space carry any additional structures, like a topology or a norm, representations are often ought to be compatible with these structures.

Since we are going to work with locally compact groups that act on Banach spaces and Hilbert spaces, we consider representations that are compatible with norms and inner products. This section introduces some elementary concepts of representations that only depend on a norm but not on an inner product.

**Definition 2.4.1.** Let $\mathcal{G}$ be a locally compact group and $V$ be a normed space. A *representation of $\mathcal{G}$ on $V$* is a group homomorphism $\pi : \mathcal{G} \to \mathrm{GL}(V)$ into the group of bounded and boundedly invertible operators on $V$ such that for every $v \in V$ the map $\mathcal{G} \ni x \mapsto \pi(x)v \in V$ is continuous.

A representation $\pi$ defines a linear group action of $\mathcal{G}$ on $V$ by $(x, v) \mapsto \pi(x)v$. The condition of continuity means that the mapping $\pi : \mathcal{G} \to \mathrm{GL}(V)$ is continuous with respect to the strong operator topology on $\mathrm{GL}(V)$. Hence, this property is called the *strong continuity of $\pi$*.

Examples of linear group actions are the translation operators $L_x$ and $R_x$ for any locally compact group $\mathcal{G}$. Since we have

$$L_{xy}F(z) = F((xy)^{-1}z) = F(y^{-1}x^{-1}z) = (L_yF)(x^{-1}z) = L_xL_yF(z)$$

and

$$R_{xy}F(z) = F(zxy) = (R_yF)(zx) = R_xR_yF(z)$$

for all functions $F : \mathcal{G} \to \mathbb{C}$ and all $x, y, z \in \mathcal{G}$, these translations are compatible with the structure of $\mathcal{G}$. To get a representation, we need $L$ and $R$ to act boundedly and strongly continuous on a function space on $\mathcal{G}$. As we have seen in section 2.2, this is, for instance, the case for the weighted $L^p$-spaces $L_m^p(\mathcal{G})$ with moderate $m$ and $p \in [1, \infty)$. Thus,

$$L, R : \mathcal{G} \to \mathrm{GL}(L_m^p(\mathcal{G})), \quad L(x) = L_x, \quad R(x) = R_x$$

are representations of $\mathcal{G}$ on $L_m^p(\mathcal{G})$. For $p = 2$ and $m = 1$, these are called the *left regular* respectively the *right regular representation of $\mathcal{G}$*.

**Example 2.4.2.** We can generalize the idea behind regular representations. Let $H \subset \mathrm{GL}(\mathbb{R}^n)$ be a closed matrix group. We then define the *semidirect product* $\mathbb{R}^n \rtimes H$ as the topological space $\mathbb{R}^n \times H$ which is equipped with the group structure

$$(b_1, A_1)(b_2, A_2) = (b_1 + A_1 b_2, A_1 A_2).$$

The semidirect product $\mathbb{R} \rtimes H$ is a locally compact group.

Now $\mathbb{R}^n \rtimes H$ acts on $\mathbb{R}^n$ by $(b, A)x = Ax + b$ for $(b, A) \in \mathbb{R} \times H$ and $x \in \mathbb{R}^n$. This induces the representation $\pi : \mathbb{R}^n \rtimes H \to \mathrm{GL}(L^p(\mathbb{R}^n))$,

$$\pi(b, A)f(x) = |\det(A)|^{-1/p} f((b, A)^{-1}x) = |\det(A)|^{-1/p} f(A^{-1}(x - b))$$

for $p \in [1, \infty)$. For $p = 2$, it is called the *quasiregular representation*.

In example 2.1.3, we have encountered the affine group $\mathcal{A\!f\!f}$. Using the semidirect product, it can be written as $\mathcal{A\!f\!f} = \mathbb{R} \rtimes \mathbb{R}^*$. Its quasiregular representation

$$\pi : \mathcal{A\!f\!f} \to \mathrm{GL}(L^2(\mathbb{R})), \quad \pi(b,a)f(x) = \frac{1}{\sqrt{|a|}} f\left(\frac{x-b}{a}\right)$$

is called the *wavelet representation*. ◁

One of the central concepts in representation theory is *irreducibility*. Similar to other algebraic structures, it means that a representation has no proper non-trivial 'subrepresentations'. The following definition makes this concept rigid.

**Definition 2.4.3.** Let $\pi$ be a representation of $\mathcal{G}$ on $V$.

a) A vector $v \in V$ is called *cyclic*, if the space

$$\mathcal{E}_v = \mathrm{span}\{\pi(x)v \mid x \in \mathcal{G}\}$$

is dense in $V$.

b) A subspace $U \subset V$ is called $\pi$-*invariant*, if $\pi(x)U \subset U$ for every $x \in \mathcal{G}$.

c) The representation $\pi$ is called *irreducible*, if one the following equivalent conditions hold:

   i) Every $v \in V \backslash \{0\}$ is cyclic.
   ii) All non-trivial $\pi$-invariant subspaces of $V$ are dense in $V$.

   The trivial representation $V = \{0\}$ is explicitly defined to be not irreducible.

It is easy to see that the conditions i) and ii) are indeed equivalent. If we assume the first condition to be true and take any non-trivial $\pi$-invariant subspace $U \subset V$, then there exists some $v \in U \backslash \{0\}$. Now $\mathcal{E}_v \subset U$ is dense in $V$ since $v$ is cyclic, thus $U$ is dense itself. Vice versa, if every non-trivial $\pi$-invariant subspace is dense in $V$, this is especially true for all the $\pi$-invariant spaces $\mathcal{E}_v$, $v \in V \backslash \{0\}$.

## 2.5 Unitary Representations and the Voice Transform

In this chapter, we will consider representations that are compatible with inner products. Associated to those representations is the abstract *voice transform*, which generalizes, for instance, the wavelet transform and the short-time Fourier transform.

We use the convention that inner products are linear in the first argument.

**Definition 2.5.1.** Let $\mathcal{G}$ be a locally compact group and $\mathcal{H}$ be a Hilbert space. A *unitary representation of $\mathcal{G}$ on $\mathcal{H}$* is a representation $\pi : \mathcal{G} \to \mathcal{U}(\mathcal{H}) \subset \mathrm{GL}(\mathcal{H})$, that is, a strongly continuous group homomorphism from $\mathcal{G}$ into the group of unitary operators on $\mathcal{H}$.

A unitary representation therefore satisfies

$$\langle \pi(x)f, \pi(x)g \rangle_{\mathcal{H}} = \langle f, g \rangle_{\mathcal{H}}$$

for all $f, g \in \mathcal{H}$ and $x \in \mathcal{G}$, or similarly

$$\langle f, \pi(x)g \rangle_{\mathcal{H}} = \langle \pi(x^{-1})f, g \rangle_{\mathcal{H}}.$$

The left regular representation $L$ of $\mathcal{G}$ on $L^2(\mathcal{G})$ is a prominent example of a unitary representation, as the left invariance of the Haar measure implies

$$\langle F, G \rangle_{L^2(\mathcal{G})} = \int_{\mathcal{G}} F(x)\overline{G(x)}\, dx = \int_{\mathcal{G}} F(y^{-1}x)\overline{G(y^{-1}x)}\, dx = \langle L_y F, L_y G \rangle_{L^2(\mathcal{G})}$$

for all $F, G \in L^2(\mathcal{G})$ and $y \in \mathcal{G}$.

The right regular representation $R$ does *not* act unitarily on $L^2(\mathcal{G})$. Actually, we have

$$\langle F, G \rangle_{\mathcal{H}} = \int_{\mathcal{G}} F(x)\overline{G(x)}\, dx = \int_{\mathcal{G}} F(xy)\overline{G(xy)}\Delta(y)\, dx = \Delta(y)\langle R_y F, R_y G \rangle_{\mathcal{H}},$$

where we used the substitution rule $d(xy) = \Delta(y)\, dx$. Therefore, $\Delta^{1/2}R$ is a unitary action and indeed a unitary representation of $\mathcal{G}$ on $L^2(\mathcal{G})$. For this reason, $\Delta^{1/2}R$ is sometimes defined to be the right regular representation.

We now take a look at the voice transform.

**Definition 2.5.2.** Let $\pi$ be a unitary representation of $\mathcal{G}$ on $\mathcal{H}$. We define the *voice transform* (also called *generalized wavelet transform*) for $f, g \in \mathcal{H}$ by

$$\mathcal{V}_g f : \mathcal{G} \to \mathbb{C}, \quad \mathcal{V}_g f(x) = \langle f, \pi(x)g \rangle_{\mathcal{H}}.$$

The voice transform is obviously linear in $f$ and antilinear in $g$.

For fixed $f, g \in \mathcal{H}$, the function $\mathcal{V}_g f$ is continuous and bounded. The continuity is an immediate consequence of the strong continuity of $\pi$, while the boundedness follows from the Cauchy-Schwarz inequality

$$|\mathcal{V}_g f(x)| = \langle f, \pi(x)g \rangle_{\mathcal{H}} \leq \|f\|_{\mathcal{H}} \|\pi(x)g\|_{\mathcal{H}} = \|f\|_{\mathcal{H}} \|g\|_{\mathcal{H}}$$

for all $x \in \mathcal{G}$. Thus, for fixed $g$, the voice transform defines an operator $\mathcal{V}_g : \mathcal{H} \to L^\infty(\mathcal{G})$ which is bounded by $\|g\|_{\mathcal{H}}$.

*Example 2.5.3.*

i. The integers $\mathbb{Z}$ become a locally compact abelian group if equipped with the discrete topology. Their Haar measure is the counting measure, so the integral over $\mathbb{Z}$ is just the series over $\mathbb{Z}$.

   Now $\mathbb{Z}$ acts on $L^2(0,1)$ by modulation, that is

   $$\pi(k)f(t) = e^{2\pi i k t} f(t), \quad t \in [0,1)$$

   for $k \in \mathbb{Z}$ and $f \in L^2(0,1)$. This action is compatible with the group structure of $\mathbb{Z}$. In addition, $k \mapsto \pi(k)f$ is continuous for all $f$ since $\mathbb{Z}$ is discrete. Thus, $\pi$ is a representation of $\mathbb{Z}$ on $L^2(0,1)$.

   The calculation

   $$\|\pi(k)f\|_{L^2} = \left( \int_0^1 |e^{2\pi i k t} f(t)|^2\, dt \right)^{1/2} = \left( \int_0^1 |f(t)|^2\, dt \right)^{1/2} = \|f\|_{L^2}$$

shows that $\pi$ is a unitary representation. Its voice transform

$$\mathcal{V}_g f(k) = \langle f, \pi(k)g \rangle_{L^2(0,1)} = \int_0^1 f(t)\overline{g(t)} e^{-2\pi i k t} \, dt$$

is a generalization of the periodic Fourier transform. By setting $g = 1$, we get the usual Fourier coefficients of $f$.

ii. We have already encountered the wavelet representation

$$\pi : \mathcal{A}\!f\!f \to \mathcal{U}(L^2(\mathbb{R})), \quad \pi(b,a)f(t) = \frac{1}{\sqrt{|a|}} f\left(\frac{b-t}{a}\right)$$

of the affine group $\mathcal{A}\!f\!f = \mathbb{R} \rtimes \mathbb{R}^*$ on $L^2(\mathbb{R})$ in example 2.4.2. It is unitary as can be seen by a short computation. Its voice transform is given by

$$\mathcal{V}_g f(b,a) = \langle f, \pi(b,a)g \rangle_{L^2(\mathbb{R})} = \frac{1}{\sqrt{|a|}} \int_{\mathbb{R}} f(t)\overline{g\left(\frac{t-b}{a}\right)} \, dt$$

for $f, g \in L^2(\mathbb{R})$. This is exactly the wavelet transform $W_g f$. We will take a closer look at this example in section 4.1.

iii. Similar to the preceding example, we would like to identify the short-time Fourier transform as the voice transform of some unitary representation. The short-time Fourier transform is defined for $f, g \in L^2(\mathbb{R}^d)$ and $x, \omega \in \mathbb{R}^d$ as

$$V_g f(x, \omega) = \langle f, M_\omega T_x g \rangle_{L^2(\mathbb{R}^d)} = \int_{\mathbb{R}^d} f(t)\overline{g(t-x)} e^{-2\pi i t \omega} \, dt,$$

where we used the translation and modulation operators $T_x, M_\omega : L^2(\mathbb{R}^d) \to L^2(\mathbb{R}^d)$,

$$T_x g(t) = f(t-x), \quad M_\omega g(t) = e^{2\pi i t \omega} f(t), \quad t \in \mathbb{R}^d.$$

By $t\omega$ we mean the dot product $t\omega = \sum_{j=1}^d t_j \omega_j$.

If we want to associate the short-time Fourier transform with some group structure, we need to find a group operation $(x, \omega)(x', \omega') = (x_0, \omega_0)$ on $\mathbb{R}^d \times \mathbb{R}^d$ that is compatible with the concatenation of the operators $(M_\omega T_x)(M_{\omega'} T_{x'})$. Calculating this concatenation for some $g \in L^2(\mathbb{R}^d)$ shows that

$$\begin{aligned}
(M_\omega T_x)(M_{\omega'} T_{x'})g(t) &= e^{2\pi i t \omega}(M_{\omega'} T_{x'})g(t-x) \\
&= e^{2\pi i (t-x)\omega'} e^{2\pi i t \omega} g(t - x - x') \\
&= e^{-2\pi i \omega' x} e^{2\pi i (\omega + \omega')t} g(t - (x + x')) \\
&= e^{-2\pi i \omega' x} M_{\omega + \omega'} T_{x + x'} g(t).
\end{aligned}$$

This concatenation cannot be written as $M_{\omega_0} T_{x_0}$ for any tuple $(x_0, \omega_0) \in \mathbb{R}^d \times \mathbb{R}^d$. Thus, there exists no group-structure that is compatible with $M_\omega T_x$, which means that the short-time Fourier transform *is not the voice transform of any representation* (in other words: the set of time-frequency shifts is not a subgroup of $\mathcal{U}(L^2(\mathbb{R}^d))$).

We will solve this problem in section 4.2 by extending the space $\mathbb{R}^d \times \mathbb{R}^d$ by a third parameter from $S^1$ that takes care of the interfering factor $e^{-2\pi i \omega' x}$. The group constructed that way is the *reduced Heisenberg group* $\mathbb{H}_r^d$.

◁

The voice transform fulfils some equations which we will use quite frequently. For $f, g \in \mathcal{H}$ and $x \in \mathcal{G}$ we have

$$\mathcal{V}_g f(x) = \langle f, \pi(x)g \rangle_\mathcal{H} = \overline{\langle g, \pi(x^{-1})f \rangle_\mathcal{H}} = \overline{\mathcal{V}_f g(x^{-1})},$$

so in short $\mathcal{V}_g f = \overline{\mathcal{V}_f g^\vee}$. In particular, $\mathcal{V}_g g$ is invariant under the involution $\overline{\mathcal{V}_g g^\vee}$. Another calculation shows

$$\mathcal{V}_g(\pi(y)f)(x) = \langle \pi(y)f, \pi(x)g \rangle_\mathcal{H} = \langle f, \pi(y^{-1}x)g \rangle_\mathcal{H} = \mathcal{V}_g f(y^{-1}x) = L_y \mathcal{V}_g f(x)$$

for $y \in \mathcal{G}$, thus $\mathcal{V}_g(\pi(y)f) = L_y \mathcal{V}_g f$. Therefore, the voice transform associates the representation $\pi$ with the left translation $L$.

We now turn our attention to the orthogonality relations of the voice transform. To formulate such relations, the image of the transform $\mathcal{V}_g(\mathcal{H}) \subset L^\infty(\mathcal{G})$ has be contained in a Hilbert space, namely $L^2(\mathcal{G})$. Thus, we need a criteria that tells us when $\mathcal{V}_g : \mathcal{H} \to L^2(\mathcal{G})$ is a well-defined operator.

If we assume the representation to be irreducible, this leads us to the condition $\mathcal{V}_g g \in L^2(\mathcal{G})$.

**Proposition 2.5.4.** *[cf. 9, Thm. 2] Let $\pi$ be a unitary and irreducible representation of $\mathcal{G}$ on $\mathcal{H}$. Suppose $g \in \mathcal{H}$ satisfies $\mathcal{V}_g g \in L^2(\mathcal{G})$, that is*

$$\int_\mathcal{G} |\mathcal{V}_g g(x)|^2 \, dx = \int_\mathcal{G} |\langle g, \pi(x)g \rangle_\mathcal{H}|^2 \, dx < \infty.$$

*Then $\mathcal{V}_g f \in L^2(\mathcal{G})$ for all $f \in \mathcal{H}$ and the operator $\mathcal{V}_g : \mathcal{H} \to L^2(\mathcal{G})$ is the multiple of an isometry. Thus, there is a constant $C_g \in [0, \infty)$ such that*

$$\langle \mathcal{V}_g f_1, \mathcal{V}_g f_2 \rangle_{L^2(\mathcal{G})} = C_g^2 \langle f_1, f_2 \rangle_\mathcal{H}$$

*for all $f_1, f_2 \in \mathcal{H}$.*

The fact that $\mathcal{V}_g$ is the multiple of an isometry is very important to us. To make use of the preceding proposition, we will always use representations that are irreducible and for which some $g$ with $\mathcal{V}_g g \in L^2(\mathcal{G})$ exists. We write down these conditions in the next definition. To emphasize that a vector has such a special property, we will note it from now on as $\psi$, $\phi$ and so on.

**Definition 2.5.5.** Let $\pi$ be a unitary and irreducible representation of $\mathcal{G}$ on $\mathcal{H}$.

  a) We call a vector $\psi \in \mathcal{H}$ *square-integrable* if $\mathcal{V}_\psi \psi \in L^2(\mathcal{G})$. We denote the set of all square-integrable vectors by $\mathcal{Q}$. The representation $\pi$ is called square-integrable if at least one non-trivial square-integrable vector exists, i.e. if $\mathcal{Q} \neq \{0\}$.

  b) Let $\psi$ be square-integrable. The constant $C_\psi$ from proposition 2.5.4 is called the *admissibility constant of $\psi$*. The vector $\psi$ is called *admissible* if $C_\psi = 1$, i.e. if $\mathcal{V}_\psi : \mathcal{H} \to L^2(\mathcal{G})$ is an isometry.

**Remark 2.5.6.** We will often assume that the vector $\psi \in \mathcal{Q} \backslash \{0\}$ is admissible so we do not need to keep track of the constant $C_\psi$. But, since any $\psi_0 \in \mathcal{Q} \backslash \{0\}$ can be normalized by $\psi = \psi_0/C_{\psi_0}$ to obtain an admissible vector, most statements can be generalized to non-admissible vectors as well.

◁

There is an important extension of proposition 2.5.4 regarding the isometry relation. This extension establishes an orthogonality relation between the voice transforms of two *different* square-integrable vectors.

**Theorem 2.5.7** (Duflo-Moore). *[9, Thm. 3] Let $\pi$ be a square-integrable representation of $\mathcal{G}$ on $\mathcal{H}$. Then there exists a unique densely defined operator $D : \mathcal{D}(D) \subset \mathcal{H} \to \mathcal{H}$ with the following properties:*

  i. *The operator $D$ is self-adjoint, positive, injective and has a densely defined inverse.*

  ii. *The domain of $D$ consists exactly of the square-integrable vectors, i.e. $\mathcal{D}(D) = \mathcal{Q}$.*

  iii. *Given $\psi_1, \psi_2 \in \mathcal{Q}$ and $f_1, f_2 \in \mathcal{H}$, the orthogonality relation*

$$\langle \mathcal{V}_{\psi_1} f_1, \mathcal{V}_{\psi_2} f_2 \rangle_{L^2(\mathcal{G})} = \overline{\langle D\psi_1, D\psi_2 \rangle_{\mathcal{H}}} \langle f_1, f_2 \rangle_{\mathcal{H}} \tag{2.12}$$

  *holds.*

  iv. *For $x \in \mathcal{G}$ we have*
$$\pi(x)D = \sqrt{\Delta(x)} D\pi(x).$$

*We call this operator $D$ the* Duflo-Moore operator.

If $\psi$ is square-integrable, we now have for $f \in \mathcal{H}$

$$C_\psi^2 \|f\|_{\mathcal{H}}^2 = \|\mathcal{V}_\psi f\|_{L^2}^2 = \|D\psi\|_{\mathcal{H}}^2 \|f\|_{\mathcal{H}}^2.$$

Thus, we can write the admissibility constant of $\psi$ as $C_\psi = \|D\psi\|_{\mathcal{H}}$.

Another immediate consequence of item ii. is that the set of square-integrable vectors $\mathcal{Q}$ is dense in $\mathcal{H}$, as $D$ is densely defined.

If $\mathcal{G}$ is unimodular, the Duflo-Moore operator is a (positive) multiple of the identity $D = c\,\mathrm{id}_{\mathcal{H}}$ on $\mathcal{H}$ (in particular, $D$ is defined everywhere). In that case, every vector in $\mathcal{H}$ is square-integrable and we have the simplified orthogonality relation

$$\langle \mathcal{V}_{g_1} f_1, \mathcal{V}_{g_2} f_2 \rangle_{L^2(\mathcal{G})} = c^2 \overline{\langle g_1, g_2 \rangle_{\mathcal{H}}} \langle f_1, f_2 \rangle$$

for all $g_1, g_2, f_1, f_2 \in \mathcal{H}$.

At the end of this section, we state two facts regarding the injectivity of the voice transform.

**Lemma 2.5.8.** *[1, Lemma 2.16] Let $\pi$ be a unitary representation of $\mathcal{G}$ on $\mathcal{H}$. Then for $g \in \mathcal{H}$ the following are equivalent.*

  i) *$g$ is cyclic.*

  ii) *$\mathcal{V}_g$ is injective.*

*Proof.* We denote $\mathcal{E}_g = \mathrm{span}\{\pi(x)g \mid x \in \mathcal{G}\}$. The kernel of $\mathcal{V}_g$ is exactly the orthogonal complement of $\mathcal{E}_g$ since we have

$$\begin{aligned}
\mathcal{E}_g^\perp &= \{f \in \mathcal{H} \mid \langle f, h \rangle_{\mathcal{H}} = 0 \text{ for all } h \in \mathcal{E}_g\} \\
&= \{f \in \mathcal{H} \mid \langle f, \pi(x)g \rangle_{\mathcal{H}} = 0 \text{ for all } x \in \mathcal{G}\} \\
&= \ker(\mathcal{V}_g).
\end{aligned}$$

Thus we have the equivalences

$$g \text{ is cyclic} \iff \overline{\mathcal{E}}_g = \mathcal{H} \iff \mathcal{E}_g^\perp = \{0\} \iff \mathcal{V}_g \text{ is injective}.$$

$\square$

By definition, the representation $\pi$ is irreducible whenever all $g \in \mathcal{H}\backslash\{0\}$ are cyclic. Therefore, we can characterize the irreducibility of $\pi$ by means of its voice transform.

**Corollary 2.5.9.** *A (non-trivial) unitary representation $\pi : \mathcal{G} \to \mathcal{U}(\mathcal{H})$ is irreducible if and only if $\mathcal{V}_g : \mathcal{H} \to L^\infty(\mathcal{G})$ is injective for all $g \in \mathcal{H}\backslash\{0\}$.*

## 2.6 Reconstruction and Reproducing Formulas

In this section, we want to use the orthogonality relations from proposition 2.5.4 and the Duflo-Moore theorem 2.5.7 to derive two kinds of formulas, namely *reconstruction* and *reproducing formulas*. For that, we need a square-integrable representation $\pi$ of a locally compact group $\mathcal{G}$ on a Hilbert space $\mathcal{H}$.

Both kinds of formulas are closely related to the fact that $\mathcal{V}_\psi$ is an isometry for admissible $\psi$. Let $A : \mathcal{H}_1 \to \mathcal{H}_2$ be an isometry between Hilbert spaces and $A^*$ its adjoint. Then $A^*A : \mathcal{H}_1 \to \mathcal{H}_1$ is the identity, while $AA^* : \mathcal{H}_2 \to \mathcal{H}_2$ is the orthogonal projection from $\mathcal{H}_2$ onto the image of $A$. For admissible $\psi$, the reconstruction formula is a weak integral formula that corresponds to the concatenation $\mathcal{V}_\psi^* \mathcal{V}_\psi : \mathcal{H} \to \mathcal{H}$, and therefore expresses the fact that every $f \in \mathcal{H}$ can be reconstructed from its voice transform via weak integration. The reproducing formula is a convolution identity that corresponds to $\mathcal{V}_\psi \mathcal{V}_\psi^* : L^2(\mathcal{G}) \to L^2(\mathcal{G})$, thus it describes the orthogonal projection onto the image of $\mathcal{V}_\psi$. Like the Duflo-Moore theorem 2.5.7, generalized versions of both formulas are also able to relate the voice transforms of distinct square-integrable vectors.

To derive these formulas, we use similar methods as in [1, Sec. 2.5] or [31, Rem. 2.4.4].

We begin with the reconstruction formula. We want to write the operator $\mathcal{V}_\psi^* : L^2(\mathcal{G}) \to \mathcal{H}$ for square-integrable $\psi$ as a weak integral. For $F \in L^2(\mathcal{G})$ and $g \in \mathcal{H}$ we have by definition

$$\begin{aligned}
\langle \mathcal{V}_\psi^* F, g \rangle_\mathcal{H} &= \langle F, \mathcal{V}_\psi g \rangle_{L^2(\mathcal{G})} \\
&= \int_\mathcal{G} F(y) \overline{\mathcal{V}_\psi g(y)} \, dy \\
&= \int_\mathcal{G} F(y) \overline{\langle g, \pi(y)\psi \rangle_\mathcal{H}} \, dy \\
&= \int_\mathcal{G} \langle F(y)\, \pi(y)\psi, g \rangle_\mathcal{H} \, dy. \quad (2.13)
\end{aligned}$$

Now the function $y \mapsto \langle F(y)\pi(y)\psi, g \rangle_\mathcal{H} = F(y)\overline{\mathcal{V}_\psi g(y)}$ is $\mu$-measurable for all $g \in \mathcal{H}$, since $F$ is $\mu$-measurable and $\mathcal{V}_\psi g$ continuous. This means that the map $y \mapsto F(y)\, \pi(y)\psi \in \mathcal{H}$ is weakly $\mu$-measurable. Moreover, we have for all $g \in \mathcal{H}$ the estimate

$$\int_\mathcal{G} |\langle F(y)\, \pi(y)\psi, g \rangle_\mathcal{H}| \, dy = \int_\mathcal{G} |F(y)| |\mathcal{V}_\psi g(y)| \, dy \leq \|F\|_{L^2} \|\mathcal{V}_\psi g\|_{L^2} < \infty,$$

which shows that the map is weakly integrable as well. We can therefore use (2.13) to write $\mathcal{V}_\psi^* F$ as a weak integral.

**Lemma 2.6.1.** *If $\psi$ is square-integrable, the adjoint of the voice transform $\mathcal{V}_\psi : \mathcal{H} \to L^2(\mathcal{G})$ is given by the weak integral*

$$\mathcal{V}_\psi^* : L^2(\mathcal{G}) \to \mathcal{H}, \quad \mathcal{V}_\psi^*(F) = \int_\mathcal{G} F(y)\,\pi(y)\psi\,dy.$$

Using the Duflo-Moore theorem 2.5.7 we get

$$\langle \mathcal{V}_\psi^* \mathcal{V}_\phi f, g \rangle_\mathcal{H} = \langle \mathcal{V}_\phi f, \mathcal{V}_\psi g \rangle_{L^2(\mathcal{G})} = \langle \overline{\langle D\phi, D\psi \rangle_\mathcal{H}} f, g \rangle_\mathcal{H}$$

for square-integrable $\psi, \phi \in \mathcal{Q}$ and arbitrary $f, g \in \mathcal{H}$. Since this equation is true for every $g \in \mathcal{H}$, we have

$$\mathcal{V}_\psi^* \mathcal{V}_\phi f = \overline{\langle D\phi, D\psi \rangle_\mathcal{H}} f.$$

Combining this observation with lemma 2.6.1, we obtain the (generalized) reconstruction formula.

**Proposition 2.6.2** (Weak reconstruction formula). *Let $\psi, \phi$ be square-integrable vectors. Then we have for all $f \in \mathcal{H}$*

$$\overline{\langle D\phi, D\psi \rangle_\mathcal{H}} f = \int_\mathcal{G} \mathcal{V}_\phi f(y)\,\pi(y)\psi\,dy, \tag{2.14}$$

*where the integral has to be understood in the weak sense.*

*In particular, we have for admissible $\psi$*

$$f = \int_\mathcal{G} \mathcal{V}_\psi f(y)\,\pi(y)\psi\,dy \tag{2.15}$$

*for all $f \in \mathcal{H}$.*

We now derive the reproducing formula. Let $\psi, \phi$ be square-integrable and $f, g \in \mathcal{H}$ arbitrary. Then an application of Duflo-Moore's theorem 2.5.7 yields for $x \in \mathcal{G}$

$$\overline{\langle D\psi, D\phi \rangle_\mathcal{H}} \mathcal{V}_g f(x) = \overline{\langle D\psi, D\phi \rangle_\mathcal{H}} \langle f, \pi(x)g \rangle_\mathcal{H} = \langle \mathcal{V}_\psi f, \mathcal{V}_\phi(\pi(x)g) \rangle_{L^2(\mathcal{G})}.$$

Thus we have

$$\begin{aligned} \mathcal{V}_g f(x) \overline{\langle D\psi, D\phi \rangle_\mathcal{H}} &= \int_\mathcal{G} \mathcal{V}_\psi f(y) \overline{\mathcal{V}_\phi(\pi(x)g)(y)}\,dy \\ &= \int_\mathcal{G} \mathcal{V}_\psi f(y) \overline{\mathcal{V}_\phi g(x^{-1}y)}\,dy \\ &= \int_\mathcal{G} \mathcal{V}_\psi f(y) \mathcal{V}_\phi g^\triangledown(y^{-1}x)\,dy \\ &= (\mathcal{V}_\psi f * \mathcal{V}_\phi g^\triangledown)(x). \end{aligned}$$

All the above integrals are absolutely convergent due to the Cauchy-Schwarz inequality.

**Proposition 2.6.3.** *Let $\psi, \phi$ be square-integrable and $f, g \in \mathcal{H}$. Then we have*

$$\overline{\langle D\psi, D\phi \rangle_\mathcal{H}} \mathcal{V}_g f = \mathcal{V}_\psi f * \mathcal{V}_\phi g^\triangledown = \mathcal{V}_\psi f * \mathcal{V}_g \phi,$$

*where the convolutions are pointwise absolutely convergent.*

We explicitly state the three most important special cases of this general formula.

**Corollary 2.6.4.** *Let $\psi, \phi$ be square-integrable and $f, g \in \mathcal{H}$. Then we have*

$$\overline{\langle D\psi, D\phi \rangle}_\mathcal{H} \mathcal{V}_\phi \psi = \mathcal{V}_\psi \psi * \mathcal{V}_\phi \phi \tag{2.16}$$

*as well as*

$$C_\psi^2 \mathcal{V}_g f = \mathcal{V}_\psi f * \mathcal{V}_g \psi. \tag{2.17}$$

*If $\psi$ is admissible, we particularly have the* reproducing formula

$$\mathcal{V}_\psi f = \mathcal{V}_\psi f * \mathcal{V}_\psi \psi. \tag{2.18}$$

Equation (2.17) can be seen as a change of basis. The voice transforms $\mathcal{V}_\psi$ and $\mathcal{V}_g$ associate to a vector $f \in \mathcal{H}$ a continuous family of 'coefficients' $\mathcal{V}_\psi f$ and $\mathcal{V}_g f$. We can represent $f$ through these coefficients by applying the reconstruction formula 2.6.2. The function $\mathcal{V}_g \psi$ acts like a change of basis on the coefficients. Through this, we are able to show that certain properties of $f$ and its voice transform are independent from the used vector $\psi$ and $g$. Note that this formula is true even if $g$ is not square-integrable.

The above corollary also implies that for admissible $\psi$, the voice transform $\mathcal{V}_\psi \psi$ is convolution idempotent, that is

$$\mathcal{V}_\psi \psi * \mathcal{V}_\psi \psi = \mathcal{V}_\psi \psi.$$

The reproducing formula (2.18) can be understood in the sense of *reproducing kernel Hilbert spaces*. A reproducing kernel Hilbert space is a Hilbert space $H$ of functions $F : X \to \mathbb{C}$ on some set $X$, such that the evaluation functionals

$$E_x : H \to \mathbb{C}, \quad F \mapsto F(x)$$

are bounded for all $x \in X$. In that case, every functional $E_x$ has a Riesz representation $k_x \in H$, hence we have $F(x) = \langle F, k_x \rangle_H$ for all $F \in H$. The function

$$K(x, y) = \langle k_y, k_x \rangle_H = k_y(x) = \overline{k_x(y)}$$

is then called the *reproducing kernel of $H$*.

Given a square-integrable representation $\pi$ and an admissible vector $\psi$, the Hilbert space $\mathcal{M}^2 = \mathcal{V}_\psi(\mathcal{H}) \subset L^2(\mathcal{G})$ consists of continuous functions on $\mathcal{G}$. The evaluation functionals are bounded since

$$|\mathcal{V}_\psi f(x)| = |\langle f, \pi(x)\psi \rangle_\mathcal{H}| \leq \|f\|_\mathcal{H} \|\psi\|_\mathcal{H} = \|\mathcal{V}_\psi f\|_{L^2} \|\psi\|_\mathcal{H}$$

for all $x \in \mathcal{G}$ and $f \in \mathcal{H}$. Thus, $\mathcal{M}^2$ is a reproducing kernel Hilbert space. Now the reproducing formula $F = F * \mathcal{V}_\psi \psi$ holds for all $F \in \mathcal{M}^2$, so it follows for $x \in \mathcal{G}$

$$F(x) = (F * \mathcal{V}_\psi \psi)(x) = \int_\mathcal{G} F(y) \mathcal{V}_\psi \psi(y^{-1}x) \, dy$$
$$= \int_\mathcal{G} F(y) \overline{\mathcal{V}_\psi \psi(x^{-1}y)} \, dy$$
$$= \langle F, L_x \mathcal{V}_\psi \psi \rangle_{L^2(\mathcal{G})}.$$

In this case we therefore have $k_x = L_x \mathcal{V}_\psi \psi$ and $K(x, y) = L_y \mathcal{V}_\psi \psi(x) = \mathcal{V}_\psi \psi(y^{-1}x)$.

As already mentioned, the reproducing formula is associated to the concatenation $\mathcal{V}_\psi \mathcal{V}_\psi^*$. Let $\psi$ be admissible, $F \in L^2(\mathcal{G})$ and $x \in \mathcal{G}$. Then we can use the weak integral formula from lemma 2.6.1 to obtain

$$\begin{aligned} \mathcal{V}_\psi \mathcal{V}_\psi^* F(x) &= \langle \mathcal{V}_\psi^* F, \pi(x)\psi \rangle_\mathcal{H} \\ &= \int_\mathcal{G} F(y) \langle \pi(y)\psi, \pi(x)\psi \rangle_\mathcal{H} \, dy \\ &= \int_\mathcal{G} F(y) \langle \psi, \pi(y^{-1}x)\psi \rangle_\mathcal{H} \, dy \\ &= (F * \mathcal{V}_\psi \psi)(x). \end{aligned}$$

Therefore, right convolution with the kernel $\mathcal{V}_\psi \psi$ acts on $L^2(\mathcal{G})$ the same as applying $\mathcal{V}_\psi \mathcal{V}_\psi^*$. But since $\mathcal{V}_\psi$ is an isometry, this is precisely the orthogonal projection onto the image space $\mathcal{M}^2 = \mathcal{V}_\psi(\mathcal{H})$.

We conclude that the image space $\mathcal{M}^2$ contains exactly those functions $F \in L^2(\mathcal{G})$ which satisfy the reproducing formula $F = F * \mathcal{V}_\psi \psi$. The voice transform is an isometric isomorphism between $\mathcal{H}$ and $\mathcal{M}^2$. This fact is one of the key features of square-integrable representations.

**Theorem 2.6.5** (Correspondence principle, unitary version). *Let $\pi$ be a square-integrable representation of $\mathcal{G}$ on $\mathcal{H}$. Suppose $\psi \in \mathcal{H}$ is admissible. Then*

$$P : L^2(\mathcal{G}) \to \mathcal{V}_\psi(\mathcal{H}) \subset L^2(\mathcal{G}), \quad F \mapsto F * \mathcal{V}_\psi \psi$$

*is the orthogonal projection onto the image of $\mathcal{V}_\psi$. In particular, we have the equivalence for $F \in L^2(\mathcal{G})$*

$$F = \mathcal{V}_\psi f \text{ for some } f \in \mathcal{H} \quad \Longleftrightarrow \quad F = F * \mathcal{V}_\psi \psi.$$

CHAPTER 3

Construction of Coorbit Spaces

In this chapter, we will construct the coorbit spaces $\mathcal{C}o_m^p$ and note some of their essential properties. Starting with a square-integrable representation, these spaces are defined to contain all vectors for which the voice transform satisfies a decay condition described by a weighted $L^p$-norm. This way, the voice transform $\mathcal{V}_\psi : \mathcal{H} \to L^2(\mathcal{G})$ is extended to isometries of the form $\mathcal{V}_\psi : \mathcal{C}o_m^p \to L_m^p(\mathcal{G})$.

There are two technicalities that have to be taken care of. First, not every square-integrable vector $\psi \in \mathcal{Q}$ is able to give rise to such an extension of the voice transform. Square-integrability is 'good enough' to obtain a well-defined operator $\mathcal{V}_\psi : \mathcal{H} \to L^2(\mathcal{G})$, but for the more general coorbit spaces we need stronger integrability conditions. Second, if we defined the coorbit spaces as subspaces of $\mathcal{H}$, they would not be complete in general. We solve this problem by extending the voice transform to a suitable distribution space beforehand.

With these considerations in mind, we proceed in this chapter as follows.

1. We begin by defining the set of analyzing vectors $\mathcal{A}_w \subset \mathcal{Q} \subset \mathcal{H}$, which elements satisfy strong enough integrability conditions to be used in the construction.

2. Using the analyzing vectors, we define the space of test vectors $\mathcal{H}_w^1 \subset \mathcal{H}$. It is a Banach space that lies dense in $\mathcal{H}$.

3. We then define the reservoir $\mathcal{R}_w$ as the antidual space of $\mathcal{H}_w^1$. The reservoir acts like a space of distributions that contains $\mathcal{H}$ as a dense subspace, allowing us to extend the voice transform to $\mathcal{R}_w$.

4. We finally define the coorbit spaces as subspaces of the reservoir. We do this by assigning to each space $L_m^p(\mathcal{G})$ the set of all distributions $f \in \mathcal{R}_w$ for which the voice transform (with respect to an analyzing vector) is contained in $L_m^p(\mathcal{G})$.

We mainly follow the structure of Voigtlaender [31, sec. 2.4], though many statements and proofs can be replaced by simpler versions, as we are only working with weighted $L^p$-spaces in contrast to the solid quasi-Banach spaces used there. Similar ways of construction can be found in [1] or [5], where the former omits several technical details to give a more compact overview. Of course, we will also make use of the original articles by Feichtinger and Gröchenig [14; 15; 16].

For the rest of this chapter, we assume $\mathcal{G}$ is a locally compact (and $\sigma$-compact) group, and $\pi$ is a square-integrable unitary representation of $\mathcal{G}$ on the Hilbert space $\mathcal{H}$. We also assume that $w$ is a submultiplicative weight on $\mathcal{G}$ which is bounded away from zero.

## 3.1 Analyzing Vectors

We begin with the definition of analyzing vectors.

**Definition 3.1.1.** We call a vector $\psi \in \mathcal{H}$ an *analyzing vector* or *w-integrable* if $\mathcal{V}_\psi \psi \in L^1_w(\mathcal{G})$, that is if
$$\int_\mathcal{G} |\langle \psi, \pi(x)\psi \rangle_\mathcal{H}| \, w(x) \, dx < \infty.$$
We denote the set of all $w$-integrable vectors by $\mathcal{A}_w$.

The representation $\pi$ is called *w-integrable* if there exists a non-zero analyzing vector, that is if $\mathcal{A}_w \neq \{0\}$.

The $w$-integrability of $\psi \in \mathcal{H}$ is indeed stronger than the square-integrability. Since $w$ is bounded away from zero, $L^1_w(\mathcal{G})$ lies within $L^1(\mathcal{G})$. Together with the fact that $\mathcal{V}_\psi \psi$ is a bounded function on $\mathcal{G}$, it follows immediately that
$$\mathcal{V}_\psi \psi \in L^1_w(\mathcal{G}) \cap L^\infty(\mathcal{G}) \subset L^1(\mathcal{G}) \cap L^\infty(\mathcal{G}) \subset L^2(\mathcal{G}).$$

We first show that the set of analyzing vectors $\mathcal{A}_w$ is invariant under the action of $\pi$.

**Lemma 3.1.2.** *Let $\psi \in \mathcal{A}_w$. Then $\pi(x)\psi \in \mathcal{A}_w$ for all $x \in \mathcal{G}$.*

*Proof.* For $\psi \in \mathcal{A}_w$ and $x, y \in \mathcal{G}$ we have
$$\mathcal{V}_{\pi(x)\psi}(\pi(x)\psi)(y) = \langle \pi(x)\psi, \pi(y)\pi(x)\psi \rangle_\mathcal{H} = \langle \psi, \pi(x^{-1}yx)\psi \rangle_\mathcal{H}$$
$$= \mathcal{V}_\psi \psi(x^{-1}yx) = L_x R_x \mathcal{V}_\psi \psi(y).$$

Since the space $L^1_w(\mathcal{G})$ is translation invariant by proposition 2.2.8, the inequalities
$$\left\| \mathcal{V}_{\pi(x)\psi}(\pi(x)\psi) \right\|_{L^1_w} = \left\| L_x R_x \mathcal{V}_\psi \psi \right\|_{L^1_w} \leq w(x) \left\| R_x \mathcal{V}_\psi \psi \right\|_{L^1_w}$$
$$\leq w(x) w(x^{-1}) \Delta(x)^{-1} \left\| \mathcal{V}_\psi \psi \right\|_{L^1_w} < \infty$$
hold, so $\pi(x)\psi \in \mathcal{A}_w$. □

In the following theory, many definitions will depend on a previously chosen analyzing vector. To make sure all these definitions are meaningful, we assume from now on that the representation $\pi$ is indeed $w$-integrable.

We also want to show that the following theory is independent from the chosen analyzing vector, which means we will need to be able to 'switch' between two of those. The next proposition allows us to do that.

**Proposition 3.1.3.** *Let $\psi, \phi \in \mathcal{A}_w$ be two analyzing vectors. Then $\mathcal{V}_\psi \phi \in L^1_w(\mathcal{G})$.*

Before we can prove this statement, we need an auxiliary lemma.

**Lemma 3.1.4.** *[1, Lemma 3.3] Let $\psi \in \mathcal{A}_w$ and $D$ the Duflo-Moore operator from theorem 2.5.7. Then $D\psi$ is square-integrable.*

*Proof.* The set $\mathcal{Q}$ of all square-integrable vectors in $\mathcal{H}$ is exactly the domain of $D$, i.e. $\mathcal{Q} = \mathcal{D}(D)$. The domain of the adjoint operator $D^*$ of $D$ is therefore defined as

$$\mathcal{D}(D^*) = \{f \in \mathcal{H} \mid \mathcal{Q} \ni g \mapsto \langle Dg, f \rangle_\mathcal{H} \in \mathbb{C} \text{ is bounded}\}.$$

Since $D$ is self-adjoint, we have $\mathcal{Q} = \mathcal{D}(D^*)$. Hence, the square-integrability of $D\psi$ follows if we show that the anti-functional

$$\mathcal{Q} \to \mathbb{C}, \quad g \mapsto \langle Dg, D\psi \rangle_\mathcal{H}$$

is bounded.

If $D\psi = 0$, there is nothing to show, so we may assume that $D\psi \neq 0$ and therefore $\psi \neq 0$. From the Duflo-Moore theorem it follows

$$\begin{aligned}
|\langle Dg, D\psi \rangle_\mathcal{H}| &= \frac{1}{\|\psi\|_\mathcal{H}^2} |\langle \mathcal{V}_\psi \psi, \mathcal{V}_g \psi \rangle_{L^2(\mathcal{G})}| \\
&\leq \frac{1}{\|\psi\|_\mathcal{H}^2} \|\mathcal{V}_\psi \psi\|_{L^1} \|\mathcal{V}_g \psi\|_{L^\infty} \\
&\leq C \frac{\|\mathcal{V}_\psi \psi\|_{L^1_w}}{\|\psi\|_\mathcal{H}} \|g\|_\mathcal{H},
\end{aligned}$$

where we used the inequalities $\|\mathcal{V}_\psi \psi\|_{L^1} \leq C \|\mathcal{V}_\psi \psi\|_{L^1_w}$ and $\|\mathcal{V}_g \psi\|_{L^\infty} \leq \|g\|_\mathcal{H} \|\psi\|_\mathcal{H}$. This implies that the anti-functional $g \mapsto \langle Dg, D\psi \rangle_\mathcal{H}$ is indeed bounded, so $D\psi \in \mathcal{Q}$. □

We now know that the expression $D^2\psi = D(D\psi)$ is well-defined for all $\psi \in \mathcal{A}_w$. This fact is needed for the proof of proposition 3.1.3. A similar proof can be found in [31, Lemma 2.4.5].

*Proof of 3.1.3.* We may assume $\psi \neq 0 \neq \phi$, as otherwise the statement is trivial.

Consider the voice transform $\mathcal{V}_\psi(D^2\phi)$. Since the Duflo-Moore operator $D$ is injective and $\phi \neq 0$, $D^2\phi \neq 0$ as well. Furthermore, the voice transform $\mathcal{V}_\psi$ is injective, hence $\mathcal{V}_\psi(D^2\phi)$ does no vanish identically. Therefore there exists a point $x_0 \in \mathcal{G}$ for which $\mathcal{V}_\psi(D^2\phi)(x_0) \neq 0$. We define $\psi_0 = \pi(x_0)\psi$.

By lemma 3.1.2, $\psi_0$ is $w$-integrable and in particular square-integrable. Thus we can use the convolution relations from corollary 2.6.4 to obtain

$$\overline{\langle D\phi, D\psi_0 \rangle}_\mathcal{H} \mathcal{V}_{\psi_0} \phi = \mathcal{V}_\phi \phi * \mathcal{V}_{\psi_0} \psi_0$$

as well as

$$C_{\psi_0}^2 \mathcal{V}_\psi \phi = \mathcal{V}_{\psi_0} \phi * \mathcal{V}_\psi \psi_0.$$

Combining these two equations yields

$$\overline{\langle D\phi, D\psi_0 \rangle}_\mathcal{H} C_{\psi_0}^2 \mathcal{V}_\psi \phi = \mathcal{V}_\phi \phi * \mathcal{V}_{\psi_0} \psi_0 * \mathcal{V}_\psi \psi_0.$$

Now the voice transforms $\mathcal{V}_\phi \phi$ and $\mathcal{V}_{\psi_0} \psi_0$ are contained in $L^1_w(\mathcal{G})$ because $\phi$ and $\psi_0$ are analyzing vectors. Equally, $\mathcal{V}_\psi \psi_0$ is $w$-integrable, because $L^1_w(\mathcal{G})$ is translation invariant and we have $\mathcal{V}_\psi \psi_0 = \mathcal{V}_\psi(\pi(x_0)\psi) = L_{x_0} \mathcal{V}_\psi \psi$. We can therefore apply Young inequality (2.5) to see that

$$|\langle D\phi, D\psi_0 \rangle_\mathcal{H}| C_{\psi_0}^2 \|\mathcal{V}_\psi \phi\|_{L^1_w} \leq \|\mathcal{V}_\phi \phi\|_{L^1_w} \|\mathcal{V}_{\psi_0} \psi_0\|_{L^1_w} \|\mathcal{V}_\psi \psi_0\|_{L^1_w} < \infty.$$

The constants on the left side do not vanish. First, $C_{\psi_0}^2 = \langle D\psi_0, D\psi_0\rangle_\mathcal{H}$ is positive as $\psi_0$ is non-zero and $D$ is injective. Second, we have by the construction of $\psi_0$

$$\langle D\phi, D\psi_0\rangle_\mathcal{H} = \langle D^2\phi, \pi(x_0)\psi\rangle_\mathcal{H} = \mathcal{V}_\psi(D^2\phi)(x_0) \neq 0.$$

Thus it follows that $\|\mathcal{V}_\psi\phi\|_{L^1_w} < \infty$. □

*Remark 3.1.5.* Proposition 3.1.3 implies in particular that $\mathcal{A}_w$ is a vector space. $\mathcal{A}_w$ is obviously closed under multiplication with scalars. If $\psi$ and $\phi$ are $w$-integrable vectors, the proposition tells us that

$$\mathcal{V}_{\psi+\phi}(\psi + \phi) = \mathcal{V}_\psi\psi + \mathcal{V}_\psi\phi + \mathcal{V}_\phi\psi + \mathcal{V}_\phi\phi$$

is contained in $L^1_w(\mathcal{G})$, which means $\psi + \phi \in \mathcal{A}_w$.

This also shows that $\mathcal{A}_w$ lies dense in $\mathcal{H}$, as it is a non-trivial, $\pi$-invariant subspace and $\pi$ is assumed to be irreducible. ◁

## 3.2 Test Vectors

We now define the space of test vectors as a subspace of $\mathcal{H}$.

**Definition 3.2.1.** Let $\psi \in \mathcal{A}_w\backslash\{0\}$ be an analyzing vector. We define the space of *test vectors* to be
$$\mathcal{H}^1_w = \{f \in \mathcal{H} \mid \mathcal{V}_\psi f \in L^1_w(\mathcal{G})\}$$
equipped with the norm $\|f\|_{\mathcal{H}^1_w} = \|\mathcal{V}_\psi f\|_{L^1_w}$.

It is easy to see that $\|\cdot\|_{\mathcal{H}^1_w}$ indeed defines a norm. Its subadditivity and absolute homogeneity follow immediately from the linearity of $\mathcal{V}_\psi$; its positiveness follows from the injectivity of $\mathcal{V}_\psi$.

Despite the fact that the definition of $\mathcal{H}^1_w$ depends on the chosen vector $\psi$, we do not mention that $\psi$ in our notation. That is because the space of test vectors is indeed in some sense independent from $\psi$. We state this independence more precisely in the next lemma, where we write $\mathcal{H}^1_w(\psi)$ if we explicitly mean the space of test vectors characterized by $\psi$.

**Lemma 3.2.2.** *The space $\mathcal{H}^1_w$ is independent from the choice of the analyzing vector. That means: If $\psi, \phi \in \mathcal{A}_w\backslash\{0\}$ are two analyzing vectors, then $\mathcal{H}^1_w(\psi)$ and $\mathcal{H}^1_w(\phi)$ are equal as subspaces of $\mathcal{H}$ and their norms are equivalent.*

*Proof.* Both vectors $\psi$ and $\phi$ are square-integrable, which is why the convolution identity

$$C_\phi^2 \mathcal{V}_\psi f = \mathcal{V}_\phi f * \mathcal{V}_\psi\phi$$

is true for all $f \in \mathcal{H}$. It follows from Young inequality (2.5) that

$$\|\mathcal{V}_\psi f\|_{L^1_w} \leq \frac{1}{C_\phi^2}\|\mathcal{V}_\psi\phi\|_{L^1_w}\|\mathcal{V}_\phi f\|_{L^1_w}. \tag{3.1}$$

According to proposition 3.1.3 we have $\|\mathcal{V}_\psi\phi\|_{L^1_w} < \infty$, which implies the inclusion $\mathcal{H}^1_w(\phi) \subset \mathcal{H}^1_w(\psi)$.

Similarly we have the inequality

$$\left\|\mathcal{V}_\phi f\right\|_{L^1_w} \leq \frac{1}{C^2_\psi}\left\|\mathcal{V}_\phi \psi\right\|_{L^1_w}\left\|\mathcal{V}_\psi f\right\|_{L^1_w}, \qquad (3.2)$$

which implies the inclusion $\mathcal{H}^1_w(\psi) \subset \mathcal{H}^1_w(\phi)$.

Therefore, the spaces $\mathcal{H}^1_w(\psi)$ and $\mathcal{H}^1_w(\phi)$ are equal as vectors spaces. Their norm equivalence is given by the inequalities (3.1) and (3.2). □

From now on we assume $\psi \in \mathcal{A}_w \setminus \{0\}$ to be an analyzing vector and the space $\mathcal{H}^1_w$ of test vectors to be characterized by $\psi$.

Next we show some basic properties of $\mathcal{H}^1_w$. Similar statements can be found in [5, Prop. 3.2] and [1, Prop. 3.1f].

**Lemma 3.2.3.**

(i) It is $\mathcal{A}_w \subset \mathcal{H}^1_w$.

(ii) The space $\mathcal{H}^1_w$ is $\pi$-invariant and $\pi(x)$ operates bounded on $\mathcal{H}^1_w$ for every $x \in \mathcal{G}$ with the inequality $\left\|\pi(x)f\right\|_{\mathcal{H}^1_w} \leq w(x)\|f\|_{\mathcal{H}^1_w}$ for $f \in \mathcal{H}^1_w$. Additionally, the mapping $x \mapsto \pi(x)f$ is continuous for every $f \in \mathcal{H}^1_w$. In short: $\pi$ restricts to a representation of $\mathcal{G}$ on $\mathcal{H}^1_w$.

(iii) The space $\mathcal{H}^1_w$ lies dense in $\mathcal{H}$.

(iv) The injection $\mathcal{H}^1_w \hookrightarrow \mathcal{H}$ is continuous.

(v) $\mathcal{H}^1_w$ is a Banach space.

*Proof.*

(i) This immediately follows from proposition 3.1.3.

(ii) Let $x \in \mathcal{G}$ and $f \in \mathcal{H}^1_w$. According to proposition 2.2.8, the left translation $L_x$ is bounded on $L^1_w(\mathcal{G})$ by $w(x)$. Therefore we have

$$\left\|\pi(x)f\right\|_{\mathcal{H}^1_w} = \left\|\mathcal{V}_\psi(\pi(x)f)\right\|_{L^1_w} = \left\|L_x(\mathcal{V}_\psi f)\right\|_{L^1_w} \leq w(x)\left\|\mathcal{V}_\psi f\right\|_{L^1_w} = w(x)\|f\|_{\mathcal{H}^1_w},$$

which proofs the $\pi$-invariance of $\mathcal{H}^1_w$ and the boundedness of $\pi(x) : \mathcal{H}^1_w \to \mathcal{H}^1_w$ by $w(x)$.

By definition, the voice transform is an isometry $\mathcal{V}_\psi : \mathcal{H}^1_w \to L^1_w(\mathcal{G})$ which maps the representation $\pi$ to the left translation $L$. As a consequence, for any $f \in \mathcal{H}^1_w$, the continuity of $\mathcal{G} \ni x \mapsto \pi(x)f \in \mathcal{H}^1_w$ is equivalent to the continuity of

$$\mathcal{G} \ni x \mapsto \mathcal{V}_\psi(\pi(x)f) = L_x(\mathcal{V}_\psi f) \in L^1_w(\mathcal{G}).$$

But the left translation is continuous according to lemma 2.2.15, so the strong continuity of $\pi$ follows.

(iii) In remark 3.1.5 we have seen that the analyzing vectors $\mathcal{A}_w$ are dense in $\mathcal{H}$. The statement now follows from the inclusion $\mathcal{A}_w \subset \mathcal{H}^1_w$.

(iv) The weight $w$ is bounded from below by some constant $C > 0$. By using the Cauchy-Schwarz inequality, we see that

$$\|f\|_{\mathcal{H}}^2 = C_\psi^{-2} \|\mathcal{V}_\psi f\|_{L^2}^2 = C_\psi^{-2} \int_{\mathcal{G}} |\langle f, \pi(x)\psi\rangle_{\mathcal{H}}| \, |\mathcal{V}_\psi f(x)| \, dx$$
$$\leq C_\psi^{-2} C^{-1} \int_{\mathcal{G}} \|f\|_{\mathcal{H}} \|\psi\|_{\mathcal{H}} \, |\mathcal{V}_\psi f(x)| w(x) \, dx$$
$$= C_\psi^{-2} C^{-1} \|f\|_{\mathcal{H}} \|\psi\|_{\mathcal{H}} \|\mathcal{V}_\psi f\|_{L^1_w}$$
$$= C_\psi^{-2} C^{-1} \|f\|_{\mathcal{H}} \|\psi\|_{\mathcal{H}} \|f\|_{\mathcal{H}^1_w}.$$

Hence it is $\|f\|_{\mathcal{H}} \leq C_\psi^{-2} C^{-1} \|\psi\|_{\mathcal{H}} \|f\|_{\mathcal{H}^1_w}$, which implies the continuity of the injection $\mathcal{H}^1_w \hookrightarrow \mathcal{H}$.

(v) Let $(f_n)_{n\in\mathbb{N}} \subset \mathcal{H}^1_w$ be Cauchy. By definition, $(\mathcal{V}_\psi f_n)_{n\in\mathbb{N}}$ is then a Cauchy sequence in $L^1_w(\mathcal{G})$. In this Banach space the sequence converges to some $F \in L^1_w(\mathcal{G})$. By taking a subsequence, we may assume that this convergence also happens pointwise almost everywhere.

Similarly, $(f_n)_{n\in\mathbb{N}}$ is Cauchy in $\mathcal{H}$, as we have already shown that the injection $\mathcal{H}^1_w \hookrightarrow \mathcal{H}$ is continuous. Therefore, the sequence converges in $\mathcal{H}$ to some $f \in \mathcal{H}$. This also implies that $(\mathcal{V}_\psi f_n)_{n\in\mathbb{N}}$ converges to $\mathcal{V}_\psi f$ as elements of $L^2(\mathcal{G})$. Again, we may assume this convergence to happen pointwise almost everywhere by taking a subsequence.

Now the sequence $(\mathcal{V}_\psi f_n)_{n\in\mathbb{N}}$ converges pointwise almost everywhere to $F \in L^1_w(\mathcal{G})$ and $\mathcal{V}_\psi f \in L^2(\mathcal{G})$, so we have $\mathcal{V}_\psi f = F$ almost everywhere. This implies $f \in \mathcal{H}^1_w$, and that $(\mathcal{V}_\psi f_n)_{n\in\mathbb{N}}$ converges to $\mathcal{V}_\psi f$ in the $L^1_w$-norm. This is equivalent to the convergence of $(f_n)_{n\in\mathbb{N}}$ to $f$ in the $\mathcal{H}^1_w$ norm. The space $\mathcal{H}^1_w$ is therefore complete.

$\square$

By definition, the space of test vectors contains exactly those $f \in \mathcal{H}$ for which $\mathcal{V}_\psi f$ is a $w$-integrable function. We can use this integrability condition to interpret the reconstruction formula

$$C_\psi^2 f = \int_{\mathcal{G}} \mathcal{V}_\psi f(x) \, \pi(x)\psi \, dx,$$

which is only well-defined in the weak sense in general, as a strong integral for $f \in \mathcal{H}^1_w$.

**Proposition 3.2.4** (strong reconstruction formula, [31, Lemma 2.4.7(4)]). *For every $f \in \mathcal{H}^1_w$, the function*

$$\varphi : \mathcal{G} \to \mathcal{H}^1_w, \quad \varphi(x) = \mathcal{V}_\psi f(x) \, \pi(x)\psi$$

*is strongly integrable, and we have the strong reconstruction formula*

$$C_\psi^2 f = \int_{\mathcal{G}} \varphi(x) \, dx = \int_{\mathcal{G}} \mathcal{V}_\psi f(x) \, \pi(x)\psi \, dx. \quad (3.3)$$

*Proof.* We first show that the integral (3.3) is indeed well-defined. We notice that $\varphi$ is continuous, as it is the product of continuous functions. According to remark 2.3.2, $\varphi$ is therefore $\mu$-measurable. Now we have for $x \in \mathcal{G}$ the inequality

$$\|\varphi(x)\|_{\mathcal{H}^1_w} = |\mathcal{V}_\psi f(x)| \, \|\pi(x)\psi\|_{\mathcal{H}^1_w} \leq w(x) \|\psi\|_{\mathcal{H}^1_w} |\mathcal{V}_\psi f(x)|.$$

It follows from $\mathcal{V}_\psi f \in L^1_w(\mathcal{G})$ that the function $\|\varphi\|_{\mathcal{H}^1_w}$ is integrable. By proposition 2.3.1, this is sufficient for $\varphi$ to be strongly integrable.

The strong integral $\Phi = \int_\mathcal{G} \varphi(x)\,dx$ is now a well-defined element of $\mathcal{H}^1_w$. Additionally, the same integral is strongly convergent in $\mathcal{H}$, as the inclusion map $\mathcal{H}^1_w \hookrightarrow \mathcal{H}$ is continuous and the strong integral is compatible with continuous maps (proposition 2.3.3). But the weak reconstruction formula from proposition 2.6.2 already tells us that the value of this integral has to be $C^2_\psi f$, as the values of the strong and weak integrals agree if the integral converges in both senses. Hence we have $\Phi = C^2_\psi f$. □

According to lemma 3.2.3, the restriction of $\pi$ to $\mathcal{H}^1_w$ is again a representation. This restriction is in general not irreducible, which means there are in general non-cyclic vectors in $\mathcal{H}^1_w$. Nonetheless, we can use the strong reconstruction formula to prove that at least all non-trivial analyzing vectors are cyclic in $\mathcal{H}^1_w$.

**Corollary 3.2.5.** *[31, Lemma 2.4.7(5)] Every w-integrable vector $\psi \in \mathcal{A}_w\setminus\{0\}$ is cyclic in $\mathcal{H}^1_w$. In particular, $\mathcal{A}_w$ is dense in $\mathcal{H}^1_w$.*

*Proof.* Every $f \in \mathcal{H}^1_w$ can be written as a strong integral
$$f = \frac{1}{C^2_\psi} \int_\mathcal{G} \mathcal{V}_\psi f(x)\,\pi(x)\psi\,dx.$$
The integrand is pointwise contained in
$$\overline{\mathcal{E}}_\psi = \overline{\text{span}\{\pi(y)\psi \mid y \in \mathcal{G}\}},$$
which itself is a Banach space as a closed subspace of $\mathcal{H}^1_w$. Now the strong integral is defined on this space $\overline{\mathcal{E}}_\psi$, which means its value $f$ is also contained in $\overline{\mathcal{E}}_\psi$.

Therefore, we have $\mathcal{H}^1_w \subset \overline{\mathcal{E}}_\psi$, so $\psi$ is cyclic in $\mathcal{H}^1_w$. Since $\mathcal{A}_w$ is a $\pi$-invariant vector space, we also have $\mathcal{E}_\psi \subset \mathcal{A}_w$, thus $\mathcal{A}_w$ is dense in $\mathcal{H}^1_w$. □

We end this section with an interesting observation regarding the relation between analyzing vectors and test vectors.

**Proposition 3.2.6.** *Suppose the weight $w$ satisfies $w(x) = \Delta(x)^{-1}w(x^{-1})$ for all $x \in \mathcal{G}$. Then $\mathcal{A}_w = \mathcal{H}^1_w$ as sets.*

*Proof.* This equality comes down to the fact that $L^1_v(\mathcal{G})^\vee = L^1_{\Delta^{-1}v^\vee}(\mathcal{G})$ for general weights $v$, since the substitution $x \mapsto x^{-1}$ with $d(x^{-1}) = \Delta(x)^{-1}\,dx$ in the $L^1_v$ norm yields
$$\|G^\vee\|_{L^1_v} = \int_\mathcal{G} |G(x^{-1})|v(x)\,dx = \int_\mathcal{G} |G(x)|\Delta(x)^{-1}v(x^{-1})\,dx = \|G\|_{L^1_{\Delta^{-1}v^\vee}}$$
for any measurable $G: \mathcal{G} \to \mathbb{C}$.

The requirement $w = \Delta^{-1}w^\vee$ therefore implies that $\|G\|_{L^1_w} = \|G^\vee\|_{L^1_w}$ for all $G \in L^1_w(\mathcal{G}) = L^1_w(\mathcal{G})^\vee$, and in particular
$$\|\mathcal{V}_\psi g\|_{L^1_w} = \|\mathcal{V}_\psi g^\triangledown\|_{L^1_w} = \|\mathcal{V}_g \psi\|_{L^1_w}$$
for all $\psi \in \mathcal{A}_w$ and $g \in \mathcal{H}^1_w$. Using the reproducing formula (2.17) and Young inequality (2.5), it follows that
$$\|\mathcal{V}_g g\|_{L^1_w} = C^{-2}_\psi \|\mathcal{V}_\psi g * \mathcal{V}_g \psi\|_{L^1_w} \leq C^{-2}_\psi \|\mathcal{V}_\psi g\|_{L^1_w} \|\mathcal{V}_g \psi\|_{L^1_w} = C^{-2}_\psi \|\mathcal{V}_\psi g\|^2_{L^1_w} < \infty.$$

Thus we have $\mathcal{H}^1_w \subset \mathcal{A}_w$. The other inclusion is already known. □

3.3. The Reservoir and the Extended Voice Transform    36## 3.3 The Reservoir and the Extended Voice Transform

We can now define the distribution space that is suitable in our context.

**Definition 3.3.1.** We define the *reservoir* as the antidual space $\mathcal{R}_w = (\mathcal{H}_w^1)^\sim$ of $\mathcal{H}_w^1$.

The antidual space is the space of all antilinear bounded functionals (*anti-functionals*), that is, the space of all bounded operators $\alpha : \mathcal{H}_w^1 \to \mathbb{C}$ for which $\alpha(\lambda g_1 + g_2) = \overline{\lambda}\alpha(g_1) + \alpha(g_2)$ holds for all $\lambda \in \mathbb{C}$ and $g_1, g_2 \in \mathcal{H}_w^1$.

Since the definition of the space $\mathcal{H}_w^1$ depends on the choice of an analyzing vector, so does the reservoir $\mathcal{R}_w$. However, given two analyzing vectors $\psi, \phi \in \mathcal{A}_w\setminus\{0\}$, the spaces $\mathcal{R}_w(\psi)$ and $\mathcal{R}_w(\phi)$ are isomorphic as normed spaces, as the same is already true for $\mathcal{H}_w^1(\psi)$ and $\mathcal{H}_w^1(\phi)$. The two reservoirs can even be constructed to be equal as vector spaces, equipped with equivalent norms.

The spaces $\mathcal{H}_w^1(\psi)$ and $\mathcal{H}_w^1(\phi)$ are equal as vector space. Thus $\mathcal{R}_w(\psi)$ and $\mathcal{R}_w(\phi)$ are both subspaces of the algebraic antidual

$$L = \{\alpha : \mathcal{H}_w^1 \to \mathbb{C} \mid \alpha \text{ is antilinear}\}.$$

Now the norm equivalence of $\mathcal{H}_w^1(\psi)$ and $\mathcal{H}_w^1(\phi)$ implies that an antilinear map $\alpha \in L$ is bounded with respect to the $\mathcal{H}_w^1(\psi)$-norm if and only if it is bounded with respect to the $\mathcal{H}_w^1(\phi)$-norm. Therefore, the spaces $\mathcal{R}_w(\psi)$ and $\mathcal{R}_w(\phi)$ are equal as vector spaces (and their norms are equivalent).

We will always assume that both spaces $\mathcal{H}_w^1$ and $\mathcal{R}_w$ are characterized by the same analyzing vector $\psi$. Then we can use the usual duality relations

$$\|\alpha\|_{\mathcal{R}_w} = \sup_{g \in \mathcal{H}_w^1\setminus\{0\}} \frac{|\alpha(g)|}{\|g\|_{\mathcal{H}_w^1}}$$

and

$$\|g\|_{\mathcal{H}_w^1} = \sup_{\alpha \in \mathcal{R}_w\setminus\{0\}} \frac{|\alpha(g)|}{\|\alpha\|_{\mathcal{R}_w}}.$$

We will now show that the Hilbert space $\mathcal{H}$ can be embedded into $\mathcal{R}_w$ continuously, and can therefore be understood to be a subspace of $\mathcal{R}_w$. To do that, we need to find a linear map $\mathcal{H} \hookrightarrow \mathcal{R}_w$ that is injective and bounded.

Given a vector $f \in \mathcal{H}$, we define the anti-functional $\alpha_f : \mathcal{H}_w^1 \to \mathbb{C}$ by $\alpha_f(g) = \langle f, g\rangle_\mathcal{H}$. Since by lemma 3.2.3 the injection $\mathcal{H}_w^1 \hookrightarrow \mathcal{H}$ is bounded by some $C > 0$, we have

$$|\alpha_f(g)| = |\langle f, g\rangle_\mathcal{H}| \leq \|f\|_\mathcal{H}\|g\|_\mathcal{H} \leq C\|f\|_\mathcal{H}\|g\|_{\mathcal{H}_w^1}, \tag{3.4}$$

so $\alpha_f$ is bounded and an element of $\mathcal{R}_w$.

The map $\imath : \mathcal{H} \to \mathcal{R}_w$, $f \mapsto \alpha_f$ is now well-defined and obviously linear. The inequality (3.4) implies

$$\|\alpha_f\|_{\mathcal{R}_w} = \sup_{g \in \mathcal{H}_w^1\setminus\{0\}} \frac{|\alpha_f(g)|}{\|g\|_{\mathcal{H}_w^1}} \leq C\|f\|_\mathcal{H},$$

thus $\imath$ is bounded. To show that $\imath$ is injective, we assume that for some $f \in \mathcal{H}$ the evaluation $\alpha_f(g) = \langle f, g\rangle_\mathcal{H} = 0$ for every $g \in \mathcal{H}_w^1$. But $\mathcal{H}_w^1$ is dense in $\mathcal{H}$ (Lemma 3.2.3), so $f$ must be zero already.

Furthermore, the image of $\imath$ is dense in $\mathcal{R}_w$. If we assume $g \in \mathcal{H}_w^1$ has the property that $\alpha_f(g) = \langle f, g\rangle_\mathcal{H} = 0$ for all $f \in \mathcal{H}$, then of course $g = 0$.

We already knew that $\mathcal{H}_w^1$ is densely embedded into $\mathcal{H}$. Overall, we have the following statement.

**Proposition 3.3.2.** *The injections $\mathcal{H}_w^1 \hookrightarrow \mathcal{H} \hookrightarrow \mathcal{R}_w$ are continuous and dense embeddings.*

From now on, we will identify vectors $f \in \mathcal{H}$ with their associated anti-functionals $\alpha_f$. We will note the elements of $\mathcal{R}_w$ as $f, h, \ldots$ and write the antidual pairing of $\mathcal{H}_w^1$ and $\mathcal{R}_w$ like an inner product
$$\alpha_f(g) = f(g) = \langle f, g \rangle_{\mathcal{R}_w \times \mathcal{H}_w^1}.$$
Then we have for $f \in \mathcal{H}$ and $g \in \mathcal{H}_w^1$
$$\langle f, g \rangle_{\mathcal{H}} = \langle f, g \rangle_{\mathcal{R}_w \times \mathcal{H}_w^1}, \tag{3.5}$$
so the antiduality of $\mathcal{R}_w$ and $\mathcal{H}_w^1$ is an extension of the inner product on $\mathcal{H}$.

*Remark 3.3.3.* It is not really important whether we use the dual $(\mathcal{H}_w^1)'$ or the antidual $(\mathcal{H}_w^1)^\sim$ of the space of test vectors as the reservoir. However, using the dual space leads to the embedding
$$\imath_0 : \mathcal{H} \to (\mathcal{H}_w^1)', \quad f \mapsto \beta_f = \langle \cdot, f \rangle_{\mathcal{H}}|_{\mathcal{H}_w^1},$$
which is antilinear. Thus the identification $f = \beta_f$ is not compatible with scalar multiplication, which means we would need to strictly differentiate between $f$ and $\beta_f$. We avoid this technical inconvenience by using the antidual space, as the identification $f = \alpha_f$ is compatible with the linear structure on $\mathcal{H}$ and $\mathcal{R}_w$. ◁

We can use the extension of the inner product of $\mathcal{H}$ to also extend the representation $\pi$ to the reservoir. For $f \in \mathcal{R}_w$ and $x \in \mathcal{G}$, we define the extension of $\pi(x)$ as
$$\langle \pi(x)f, g \rangle_{\mathcal{R}_w \times \mathcal{H}_w^1} = \langle f, \pi(x^{-1})g \rangle_{\mathcal{R}_w \times \mathcal{H}_w^1} \quad \text{for } g \in \mathcal{H}_w^1.$$
For $f \in \mathcal{H}$, this definition is compatible with the initial one since $\pi$ is unitary.

As an operator, $\pi(x) : \mathcal{R}_w \to \mathcal{R}_w$ is just the adjoint of the restriction $\pi(x^{-1})|_{\mathcal{H}_w^1} : \mathcal{H}_w^1 \to \mathcal{H}_w^1$. This restriction is bounded by lemma 3.2.3(ii). As its adjoint, $\pi(x)$ acts bounded on $\mathcal{R}_w$.

We can also extend the voice transform to the reservoir.

**Definition 3.3.4.** Let $g \in \mathcal{H}_w^1$ and $f \in \mathcal{R}_w$. We define the *extended voice transform* by
$$\mathcal{V}_g f(x) = \langle f, \pi(x)g \rangle_{\mathcal{R}_w \times \mathcal{H}_w^1}, \quad x \in \mathcal{G}.$$

The extended voice transform is well defined for all $x \in \mathcal{G}$ because $\mathcal{H}_w^1$ is $\pi$-invariant, see 3.2.3(ii). For $f \in \mathcal{H}$, it coincides with the original voice transform, as the antidual pairing coincides with the inner product.

For $f \in \mathcal{R}_w$ and $g \in \mathcal{H}_w^1$, the transform $\mathcal{V}_g f$ is a continuous function on $\mathcal{G}$ since $\pi$ operates strongly continuous on $\mathcal{H}_w^1$. Additionally, the inequality
$$|\mathcal{V}_g f(x)| \leq \|f\|_{\mathcal{R}_w} \|\pi(x)g\|_{\mathcal{H}_w^1} \leq w(x) \|f\|_{\mathcal{R}_w} \|g\|_{\mathcal{H}_w^1}$$
for $x \in \mathcal{G}$ implies that $\mathcal{V}_g f \in L_{1/w}^\infty(\mathcal{G})$. It particularly follows that
$$\mathcal{V}_g : \mathcal{R}_w \to L_{1/w}^\infty(\mathcal{G})$$
is a well-defined operator that is bounded by $\|g\|_{\mathcal{H}_w^1}$.

Many formulas that are true for the voice transform carry over to the extended version. The inner product on $\mathcal{H}$ has then to be replaced by the antidual pairing of $\mathcal{R}_w$ and $\mathcal{H}^1_w$. Similarly, the inner product on $L^2(\mathcal{G})$ has to be replaced by the antidual pairing of $L^\infty_{1/w}(\mathcal{G})$ and $L^1_w(\mathcal{G})$ that is given by

$$\langle F, G \rangle_{L^\infty_{1/w} \times L^1_w} = \int_\mathcal{G} F(x)\overline{G(x)}\,dx$$

for $F \in L^\infty_{1/w}(\mathcal{G})$ and $G \in L^1_w(\mathcal{G})$ (cf. lemma 2.2.2).

**Proposition 3.3.5.** *Let $f \in \mathcal{R}_w$, $g \in \mathcal{H}^1_w$ and $\psi, \phi \in \mathcal{A}_w$.*

(i) *For $x \in \mathcal{G}$ we have $\mathcal{V}_g(\pi(x)f) = L_x \mathcal{V}_g f$.*

(ii) *If $D$ is the Duflo-Moore operator, we have*

$$\langle \mathcal{V}_\psi f, \mathcal{V}_\phi g \rangle_{L^\infty_{1/w} \times L^1_w} = \overline{\langle D\psi, D\phi \rangle_\mathcal{H}} \langle f, g \rangle_{\mathcal{R}_w \times \mathcal{H}^1_w}. \tag{3.6}$$

*If $\psi$ is admissible, we have in particular*

$$\langle f, g \rangle_{\mathcal{R}_w \times \mathcal{H}^1_w} = \int_\mathcal{G} \mathcal{V}_\psi f(x) \overline{\mathcal{V}_\psi g(x)}\,dx.$$

(iii) *We have the convolution identity*

$$\mathcal{V}_\psi f * \mathcal{V}_\phi g^\nabla = \overline{\langle D\psi, D\phi \rangle_\mathcal{H}} \mathcal{V}_g f,$$

*where the convolution integral is absolutely convergent in every point.*

(iv) *The mapping $f \mapsto \|\mathcal{V}_\psi f\|_{L^\infty_{1/w}}$ defines an equivalent norm on $\mathcal{R}_w$. More specifically, we have the inequalities*

$$C_\psi^2 \|f\|_{\mathcal{R}_w} \le \|\mathcal{V}_\psi f\|_{L^\infty_{1/w}} \le \|\psi\|_{\mathcal{H}^1_w} \|f\|_{\mathcal{R}_w}$$

*for all $f \in \mathcal{R}_w$. The extended voice transform $\mathcal{V}_\psi : \mathcal{R}_w \to L^\infty_{1/w}(\mathcal{G})$ is in particular injective.*

*Proof.*

(i) This follows from the same calculation as for the original voice transform:

$$\mathcal{V}_g(\pi(x)f)(y) = \langle \pi(x)f, \pi(y)g \rangle_{\mathcal{R}_w \times \mathcal{H}^1_w} = \langle f, \pi(x^{-1}y)g \rangle_{\mathcal{R}_w \times \mathcal{H}^1_w} = L_x \mathcal{V}_g f(y)$$

holds for all $y \in \mathcal{G}$.

(ii) The (extended) voice transform defines a bounded operator as $\mathcal{V}_\psi : \mathcal{R}_w \to L^\infty_{1/w}(\mathcal{G})$, and $\mathcal{V}_\phi g$ is contained in $L^1_w(\mathcal{G})$. Thus the left side of (3.6) is a bounded functional in $f$. The same is true for the right side.

According to the Duflo-Moore theorem 2.5.7, the equation is certainly true for all $f \in \mathcal{H}$. Since $\mathcal{H}$ is a dense subspace of $\mathcal{R}_w$ and both sides are bounded functionals in $f$, the equality holds for all $f \in \mathcal{R}_w$.

(iii) As a consequence of (i) and (ii) we have

$$\begin{aligned}(\mathcal{V}_\psi f * \mathcal{V}_\phi g^\triangledown)(x) &= \int_\mathcal{G} \mathcal{V}_\psi f(y)\overline{\mathcal{V}_\phi g(x^{-1}y)}\,dy \\ &= \langle \mathcal{V}_\psi f, \mathcal{V}_\phi(\pi(x)g)\rangle_{L^\infty_{1/w} \times L^1_w} \\ &= \overline{\langle D\psi, D\phi\rangle}_\mathcal{H} \langle f, \pi(x)g\rangle_{\mathcal{R}_w \times \mathcal{H}^1_w} \\ &= \overline{\langle D\psi, D\phi\rangle}_\mathcal{H} \mathcal{V}_g f(x)\end{aligned}$$

for all $x \in \mathcal{G}$. The integral are absolutely convergent since $L^1_w(\mathcal{G})$ and $L^\infty_{1/w}(\mathcal{G})$ are (anti)dual to each other.

(iv) The second inequality follows from the boundedness of $\mathcal{V}_\psi$ by $\|\psi\|_{\mathcal{H}^1_w}$, so it remains to show the first one.

Let $f \in \mathcal{R}_w$ and $g \in \mathcal{H}^1_w$. Then the extension of the Duflo-Moore theorem in (ii) implies

$$C^2_\psi |\langle f, g\rangle_{\mathcal{R}_w \times \mathcal{H}^1_w}| = |\langle \mathcal{V}_\psi f, \mathcal{V}_\psi g\rangle_{L^\infty_{1/w} \times L^1_w}| \leq \|\mathcal{V}_\psi f\|_{L^\infty_{1/w}} \|\mathcal{V}_\psi g\|_{L^1_w} = \|g\|_{\mathcal{H}^1_w} \|\mathcal{V}_\psi f\|_{L^\infty_{1/w}}.$$

This estimate holds for all $g \in \mathcal{H}^1_w$, thus by duality we have $C^2_\psi \|f\|_{\mathcal{R}_w} \leq \|\mathcal{V}_\psi f\|_{L^\infty_{1/w}}$.

□

**Remark 3.3.6.** The extended voice transform $\mathcal{V}_g$ for $g \in \mathcal{H}^1_w$ is injective if and only if $g$ is cyclic in $\mathcal{H}^1_w$ (see lemma 2.5.8 for the corresponding statement on $\mathcal{H}$). Thus, we could also use corollary 3.2.5 to show the injectivity of $\mathcal{V}_\psi$ for $\psi \in \mathcal{A}_w$. Since the restriction of $\pi$ to $\mathcal{H}^1_w$ is in general no longer irreducible, the voice transform $\mathcal{V}_g$ is not injective for arbitrary $g \in \mathcal{H}^1_w$. ◁

The convolution identity from the preceding proposition implies in particular that corollary 2.6.4 is still true for the extended voice transform.

**Corollary 3.3.7.** *For $g \in \mathcal{H}^1_w$ and $f \in \mathcal{R}_w$ we have*

$$C^2_\psi \mathcal{V}_g f = \mathcal{V}_\psi f * \mathcal{V}_g \psi.$$

*If $\psi$ is admissible, we have in particular the (extended) reproducing formula*

$$\mathcal{V}_\psi f = \mathcal{V}_\psi f * \mathcal{V}_\psi \psi.$$

In section 2.6, we used the reproducing formula to identify the image of $\mathcal{V}_\psi : \mathcal{H} \to L^2(\mathcal{G})$ as the reproducing kernel Hilbert space

$$\mathcal{M}^2 = \{F \in L^2(\mathcal{G}) \mid F = F * \mathcal{V}_\psi \psi\}.$$

A similar statement is true for the extended voice transform.

**Theorem 3.3.8** (Correspondence principle for the reservoir, [31, Lemma 2.4.8(5)]).
*Let $\psi \in \mathcal{A}_w \setminus \{0\}$ be an admissible analyzing vector. Then*

$$P : L^\infty_{1/w}(\mathcal{G}) \to \mathcal{V}_\psi(\mathcal{R}_w) \subset L^\infty_{1/w}(\mathcal{G}), \quad F \mapsto F * \mathcal{V}_\psi \psi$$

*is a projection onto the image of the extended voice transform $\mathcal{V}_\psi$. In particular, we have the equivalence for $F \in L^\infty_{1/w}(\mathcal{G})$*

$$F = \mathcal{V}_\psi f \text{ for some } f \in \mathcal{R}_w \quad \Longleftrightarrow \quad F = F * \mathcal{V}_\psi \psi.$$

*Proof.* Every $F \in \mathcal{V}_\psi(\mathcal{R}_w)$ satisfies the reproducing formula as we have seen in Corollary 3.3.7. If we show that every $F \in L^\infty_{1/w}(\mathcal{G})$ that satisfies $F = F * \mathcal{V}_\psi \psi$ is of the form $\mathcal{V}_\psi f$, then it is clear that $P$ indeed defines a projection onto $\mathcal{V}_\psi(\mathcal{R}_w)$.

The restricted voice transform $\mathcal{V}_\psi|_{\mathcal{H}^1_w} : \mathcal{H}^1_w \to L^1_w(\mathcal{G})$ is an isometry. The adjoint of this operator is given by

$$S = (\mathcal{V}_\psi|_{\mathcal{H}^1_w})^\sim : L^\infty_{1/w}(\mathcal{G}) \to \mathcal{R}_w, \quad \langle SF, g\rangle_{\mathcal{R}_w \times \mathcal{H}^1_w} = \langle F, \mathcal{V}_\psi g\rangle_{L^\infty_{1/w} \times L^1_w}.$$

The operator $S$ is bounded and we have for $F \in L^\infty_{1/w}(\mathcal{G})$ and $x \in \mathcal{G}$

$$\begin{aligned}
\mathcal{V}_\psi(SF)(x) &= \langle SF, \pi(x)\psi\rangle_{\mathcal{R}_w \times \mathcal{H}^1_w} \\
&= \int_\mathcal{G} F(y)\overline{\mathcal{V}_\psi(\pi(x)\psi)(y)}\, dy \\
&= \int_\mathcal{G} F(y)\mathcal{V}_\psi\psi(y^{-1}x)\, dy \\
&= (F * \mathcal{V}_\psi\psi)(x).
\end{aligned}$$

Therefore it is $P = \mathcal{V}_\psi S$, which in particular implies that $P$ is bounded.

If we now assume that $F = F * \mathcal{V}_\psi\psi$, it follows $F = \mathcal{V}_\psi(SF)$ with $SF \in \mathcal{R}_w$. This completes the proof. □

*Remark 3.3.9.* We proved the unitary correspondence principle by considering the concatenation of $\mathcal{V}_\psi : \mathcal{H} \to L^2(\mathcal{G})$ and its adjoint $\mathcal{V}^*_\psi$. Since $\mathcal{V}_\psi$ is an isometry, $\mathcal{V}_\psi \mathcal{V}^*_\psi$ is an orthogonal projection onto its image. The preceding proof uses the same idea, though we needed to adjust the domain and codomain of the voice transform accordingly. ◁

## 3.4 Coorbit Spaces and their Properties

We are now able to define the coorbit spaces as subspaces of $\mathcal{R}_w$. Their definition depends on a weighted $L^p$ space $L^p_m(\mathcal{G})$ which is left- and right-invariant and a left and right $L^1_w(\mathcal{G})$-module. Consequently, we assume that $m$ is a moderate weight function on $\mathcal{G}$, $1 \leq p \leq \infty$, and $w$ is a $p$-control-weight of $m$. As usual, we assume that $\psi \in \mathcal{A}_w \setminus \{0\}$ is an analyzing vector as well.

**Definition 3.4.1.** We define the *coorbit space of* $L^p_m(\mathcal{G})$ as the space

$$\mathcal{C}o^p_m = \mathcal{C}o(L^p_m(\mathcal{G})) = \{f \in \mathcal{R}_w \mid \mathcal{V}_\psi f \in L^p_m(\mathcal{G})\}$$

equipped with the norm $\|f\|_{\mathcal{C}o^p_m} = \|\mathcal{V}_\psi f\|_{L^p_m}$.

We first note that $\|\cdot\|_{\mathcal{C}o^p_m}$ indeed defines a norm on $\mathcal{C}o^p_m$, as all the required properties immediately carry over from $L^p_m(\mathcal{G})$ through the injectivity and linearity of $\mathcal{V}_\psi$. We also note that $\mathcal{V}_\psi : \mathcal{C}o^p_m \to L^p_m(\mathcal{G})$ is by definition an isometry.

The coorbit space $\mathcal{C}o^p_m$ is $\pi$-invariant. This follows from the left invariance of $L^p_m(\mathcal{G})$ in proposition 2.2.5, together with the computation

$$\|\pi(x)f\|_{\mathcal{C}o^p_m} = \|\mathcal{V}_\psi(\pi(x))f\|_{L^p_m} = \|L_x \mathcal{V}_\psi f\|_{L^p_m} \leq w(x)\|\mathcal{V}_\psi f\|_{L^p_m} = w(x)\|f\|_{\mathcal{C}o^p_m}$$

for $x \in \mathcal{G}$ and $f \in \mathcal{C}o^p_m$. According to lemma 2.2.15, the mapping $x \mapsto \pi(x)f$ acts for $p < \infty$ continuously in the $\mathcal{C}o^p_m$-norm.

Similar to the space of test vectors, the coorbit spaces are independent from the choice of analyzing vector.

**Lemma 3.4.2.** *If $\phi \in \mathcal{A}_w \setminus \{0\}$ is another analyzing vector, then the spaces $\mathcal{C}o_m^p(\psi)$ and $\mathcal{C}o_m^p(\phi)$ coincide as subspaces of $\mathcal{R}_w$ and their norms are equivalent.*

*Proof.* The proof is similar to that of lemma 3.2.2.

According to proposition 3.1.3, $\mathcal{V}_\psi \phi$ and $\mathcal{V}_\phi \psi$ are contained in $L_w^1(\mathcal{G})$. By applying Young inequality (2.8) to the convolution identity from Corollary 3.3.7, we get

$$C_\psi^2 \|\mathcal{V}_\phi f\|_{L_m^p} \leq \|\mathcal{V}_\psi f\|_{L_m^p} \|\mathcal{V}_\phi \psi\|_{L_w^1}$$

and vice versa

$$C_\phi^2 \|\mathcal{V}_\psi f\|_{L_m^p} \leq \|\mathcal{V}_\phi f\|_{L_m^p} \|\mathcal{V}_\psi \phi\|_{L_w^1}$$

for all $f \in \mathcal{R}_w$. Therefore we have $\mathcal{V}_\psi f \in L_m^p(\mathcal{G})$ if and only if $\mathcal{V}_\phi f \in L_m^p(\mathcal{G})$, which implies the equality of $\mathcal{C}o_m^p(\psi)$ and $\mathcal{C}o_m^p(\phi)$ as vector spaces. The two inequalities then describe their norm equivalence. □

We want to show some important properties of the coorbit spaces in this section. This includes their completeness, their duality relations, as well as a version of the correspondence principle. In order to to that, we first show that the reproducing kernel $\mathcal{V}_\psi \psi$ is contained in $L_m^p(\mathcal{G})$. Since control-weights are compatible with duality, i.e. $w$ is also a $q$-control-weight of $1/m$, this also proves that $\mathcal{V}_\psi \psi \in L_{1/m}^q(\mathcal{G})$; this fact will be essential in all of the following theory.

The method is motivated by [15, Prop. 4.3(iii)]. The construction of the vector $\phi$ in the proof is taken from [14, Lemma 6.1].

**Lemma 3.4.3.** *We have $\mathcal{A}_w \subset \mathcal{C}o_m^p$.*

*Proof.* Let $\varphi \in C_c(\mathcal{G})$ be a continuous function with compact support for which the convolution $\varphi * \mathcal{V}_\psi \psi$ does not vanish identically. Since the space $C_c(\mathcal{G})$ is dense in $L_w^1(\mathcal{G})$ (cf. lemma 2.2.16) and the convolution $F \mapsto F * \mathcal{V}_\psi \psi$ is not the zero operator on $L_w^1(\mathcal{G})$, such a $\varphi$ does indeed exist. According to the unitary correspondence principle 2.6.5, there exists now a $\phi \in \mathcal{H} \setminus \{0\}$ such that

$$\mathcal{V}_\psi \phi = \varphi * \mathcal{V}_\psi \psi.$$

Using the convolution identity (2.17) we can write $\mathcal{V}_\phi \phi$ as

$$\mathcal{V}_\phi \phi = C_\psi^{-2} \mathcal{V}_\psi \phi * \mathcal{V}_\psi \phi^\nabla$$
$$= C_\psi^{-2} \varphi * \mathcal{V}_\psi \psi * (\varphi * \mathcal{V}_\psi \psi)^\nabla$$
$$= C_\psi^{-2} \varphi * \mathcal{V}_\psi \psi * \mathcal{V}_\psi \psi^\nabla * \varphi^\nabla$$
$$= C_\psi^{-2} \varphi * \mathcal{V}_\psi \psi * \varphi^\nabla,$$

where we used the fact that $(F * G)^\nabla = G^\nabla * F^\nabla$.

Now $\varphi$ and $\varphi^\nabla$ are contained in $C_c(\mathcal{G}) \subset L_w^1(\mathcal{G}) \cap L_m^p(\mathcal{G})$. Thus, the Young inequality (2.5) implies

$$\|\mathcal{V}_\phi \phi\|_{L_w^1} \leq C_\psi^{-2} \|\varphi\|_{L_w^1} \|\mathcal{V}_\psi \psi\|_{L_w^1} \|\varphi^\nabla\|_{L_w^1} < \infty,$$

which implies that $\phi$ is an analyzing vector. Using Young inequality (2.7) and the equation $\mathcal{V}_\phi \psi = \mathcal{V}_\psi \psi * \varphi^\nabla$, we also get

$$\|\mathcal{V}_\phi \psi\|_{L_m^p} \leq \|\mathcal{V}_\psi \psi\|_{L_w^1} \|\varphi^\nabla\|_{L_m^p} < \infty,$$

hence $\mathcal{V}_\phi \psi \in L^p_m(\mathcal{G})$. Since $\phi$ is an analyzing vector and can therefore be used to characterize coorbit spaces, it follows that $\psi \in \mathcal{C}o^p_m$. But $\psi$ is an arbitrary analyzing vector, so we have shown that $\mathcal{A}_w \subset \mathcal{C}o^p_m$. $\square$

Using the above lemma, we are immediately able to prove one of the most important theorems of coorbit theory: the general correspondence principle.

**Theorem 3.4.4** (correspondence principle, [15, Prop. 4.3]). *Let $\psi \in \mathcal{A}_w \backslash \{0\}$ be an admissible analyzing vector. Then*

$$P : L^p_m(\mathcal{G}) \to \mathcal{V}_\psi(\mathcal{C}o^p_m) \subset L^p_m(\mathcal{G}), \quad F \mapsto F * \mathcal{V}_\psi \psi$$

*is a projection onto the image of the voice transform $\mathcal{V}_\psi : \mathcal{C}o^p_m \to L^p_m(\mathcal{G})$. In particular, we have the equivalence for $F \in L^p_m(\mathcal{G})$*

$$F = \mathcal{V}_\psi f \text{ for some } f \in \mathcal{C}o^p_m \quad \iff \quad F = F * \mathcal{V}_\psi \psi.$$

In other words, the correspondence principle says that $\mathcal{V}_\psi$ defines an isometric isomorphism between the coorbit space $\mathcal{C}o^p_m$ and the *reproducing kernel space*

$$\mathcal{M}^p_m = \{F \in L^p_m(\mathcal{G}) \mid F = F * \mathcal{V}_\psi \psi\}. \tag{3.7}$$

This allows us to work with spaces of continuous functions $\mathcal{M}^p_m$ instead of the abstractly constructed coorbit spaces of distributions $\mathcal{C}o^p_m \subset (\mathcal{H}^1_w)^\sim$ later on.

*Proof.* We have already seen in Corollary 3.3.7 that every function $\mathcal{V}_\psi f$ with $f \in \mathcal{C}o^p_m \subset \mathcal{R}_w$ satisfies the reproducing formula. Thus it remains to show that the opposite direction holds, i.e. that every $F \in L^p_m(\mathcal{G})$ with $F = F * \mathcal{V}_\psi \psi$ can be written as $F = \mathcal{V}_\psi f$ for some $f \in \mathcal{C}o^p_m$.

Suppose $F \in L^p_m(\mathcal{G})$ satisfies $F = F * \mathcal{V}_\psi \psi$. By lemma 3.4.3, the reproducing kernel $\mathcal{V}_\psi \psi$ is contained in $L^q_{1/m}(\mathcal{G})$, since the $p$-control-weight $w$ of $m$ is also a $q$-control-weight of $1/m$ (see lemma 2.2.12). Using this together with Young inequality (2.9), we get

$$\|F\|_{L^\infty_{1/w}} = \|F * \mathcal{V}_\psi \psi\|_{L^\infty_{1/w}} \leq \|F\|_{L^p_m} \|\mathcal{V}_\psi \psi\|_{L^q_{1/m}} < \infty,$$

thus $F \in L^\infty_{1/w}(\mathcal{G})$.

According to the correspondence principle for the reservoir 3.3.8, there exists some $f \in \mathcal{R}_w$ such that $F = \mathcal{V}_\psi f$. But $\mathcal{V}_\psi f = F \in L^p_m(\mathcal{G})$, so it follows immediately that $f \in \mathcal{C}o^p_m$ as desired.

The boundedness of the projection $P$ is a consequence of Young inequality (2.8). $\square$

As a closed subspace of $L^p_m(\mathcal{G})$, the reproducing kernel space $\mathcal{M}^p_m$ is complete. Now $\mathcal{C}o^p_m$ is isometrically isomorphic to $\mathcal{M}^p_m$, so we get the following corollary.

**Corollary 3.4.5.** *The coorbit space $\mathcal{C}o^p_m$ is a Banach space.*

According to the correspondence principle, right sided convolution with the kernel $\mathcal{V}_\psi \psi$ acts as a projection from $L^p_m(\mathcal{G})$ to $\mathcal{M}^p_m$. For $p = 2$ and $m = 1$, this projection is the orthogonal projection as we have seen in theorem 2.6.5. But even for general $p$ and $m$, the projection acts similar to an orthogonal projection, where the inner product is replaced by the antidual pairing of weighted $L^p$ spaces. The next lemma is comparable to the self-adjointness of orthogonal projections.

**Lemma 3.4.6.** *Let $F \in L^p_m(\mathcal{G})$ and $H \in L^q_{1/m}(\mathcal{G})$. Then we have*

$$\int_{\mathcal{G}} (F * \mathcal{V}_\psi \psi)(x) \overline{H(x)} \, dx = \int_{\mathcal{G}} F(x) \overline{(H * \mathcal{V}_\psi \psi)(x)} \, dx.$$

*Proof.* This follows from the short calculation

$$\int_{\mathcal{G}} (F * \mathcal{V}_\psi \psi)(x) \overline{H(x)} \, dx = \int_{\mathcal{G}} \int_{\mathcal{G}} F(y) \mathcal{V}_\psi \psi(y^{-1}x) \, dy \, \overline{H(x)} \, dx$$

$$= \int_{\mathcal{G}} F(y) \int_{\mathcal{G}} \overline{\mathcal{V}_\psi \psi(x^{-1}y) H(x)} \, dx \, dy$$

$$= \int_{\mathcal{G}} F(y) \overline{(H * \mathcal{V}_\psi \psi)(y)} \, dy.$$

$\square$

In particular, the above lemma states that for admissible $\psi$, $F \in \mathcal{M}^p_m$ and $H \in L^q_{1/m}(\mathcal{G})$ we have

$$\int_{\mathcal{G}} F(x) \overline{H(x)} \, dy = \int_{\mathcal{G}} F(x) \overline{(H * \mathcal{V}_\psi \psi)(x)} \, dx$$

since $F = F * \mathcal{V}_\psi \psi$. Of course, a similar statement is true for $\mathcal{M}^q_{1/m}$.

Using the correspondence principle and this 'self-adjointness', we can show that the coorbit spaces $\mathcal{C}o^p_m$ and $\mathcal{C}o^q_{1/m}$ are (anti)dual to each other. Note that the dual norms that arise from this duality are only equivalent to the coorbit space norms, but not equal.

**Proposition 3.4.7.** *[14, Thm. 4.9] Let $1 \leq q < \infty$ and $p^{-1} + q^{-1} = 1$, and suppose $\psi$ is admissible. Then we can identify $\mathcal{C}o^p_m$ with the antidual of $\mathcal{C}o^q_{1/m}$ through the pairing*

$$\langle f, h \rangle_{\mathcal{C}o^p_m \times \mathcal{C}o^q_{1/m}} = \int_{\mathcal{G}} \mathcal{V}_\psi f(x) \overline{\mathcal{V}_\psi h(x)} \, dx$$

*Proof.* We first note that both coorbit spaces $\mathcal{C}o^p_m$ and $\mathcal{C}o^q_{1/m}$ are characterized by the same vector $\psi \in \mathcal{A}_w \setminus \{0\}$, since the weight $w$ is a $p$-control-weight of $m$ as well as a $q$-control-weight of $1/m$. Due to the correspondence principle, it suffices to show that $\mathcal{M}^p_m$ is the antidual space of $\mathcal{M}^q_{1/m}$ with the pairing

$$\langle F, H \rangle_{\mathcal{M}^p_m \times \mathcal{M}^q_{1/m}} = \int_{\mathcal{G}} F(x) \overline{H(x)} \, dx,$$

as this antiduality is pulled back by the isometries $\mathcal{V}_\psi : \mathcal{C}o^p_m \to \mathcal{M}^p_m$ and $\mathcal{V}_\psi : \mathcal{C}o^q_{1/m} \to \mathcal{M}^q_{1/m}$ to the coorbit spaces.

To show that $\mathcal{M}^p_m$ is the antidual of $\mathcal{M}^q_{1/m}$, we prove that the operator

$$T : \mathcal{M}^p_m \to (\mathcal{M}^q_{1/m})^\sim, \quad TF(H) = \int_{\mathcal{G}} F(x) \overline{H(x)} \, dx$$

is an isomorphism between normed spaces, i.e. that it is an isomorphism which induces a norm equivalence.

We begin by showing the injectivity. Suppose that $TF$ is the zero functional for some $F \in \mathcal{M}^p_m$. Then we have for all $H \in L^q_{1/m}(\mathcal{G})$

$$\int_{\mathcal{G}} F(x) \overline{H(x)} \, dx = \int_{\mathcal{G}} F(x) \overline{(H * \mathcal{V}_\psi \psi)(x)} \, dx = T(H * \mathcal{V}_\psi \psi) = 0,$$

where we used lemma 3.4.6. By the duality of $L_m^p(\mathcal{G})$ and $L_{1/m}^q(\mathcal{G})$, it follows $F = 0$.

To show that $T$ is surjective, suppose that $\lambda \in (\mathcal{M}_{1/m}^q)^\sim$ is an arbitrary anti-functional. Then we can extend $\lambda$ to all of $L_{1/m}^q(\mathcal{G})$ by defining

$$\widetilde{\lambda}(H) = \lambda(H * \mathcal{V}_\psi \psi).$$

This $\widetilde{\lambda}$ is bounded because

$$\widetilde{\lambda}(H) \leq \|\lambda\|_{(\mathcal{M}_{1/m}^q)^\sim} \|H * \mathcal{V}_\psi \psi\|_{L_{1/m}^q} \leq \|\lambda\|_{(\mathcal{M}_{1/m}^q)^\sim} \|H\|_{L_{1/m}^q} \|\mathcal{V}_\psi \psi\|_{L_w^1} < \infty$$

for all $H \in L_{1/m}^q(\mathcal{G})$, where we used Young inequality (2.8). Now there exists a Riesz representation $\widetilde{F} \in L_m^p(\mathcal{G})$ of $\widetilde{\lambda}$ such that

$$\widetilde{\lambda}(H) = \int_\mathcal{G} \widetilde{F}(x) \overline{H(x)} \, dx$$

for all $H \in L_{1/m}^q(\mathcal{G})$. If we define $F = \widetilde{F} * \mathcal{V}_\psi \psi \in \mathcal{M}_m^p$, we get for $H \in \mathcal{M}_{1/m}^q$

$$\lambda(H) = \int_\mathcal{G} \widetilde{F}(x) \overline{H(x)} \, dx$$
$$= \int_\mathcal{G} (\widetilde{F} * \mathcal{V}_\psi \psi)(x) \overline{H(x)} \, dx$$
$$= \int_\mathcal{G} F(x) \overline{H(x)} \, dx$$
$$= TF(H)$$

where we used lemma 3.4.6 once more. It follows $TF = \lambda$, so we have proven that $T$ is surjective.

It remains to show that $T$ in fact defines a norm equivalence. For $F \in \mathcal{M}_m^p$ and $H \in \mathcal{M}_{1/m}^q$, we first have

$$|TF(H)| \leq \int_\mathcal{G} |F(x)||H(x)| \, dx \leq \|F\|_{L_m^p} \|H\|_{L_{1/m}^q},$$

thus we obtain the inequality $\|TF\|_{(\mathcal{M}_{1/m}^q)^\sim} \leq \|F\|_{L_m^p}$ (which coincides with the boundedness of $T$). For the other inequality we calculate for $F \in \mathcal{M}_m^p$

$$\|F\|_{L_m^p} = \sup_{H \in L_{1/m}^q(\mathcal{G}) \setminus \{0\}} \frac{1}{\|H\|_{L_{1/m}^q}} \left| \int_\mathcal{G} F(x) \overline{H(x)} \, dx \right|$$
$$= \sup_{\substack{H \in L_{1/m}^q(\mathcal{G}) \\ H * \mathcal{V}_\psi \psi \neq 0}} \frac{1}{\|H\|_{L_{1/m}^q}} \left| \int_\mathcal{G} F(x) \overline{(H * \mathcal{V}_\psi \psi)(x)} \, dx \right|$$
$$= \sup_{\substack{H \in L_{1/m}^q(\mathcal{G}) \\ H * \mathcal{V}_\psi \psi \neq 0}} \frac{\|H * \mathcal{V}_\psi \psi\|_{L_{1/m}^q}}{\|H\|_{L_{1/m}^q}} \frac{|TF(H * \mathcal{V}_\psi \psi)|}{\|H * \mathcal{V}_\psi \psi\|_{L_{1/m}^q}}$$
$$\leq \|\mathcal{V}_\psi \psi\|_{L_w^1} \|TF\|_{(\mathcal{M}_{1/m}^q)^\sim}.$$

In the last step, we used the Young inequality (2.8). This completes the proof. □

Using the property $\mathcal{V}_\psi\psi \in L_m^p(\mathcal{G}) \cap L_{1/m}^q(\mathcal{G})$ from lemma 3.4.3, it is easy to see that there are continuous embeddings $\mathcal{H}_w^1 \hookrightarrow \mathcal{C}o_m^p \hookrightarrow \mathcal{R}_w$. With the above duality, we can also show that these embeddings are dense.

**Proposition 3.4.8.** *We have the continuous and dense embeddings $\mathcal{H}_w^1 \hookrightarrow \mathcal{C}o_m^p \hookrightarrow \mathcal{R}_w$.*

*Proof.* We may assume that $\psi$ is admissible by normalizing it. Then any $g \in \mathcal{H}_w^1$ satisfies the reproducing formula $\mathcal{V}_\psi g = \mathcal{V}_\psi g * \mathcal{V}_\psi\psi$. Using Young inequality (2.7) and the fact that $\mathcal{V}_\psi\psi \in L_m^p(\mathcal{G})$, we get

$$\|\mathcal{V}_\psi g\|_{L_m^p} \leq \|\mathcal{V}_\psi g\|_{L_w^1} \|\mathcal{V}_\psi\psi\|_{L_m^p} < \infty.$$

This shows that $\mathcal{H}_w^1 \subset \mathcal{C}o_m^p$ and that the injection map is bounded by $\|\mathcal{V}_\psi\psi\|_{L_m^p}$.

Similarly, for any $f \in \mathcal{C}o_m^p$ we have by proposition 3.3.5(iv) and Young inequality (2.9)

$$\|f\|_{\mathcal{R}_w} \leq C_\psi^{-2} \|\mathcal{V}_\psi f\|_{L_{1/w}^\infty} \leq C_\psi^{-2} \|\mathcal{V}_\psi f\|_{L_m^p} \|\mathcal{V}_\psi\psi\|_{L_{1/m}^q}.$$

Therefore the injection map $\mathcal{C}o_m^p \hookrightarrow \mathcal{R}_w$ is bounded by $C_\psi^{-2}\|\mathcal{V}_\psi\psi\|_{L_{1/m}^q}$.

We use duality to see that both embeddings are dense. We need the fact that the embeddings $\mathcal{H}_w^1 \hookrightarrow \mathcal{C}o_{1/m}^q \hookrightarrow \mathcal{R}_w$ are continuous, too. But that is clear since control-weights are compatible with duality, so the above arguments are equally true for $\mathcal{C}o_{1/m}^q$.

To show the density of $\mathcal{H}_w^1$ in $\mathcal{C}o_m^p$, we assume that for some $h \in \mathcal{C}o_{1/m}^q$ we have $\langle g, h\rangle_{\mathcal{C}o_m^p \times \mathcal{C}o_{1/m}^q} = 0$ for all $g \in \mathcal{H}_w^1 \subset \mathcal{C}o_m^p$. Then it follows that

$$0 = \langle g, h\rangle_{\mathcal{C}o_m^p \times \mathcal{C}o_{1/m}^q} = \int_\mathcal{G} \mathcal{V}_\psi g(x)\overline{\mathcal{V}_\psi h(x)}\,dx = \overline{\langle h, g\rangle}_{\mathcal{R}_w \times \mathcal{H}_w^1}$$

for all $g \in \mathcal{H}_w^1$, thus $h = 0$ as an element of $\mathcal{R}_w$. This implies that $\mathcal{H}_w^1$ is dense in $\mathcal{C}o_m^p$.

The density of $\mathcal{C}o_m^p \hookrightarrow \mathcal{R}_w$ can be shown similarly. □

*Remark 3.4.9.* [cf. 15, Thm. 4.2(iii)] It is worth mentioning that the coorbit space $\mathcal{C}o_m^p$ is independent from the chosen $p$-control-weight $w$. The smallest possible $p$-control-weight is $w_0$ from (2.6) (see remark 2.2.14 for the proof). By definition, we have $w \geq w_0$ and therefore $L_w^1(\mathcal{G}) \hookrightarrow L_{w_0}^1(\mathcal{G})$. This implies that $\pi$ is $w_0$-integrable as $\{0\} \neq \mathcal{A}_w \subset \mathcal{A}_{w_0}$, and in particular $\psi \in \mathcal{A}_{w_0}$. It also follows that there are continuous embeddings $\mathcal{H}_w^1 \hookrightarrow \mathcal{H}_{w_0}^1$ as well as $\mathcal{R}_{w_0} \hookrightarrow \mathcal{R}_w$. The embeddings are compatible with the antiduality of those spaces, i.e. we have $\langle f, g\rangle_{\mathcal{R}_w \times \mathcal{H}_w^1} = \langle f, g\rangle_{\mathcal{H}_{w_0}^1 \times \mathcal{R}_{w_0}}$ whenever $f \in \mathcal{R}_{w_0}$ and $g \in \mathcal{H}_w^1$.

If now $f \in \mathcal{R}_w$ is is contained in $\mathcal{C}o_m^p$, that is $\|\mathcal{V}_\psi f\|_{L_m^p} < \infty$, we can use Young inequality (2.9) and Lemma 3.4.3 to see that

$$\|\mathcal{V}_\psi f\|_{L_{1/w_0}^\infty} \leq \|\mathcal{V}_\psi f\|_{L_m^p} \|\mathcal{V}_\psi\psi\|_{L_{1/m}^q} < \infty,$$

so $\mathcal{V}_\psi f \in L_{1/w_0}^\infty(\mathcal{G})$. By the correspondence principle for the reservoir 3.3.8 there exists some $f_0 \in \mathcal{R}_{w_0}$ with $\mathcal{V}_\psi f_0 = \mathcal{V}_\psi f$, but since $\mathcal{V}_\psi$ is injective it follows that $f = f_0$ and therefore $f \in \mathcal{R}_{w_0}$.

The coorbit space $\mathcal{C}o_m^p$ is therefore always contained in the smallest distribution space $\mathcal{R}_{w_0}$. The norm $\|f\|_{\mathcal{C}o_m^p} = \|\mathcal{V}_\psi f\|_{L_m^p}$ does not depend $w$. So whatever $p$-control-weight we might use, the construction always yields the same space with the same norm as if we used the minimal weight $w_0$. ◁

# Examples of Voice Transforms and their Coorbit Spaces

In the following two sections, we will express the wavelet transform and the short-time Fourier transform in terms of unitary representations and describe their coorbit theory. This includes the definition of suitable weight functions, as well as finding sufficient conditions for the square-integrability and $w$-integrability of vectors in each case.

We will see that these sufficient conditions are closely related to those properties that are known to be important from the 'classical' theory. Regarding the wavelet transform, we will recover conditions about vanishing moments and smoothness. Regarding the short-time Fourier transform, window functions need to be well localized in the time and frequency domain.

In the following, we use the $d$-dimensional Fourier transform

$$\widehat{f}(\omega) = \mathcal{F}f(\omega) = \int_{\mathbb{R}^d} f(x) e^{-2\pi i x \omega} \, dx$$

as well as the inverse Fourier transform

$$\check{f}(x) = \mathcal{F}^{-1}f(x) = \int_{\mathbb{R}^d} f(\omega) e^{2\pi i x \omega} \, d\omega,$$

where $x\omega$ means the dot product of $x \in \mathbb{R}^d$ and $\omega \in \mathbb{R}^d$.

## 4.1 The Wavelet Transform

We have already encountered the wavelet transform in example 2.4.2, where we identified it as the voice transform of the wavelet representation of the affine group. For completeness, we repeat the necessary definitions and check that all technical prerequisites are met, before moving forward to the associated coorbit theory.

**The Wavelet Representation and its Voice Transform**

**Definition 4.1.1.** We define the *affine group* (also called *ax + b group*) as the topological space $\mathcal{A}\!f\!f = \mathbb{R} \times \mathbb{R}^*$, equipped with the operation

$$(b_1, a_1)(b_2, a_2) = (b_1 + a_1 b_2, a_1 a_2). \tag{4.1}$$

A short calculation shows that (4.1) actually defines a group structure. The neutral element in $\mathcal{A}\!f\!f$ is $(0,1)$. The inverse of $(b,a) \in \mathcal{A}\!f\!f$ is $(-b/a, 1/a)$. Clearly, the group operations are continuous and $\mathcal{A}\!f\!f$ is a Hausdorff space. Thus, $\mathcal{A}\!f\!f$ is a locally compact group.

In addition, the affine group is $\sigma$-compact, as the family

$$\Big([-n,n] \times [2^{-n}, 2^n] \cup [-n,n] \times [-2^n, -2^{-n}]\Big)_{n\in\mathbb{N}}$$

is a countable cover of $\mathcal{A}\!f\!f$. Thus, the affine group satisfies all our technical assumptions.

We have already found a Haar measure on $\mathcal{A}\!f\!f$ in example 2.1.3(iii). We write $d\mu$ for this Haar measure and $da$ resp. $db$ for the Lebesgue measure on $\mathbb{R}$ reps. $\mathbb{R}^*$.

**Proposition 4.1.2.** *A Haar measure on the affine group is given by*

$$d\mu(b,a) = \frac{db\,da}{a^2}.$$

*Its modular function is* $\Delta(b,a) = \Delta(a) = |a|$.

The wavelet representation is defined as follows.

**Definition 4.1.3.** The *wavelet representation* of $\mathcal{A}\!f\!f$ on $L^2(\mathbb{R})$ is given by

$$\pi : \mathcal{A}\!f\!f \to \mathcal{U}(L^2(\mathbb{R})), \quad \pi(b,a)f(x) = \frac{1}{\sqrt{|a|}} f\left(\frac{x-b}{a}\right).$$

Using the dilation and translation operators $D_a, T_b : L^2(\mathbb{R}) \to L^2(\mathbb{R})$,

$$D_a f(x) = \frac{1}{\sqrt{|a|}} f\left(\frac{x}{a}\right), \quad T_b f(x) = f(x-b)$$

for $a \in \mathbb{R}^*$ and $b \in \mathbb{R}$, we can express the wavelet representation as $\pi(b,a) = T_b D_a$.

We first show that $\pi$ defines an unitary representation. That is, we show that $\pi$ defines a strongly continuous group homomorphism and every operator $\pi(b,a)$ is unitary.

Let $f \in L^2(\mathbb{R})$, $x \in \mathbb{R}$ and $(b,a), (b',a') \in \mathcal{G}$. The compatibility of $\pi$ with the group structure follows from the calculation

$$\pi(b,a)\pi(b',a')f(x) = \frac{1}{\sqrt{|a|}} \pi(b',a')f((x-b)/a) = \frac{1}{\sqrt{|aa'|}} f\left(\frac{(x-b)/a - b'}{a'}\right)$$

$$= \frac{1}{\sqrt{|aa'|}} f\left(\frac{x - (ab' + b)}{aa'}\right) = \pi((b,a)(b',a'))f(x).$$

The group operation is unitary since

$$\int_\mathbb{R} |\pi(b,a)f(y)|^2 \, dy = \int_\mathbb{R} \frac{1}{|a|} \left| f\left(\frac{y-b}{a}\right) \right|^2 dy = \int_\mathbb{R} |f(y)|^2 \, dy$$

Thus, $\pi : \mathcal{A}\!f\!f \to \mathcal{U}(L^2(\mathbb{R}))$ is a well-defined group homomorphism.

To see that $\pi$ is strongly continuous, we use the fact that the continuous functions with compact support $C_c(\mathbb{R})$ are dense in $L^2(\mathbb{R})$. For $g \in C_c(\mathbb{R})$ and $a \in \mathbb{R}^*$ we have

$$\|D_y g - D_a g\|_{L^2} \le \|D_y g - D_a g\|_{L^\infty} \|D_y g - D_a g\|_{L^1} \xrightarrow{y \to a} 0$$

since the difference of the $L^\infty$-norm converges to 0 for functions in $C_c(\mathbb{R})$. This property carries over to all functions in $L^2(\mathbb{R})$ by approximating them with sequences $(g_n)_{n\in\mathbb{N}}$ from $C_c(\mathbb{R})$, which shows that the dilation $D$ acts strongly continuous on $L^2(\mathbb{R})$. The same is true for the translation operator $T_b$, $b \in \mathbb{R}$. Hence the concatenation $\pi(b,a) = T_b D_a$ is strongly continuous, too.

The wavelet transform is now the voice transform of the wavelet representation, that is

$$W_g f(b,a) = \mathcal{V}_g f(b,a) = \langle f, \pi(b,a)g \rangle_{L^2(\mathbb{R})} = \frac{1}{\sqrt{|a|}} \int_{\mathbb{R}} f(t) \overline{g\left(\frac{t-b}{a}\right)} dt.$$

We need to show that the wavelet representation is irreducible and square-integrable. According to corollary 2.5.9, the irreducibility is equivalent to the property that all wavelet transforms $W_g : L^2(\mathbb{R}) \to L^\infty(\mathcal{A}\!f\!f)$, $g \in L^2(\mathbb{R})\setminus\{0\}$ are injective. To prove the square-integrability, we explicitly show that $W_g$ is the multiple of an isometry as stated in proposition 2.5.4, which has the advantage of giving us a formula for the constant $C_g$ as well as a simple condition for square-integrability.

We follow the steps of [7, Prop. 2.4.1]. A square-integrable vector of the wavelet representation is called a *wavelet*.

**Proposition 4.1.4.** *The wavelet representation is irreducible and square-integrable. A function $\psi \in L^2(\mathbb{R})$ is a wavelet if and only if*

$$C_\psi = \left( \int_{\mathbb{R}} \frac{|\widehat{\psi}(\omega)|^2}{|\omega|} d\omega \right)^{1/2}$$

*is finite. In that case, $C_\psi$ is the admissibility constant from proposition 2.5.4, that is, $W_\psi : L^2(\mathbb{R}) \to L^2(\mathcal{A}\!f\!f)$ is $C_\psi$ times an isometry.*

*Proof.* We first derive a different way of writing the wavelet transform. Let $f, g \in L^2(\mathbb{R})$ be arbitrary. Using Plancherel's theorem, we can express the value of $W_g f$ at the point $(b,a) \in \mathcal{A}\!f\!f$ as

$$W_g f(b,a) = \langle f, \pi(b,a)g \rangle_{L^2(\mathbb{R})} = \langle \widehat{f}, \widehat{\pi(b,a)g} \rangle_{L^2(\mathbb{R})}.$$

The Fourier transform of $\pi(b,a)g = T_b D_a g$ is equal to

$$\widehat{\pi(b,a)g} = |a|^{1/2} e^{-2\pi i b \cdot} \widehat{g}(a\,\cdot),$$

so we have

$$W_g f(b,a) = |a|^{1/2} \int_{\mathbb{R}} \widehat{f}(t) e^{2\pi i b t} \overline{\widehat{g}(at)}\, dt.$$

The integral can be understood as the inverse Fourier transform of $F_a = |a|^{1/2} \widehat{f}\, \overline{\widehat{g}(a\,\cdot)}$, evaluated at the point $b$. Therefore we have

$$W_g f(b,a) = \check{F}_a(b).$$

We are now able to prove the injectivity of $W_g$ for $g \neq 0$. Suppose $f \in L^2(\mathbb{R})$ is non-zero. Then $\widehat{f}$ and $\widehat{g}$ are non-zero too, which means the sets $A_f = \{\widehat{f} \neq 0\}$ and $A_g = \{\widehat{g} \neq 0\}$ have positive measure. Hence, for $(b,a) \in \mathcal{A}\!f\!f$ with $b \in A_f$ and $ab \in A_g$, it follows $F_a(b) \neq 0$. The set of all those tuples

$$M = \{(b,a) \in \mathbb{R}^2 \mid b \in A_f, ab \in A_g, a \neq 0\} \subset \mathbb{R}^2$$

is no null set, as it has the Lebesgue measure

$$\lambda^2(M) = \int_{A_f} \int_{A_g/b} 1 \, da \, db = \int_{A_f} \lambda^1(A_g) \frac{db}{|b|} > 0.$$

Therefore, the set of all $a \in \mathbb{R}^*$ for which $F_a$ (and therefore $\check{F}_a$, too) does not vanish identically has a positive measure. This implies that

$$\{(b,a) \in \mathcal{A}\!f\!f \mid \check{F}_a \neq 0\} = \{(b,a) \in \mathcal{A}\!f\!f \mid W_g f(b,a) \neq 0\}$$

has positive Lebesgue measure, and by that a positive Haar measure. Since $W_g$ is not the zero function, $W_g$ is injective, which proves the irreducibility of $\pi$.

Next we will show that $W_\psi$ is $C_\psi$ times an isometry, provided $C_\psi$ is finite. We again write $F_a = |a|^{1/2} \widehat{f} \, \widehat{\psi}(a \cdot)$ and obtain with Plancherel's theorem

$$\|W_\psi f\|_{L^2}^2 = \int_\mathbb{R} \int_\mathbb{R} |\check{F}_a(b)|^2 \, db \, \frac{da}{a^2} = \int_\mathbb{R} \int_\mathbb{R} |F_a(b)|^2 \, db \, \frac{da}{a^2}.$$

By Fubini's theorem, the last integral evaluates to

$$\|W_\psi f\|_{L^2}^2 = \int_\mathbb{R} \int_\mathbb{R} |a| |\widehat{f}(b)|^2 |\widehat{\psi}(ab)|^2 \, db \, \frac{da}{a^2}$$
$$= \int_\mathbb{R} |\widehat{f}(b)|^2 \int_\mathbb{R} \frac{|\widehat{\psi}(ba)|^2}{|a|} \, da \, db$$
$$= \int_\mathbb{R} |\widehat{f}(b)|^2 \int_\mathbb{R} \frac{|\widehat{\psi}(a)|^2}{|a|} \, da \, db$$
$$= \|f\|_{L^2} \int_\mathbb{R} \frac{|\widehat{\psi}(a)|^2}{|a|} \, da.$$

The last integral converges if and only if $C_\psi$ is finite. By setting $f = \psi$ we see that $\psi$ is square-integrable if and only if $C_\psi < \infty$. In that case, the above calculation shows that $W_\psi$ is $C_\psi$ times an isometry. □

We have collected all necessary information about the wavelet representation and wavelet transform to move on to the associated coorbit theory. Before doing so, we note that the Duflo-Moore operator of the wavelet representation is given by

$$D : \mathcal{Q} \subset L^2(\mathbb{R}) \to L^2(\mathbb{R}), \quad D\psi = \mathcal{F}^{-1}\left[\frac{\mathcal{F}\psi}{\sqrt{|\cdot|}}\right].$$

This can be shown similar to proposition 4.1.4, together with the uniqueness statement from theorem 2.5.7.

**Wavelet Coorbit Theory**

We first have to define the weight functions on the affine group we are going to use. We use the same weights as [5, sec. 3.2.3.1].

**Definition 4.1.5.** We define the following weights on $\mathcal{A}\!f\!f$:

a) For $s \in \mathbb{R}$, we set $m_s(b,a) = |a|^{-s}$.

b) For $\rho \geq 0$, we set $w_\rho(b, a) = |a|^\rho + |a|^{-\rho}$.

Both functions $m_s$ and $w_\rho$ are obviously locally integrable and therefore indeed weight functions. They only depend on the scale factor $a$ but not on the shift parameter $b$. The function $m_s$ weights for small $s$ the coarse scales $|a| \to \infty$ stronger, and for large $s$ the fine scales $|a| \to 0$.

**Lemma 4.1.6.** *The weights on $\mathcal{A}\!f\!f$ have the following properties:*

(i) *$m_s$ is moderate with respect to $\alpha = \beta = m_s$ for all $s \in \mathbb{R}$.*

(ii) *$w_\rho$ is submultiplicative for all $\rho \geq 0$.*

(iii) *For $\rho \geq |s| + \max\{1/p, 1/q\}$, $w_\rho$ is a p-control-weight of $m_s$. In particular, $w_\rho$ is a simultaneous control-weight of $m_s$ if $\rho \geq |s| + 1$.*

*Proof.*

(i) For $(b_1, a_1), (b_2, a_2) \in \mathcal{A}\!f\!f$ we have

$$\begin{aligned}
m_s((b_1, a_1)(b_2, a_2)) &= m_s(b_1 + a_1 b_2, a_1 a_2) \\
&= |a_1 a_2|^{-s} \\
&= |a_1|^{-s} |a_2|^{-s} \\
&= m_s(b_1, a_1) m_s(b_2, a_2),
\end{aligned}$$

thus $m_s$ is multiplicative, and in particular moderate.

(ii) For $(b_1, a_1), (b_2, a_2) \in \mathcal{A}\!f\!f$ we have the inequality

$$\begin{aligned}
w_\rho((b_1, a_1)(b_2, a_2)) &= w_\rho(b_1 + a_1 b_2, a_1 a_2) \\
&= |a_1 a_2|^\rho + |a_1 a_2|^{-\rho} \\
&\leq (|a_1|^\rho + |a_1|^{-\rho})(|a_2|^\rho + |a_2|^{-\rho}) \\
&= w_\rho(b_1, a_1) w_\rho(b_2, a_2),
\end{aligned}$$

thus $w_\rho$ is submultiplicative.

(iii) A short calculation shows that the inequality

$$w_\rho(b, a) = |a|^\rho + |a|^{-\rho} \geq \max\{|a|^{-s}, |a|^s, |a|^{-s-1/q}, |a|^{s-1/p}\}$$

holds whenever $\rho \geq |s| + \max\{1/p, 1/q\}$. By (i), we have $|a|^{-s} \geq \alpha_0$ and $|a|^{-s} \geq \beta_0$, and we have $\Delta(b, a) = |a|$, so this inequality implies that the condition (2.6) is indeed satisfied by $w_\rho$.

□

We want to construct the wavelet coorbit spaces $\mathcal{C}o^p_{m_s}$ with respect to the p-control-weight $w_\rho$. For the theory of chapter 3 to be applicable, we need to show that the wavelet representation is $w_\rho$-integrable for $\rho \geq 0$, that is, that the set of analyzing vectors

$$\mathcal{A}_{w_\rho} = \{\psi \in L^2(\mathbb{R}) \mid W_\psi \psi \in L^1_{w_\rho}(\mathcal{A}\!f\!f)\}$$

contains some non-trivial wavelet.

For that we derive a sufficient condition for a function $\psi$ to be $w_\rho$-integrable, which depends on the vanishing moments and the smoothness of $\psi$. We need to prove some auxiliary lemmas beforehand. We use a similar approach as [13].

We first show that a function with $L$ vanishing moments has an antiderivative with $L-1$ vanishing moments.

**Lemma 4.1.7.** *Suppose that $f \in L^1(\mathbb{R})$ is continuous and has $L \in \mathbb{N}$ vanishing moments, that is*
$$\int_{\mathbb{R}} t^k f(t)\, dt = 0, \quad \text{for } k = 0, \ldots, L-1,$$
*with absolutely convergent integrals. Suppose further that the L-th moment is absolutely convergent, i.e. $\int_{\mathbb{R}} |t|^L |f(t)|\, dt$ is finite.*

*Then the antiderivative $h(x) = \int_{-\infty}^{x} f(t)\, dt$ of $f$ is well-defined and has $L-1$ absolutely convergent vanishing moments, as well as an absolutely convergent $L-1$-st moment.*

*Proof.* First of all, the antiderivative $h$ is well-defined since $f$ is integrable.

We begin by proving that the first $L$ moments of $h$ are absolutely convergent. As a consequence of the 0th vanishing moment of $f$, its antiderivative $h$ can be written as
$$h(x) = \int_{-\infty}^{x} f(t)\, dt = -\int_{x}^{\infty} f(t)\, dt$$
for every $x \in \mathbb{R}$, which leads to the estimate
$$|h(x)| \leq \left| \int_{|t| \geq |x|} f(t)\, dt \right| \leq \int_{|t| \geq |x|} |f(t)|\, dt.$$

Using this inequality we see that
$$\int_{\mathbb{R}} |x|^k |h(x)|\, dx \leq \int_{\mathbb{R}} \int_{|t| \geq |x|} |x|^k |f(t)|\, dt\, dx$$
$$= \int_{\mathbb{R}} \int_{|x| \leq |t|} |x|^k |f(t)|\, dx\, dt$$
$$= \frac{2}{k+1} \int_{\mathbb{R}} |t|^{k+1} |f(t)|\, dt.$$

The last integral is finite for $k = 0, \ldots, L-1$, thus the first $L$ moments of $h$ are absolutely convergent.

For $k = 0, \ldots, L-2$ we use partial integration and the vanishing moments of $f$ to obtain
$$\int_{\mathbb{R}} x^k h(x)\, dx = -\int_{\mathbb{R}} \frac{1}{k+1} x^{k+1} f(x)\, dx + \frac{1}{k+1} \lim_{y \to \infty} \left( y^{k+1} h(y) - (-y)^{k+1} h(-y) \right)$$
$$= \frac{1}{k+1} \lim_{y \to \infty} \left( y^{k+1} \int_{-\infty}^{y} f(t)\, dt - (-y)^{k+1} \int_{-\infty}^{-y} f(t)\, dt \right)$$
$$= \frac{1}{k+1} \lim_{y \to \infty} \left( -y^{k+1} \int_{y}^{\infty} f(t)\, dt - (-y)^{k+1} \int_{-\infty}^{-y} f(t)\, dt \right).$$

Taking the absolute value yields

$$\left|\int_{\mathbb{R}} x^k h(x)\, dx\right| \leq \frac{1}{k+1} \lim_{y\to\infty} \left(\int_y^\infty |y|^{k+1}|f(t)|\, dt + \int_{-\infty}^{-y} |y|^{k+1}|f(t)|\, dt\right)$$

$$\leq \frac{1}{k+1} \lim_{y\to\infty} \left(\int_y^\infty |t|^{k+1}|f(t)|\, dt + \int_{-\infty}^{-y} |t|^{k+1}|f(t)|\, dt\right).$$

Since the first $L + 1$ moments of $f$ are absolutely convergent, the right side converges to 0, proving that the first $L - 1$ moments of $h$ are vanishing. $\square$

The 0th moment of the function $h$ is just the integral $\int_{\mathbb{R}} h(t)\, dt$. The absolute convergence of the 0th moment of $h$ therefore implies that $h$ is an integrable function.

By induction we get the following corollary.

**Corollary 4.1.8.** *Suppose that $f \in L^1(\mathbb{R})$ is continuous and has $L \in \mathbb{N}$ absolutely convergent and vanishing moments as well as an absolutely convergent L-th moment. Then there exists an L-times differentiable function $h \in L^1(\mathbb{R})$ such that $f = h^{(L)}$ and all derivatives $h^{(k)}$, $k = 0, \ldots, L$ are integrable functions.*

*Additionally, the derivatives of $h$ satisfy the recursive rule*

$$h^{(k)}(x) = \int_{-\infty}^x h^{(k+1)}(t)\, dt, \quad k = 0, \ldots, L-1.$$

We also need a lemma about partial integration. It is an one-dimensional version of [13, Lemma 3.1].

**Lemma 4.1.9.** *Let $f, g : \mathbb{R} \to \mathbb{C}$ be differentiable functions such that $f, g, f', g' \in L^1(\mathbb{R})$. Then we have*

$$\int_{\mathbb{R}} f(t) g'(t+x)\, dt = -\int_{\mathbb{R}} f'(t) g(t+x)\, dt$$

*for almost all $x \in \mathbb{R}$.*

*Proof.* It follows from Fubini's theorem that

$$\int_{\mathbb{R}} \int_{\mathbb{R}} |f(t)||g'(t+x)|\, dt\, dx = \|f\|_{L^1} \|g'\|_{L^1} < \infty.$$

This implies that the inner integral $\int_{\mathbb{R}} f(t) g'(t+x)\, dt$ is absolutely convergent for almost all $x \in \mathbb{R}$. The same is true for the integral $\int_{\mathbb{R}} f'(t) g(t+x)\, dt$. If we show that $f(y)g(y+x) \to 0$ for $y \to \pm\infty$ and almost all $x \in \mathbb{R}$, the stated integral equality therefore follows from the usual partial integration formula.

The integral $\int_{\mathbb{R}} f'(t)\, dt$ is absolutely convergent by requirement. Thus, the limits

$$\lim_{y\to\infty} f(y) = f(0) + \int_0^\infty f'(t)\, dt$$

and

$$\lim_{y\to-\infty} f(y) = f(0) - \int_{-\infty}^0 f'(t)\, dt$$

do exist. But since the integral $\int_{\mathbb{R}} |f(t)|\, dt$ is finite as well, these limits have to be 0. The same is true for all functions $t \mapsto g(t+x)$, $x \in \mathbb{R}$. This finishes the proof. $\square$

We are now able to prove our sufficient condition for the $w_\rho$-integrability of a wavelet $\psi$. We do that by adapting the proof of [13, Lemma 3.3].

**Proposition 4.1.10.** *Suppose the function $\psi \in L^2(\mathbb{R})$ has $L \in \mathbb{N}$ vanishing moments as well as $L+1$ absolutely convergent moments. Suppose further that $\psi$ is $L$-times continuously differentiable such that all derivatives $\psi^{(k)}$, $k = 0, \ldots, L$ are integrable.*

*Then $\psi \in \mathcal{A}_{w_\rho}$ for $0 \leq \rho < L - 1/2$.*

*Proof.* We define the first $L$ antiderivatives of $\psi$ inductively by $\psi^{(0)} = \psi$,

$$\psi^{(-k)}(x) = \int_{-\infty}^{x} \psi^{(-k+1)}(t)\, dt, \quad k = 1, \ldots, L.$$

Since $\psi$ satisfies all requirements of corollary 4.1.8, all those antiderivatives are well-defined and integrable. Overall, we have the family $\psi^{(-L)}, \ldots, \psi^{(L)}$ of integrable functions.

We now consider the wavelet transform $W_\psi \psi$ at the point $(b, a) \in \mathcal{A\!f\!f}$. It is defined as

$$W_\psi \psi(b, a) = |a|^{-1/2} \int_{\mathbb{R}} \psi(t) \overline{\psi\left(\frac{t-b}{a}\right)}\, dt = |a|^{1/2} \int_{\mathbb{R}} \psi(at+b) \overline{\psi(t)}\, dt.$$

We apply the partial integration formula from lemma 4.1.9 to the second integral $L$ times. By differentiating the left factor $L$ times, we get

$$W_\psi \psi(b, a) = (-1)^L |a|^{1/2} \int_{\mathbb{R}} a^L \psi^{(L)}(at+b) \overline{\psi^{(-L)}(t)}\, dt.$$

If we instead differentiate the right factor $L$ times, we get

$$W_\psi \psi(b, a) = (-1)^L |a|^{1/2} \int_{\mathbb{R}} a^{-L} \psi^{(-L)}(at+b) \overline{\psi^{(L)}(t)}\, dt.$$

We use these formulas to estimate the integral of $|W_\psi \psi(b, a)|$ over $b$. The first formula yields

$$\int_{\mathbb{R}} |W_\psi \psi(b, a)|\, db \leq |a|^{1/2+L} \int_{\mathbb{R}} \int_{\mathbb{R}} |\psi^{(L)}(at+b)|\, |\psi^{(-L)}(t)|\, dt\, db$$
$$= |a|^{1/2+L} \left\|\psi^{(L)}\right\|_{L^1(\mathbb{R})} \left\|\psi^{(-L)}\right\|_{L^1(\mathbb{R})}, \tag{4.2}$$

while the second formula leads to the estimate

$$\int_{\mathbb{R}} |W_\psi \psi(b, a)|\, db \leq |a|^{1/2-L} \left\|\psi^{(-L)}\right\|_{L^1(\mathbb{R})} \left\|\psi^{(L)}\right\|_{L^1(\mathbb{R})}. \tag{4.3}$$

To estimate the $L^1_{w_\rho}(\mathcal{A\!f\!f})$-norm of $W_\psi \psi$, we have to integrate these terms with respect to $(|a|^\rho + |a|^{-\rho})\frac{da}{a^2}$. For this, we split the integral at $|a| = 1$. For $|a| \leq 1$ we apply inequality (4.2), for $|a| > 1$ we apply inequality (4.3), and obtain

$$\int_{\mathbb{R}} \int_{\mathbb{R}} |W_\psi \psi(b, a)|\, db\, (|a|^\rho + |a|^{-\rho})\frac{da}{a^2}$$
$$\leq \int_{|a| \leq 1} |a|^{1/2+L} \left\|\psi^{(L)}\right\|_{L^1(\mathbb{R})} \left\|\psi^{(-L)}\right\|_{L^1(\mathbb{R})} (|a|^\rho + |a|^{-\rho})\frac{da}{a^2}$$
$$+ \int_{|a| > 1} |a|^{1/2-L} \left\|\psi^{(-L)}\right\|_{L^1(\mathbb{R})} \left\|\psi^{(L)}\right\|_{L^1(\mathbb{R})} (|a|^\rho + |a|^{-\rho})\frac{da}{a^2}$$
$$= \left\|\psi^{(-L)}\right\|_{L^1(\mathbb{R})} \left\|\psi^{(L)}\right\|_{L^1(\mathbb{R})}$$
$$\left( \int_{|a| \leq 1} |a|^{\rho+L-3/2} + |a|^{-\rho+L-3/2}\, da + \int_{|a| > 1} |a|^{\rho-L-3/2} + |a|^{-\rho-L-3/2}\, da \right).$$

The integral over $|a| \leq 1$ converges for $-\rho + L - 3/2 > -1$, which is equivalent to $\rho < L - 1/2$. The integral over $|a| > 1$ converges for $\rho - L - 3/2 < -1$, which in turn is equivalent to $\rho < L + 1/2$.

Thus the wavelet transform $W_\psi \psi$ is contained in $L^1_{w_\rho}(\mathcal{A}\!f\!f)$ for $0 \leq \rho < L - 1/2$, making $\psi$ a $w_\rho$-integrable vector for those values of $\rho$. □

***Remark 4.1.11.*** In the preceding proof, the inequalities (4.2) and (4.3) are crucial. Both have the form
$$\int_{\mathbb{R}} |W_\psi \psi(b, a)|\, db \leq C|a|^\nu$$
for a suitable exponent $\nu \in \mathbb{R}$ and some constant $C > 0$. These estimates are related to the fact that the wavelet transform $W_\psi f(b,a)$ of a function $f \in L^2(\mathbb{R})$ is contained in the class $\mathcal{O}(|a|^\nu)$ for $|a| \to 0$, if $\psi$ has enough vanishing moments and $f$ is smooth enough [27, Thm. 6.3/6.4].

Indeed, a proof similar to that of [27, Thm. 6.4] can be used to obtain the needed inequalities. There, $\psi$ is developed into its $L$-th Taylor polynomial centred at $b$. After that, integration over $db$ yields the expected estimate. ◁

We can now use the above condition 4.1.10 to prove that the wavelet representation is $w_\rho$-integrable for all $\rho \geq 0$. The space of *Schwartz functions* is given by
$$\mathcal{S} = \left\{ f : \mathbb{R} \to \mathbb{C} \,\Big|\, f \in C^\infty(\mathbb{R}),\, \sup_{x \in \mathbb{R}} |x|^k |f^{(l)}(x)| < \infty \text{ for all } k, l \in \mathbb{N}_0 \right\}.$$

Every function from $\mathcal{S}$ satisfies the smoothness and integrability requirements of 4.1.10. Thus, any Schwartz function with $L > \rho + 1/2$ vanishing moments is $w_\rho$-integrable. This condition is equivalent to $\widehat{\psi}$ being a Schwartz function that has a zero of multiplicity at least $L$ in $\omega = 0$.

Such a function is, for instance, given by
$$\widehat{\psi}(x) = e^{-(x^2 + x^{-2})},$$
as $\widehat{\psi}$ has a zero of multiplicity $\infty$ in 0. Other examples are certain Meyer wavelets [27, sec. 7.2.2]. Thus, the space of Schwartz functions with infinitely many vanishing moments
$$\mathcal{S}_0 = \left\{ f \in \mathcal{S} \,\Big|\, \int_{\mathbb{R}} t^k f(t)\, dt = 0 \text{ for all } k \in \mathbb{N}_0 \right\}$$
is non-trivial.

**Corollary 4.1.12.** *For every $\rho \geq 0$, we have $\mathcal{S}_0 \subset \mathcal{A}_{w_\rho}$. In particular, the wavelet representation is $w_\rho$-integrable for all $\rho \geq 0$.*

We can now apply the theory of chapter 3 to the wavelet transform. For a given $w_\rho$-integrable wavelet $\psi \in \mathcal{A}_{w_\rho} \setminus \{0\}$, the space of test vectors $\mathcal{H}^1_{w_\rho}$ consists of all functions $g \in L^2(\mathbb{R})$ for which the integral
$$\int_{\mathbb{R}} \int_{\mathbb{R}} |W_\psi g(b,a)| \left(|a|^\rho + |a|^{-\rho}\right) db\, \frac{da}{a^2}$$

is finite. Based on this, the reservoir $\mathcal{R}_{w_\rho} = (\mathcal{H}^1_{w_\rho})^\sim$ is defined. If now $s \in \mathbb{R}$ such that that $\rho \geq |s| + 1$, then $w_\rho$ is a simultaneous control-weight of $m_s(b, a) = |a|^{-s}$ according to lemma 4.1.6. The wavelet coorbit spaces are then defined as

$$\mathcal{C}o^p_{m_s} = \left\{ f \in \mathcal{R}_{w_\rho} \;\Big|\; \int_\mathbb{R} \int_\mathbb{R} |W_\psi f(b, a)|^p |a|^{-sp} \, db \, \frac{da}{a^2} < \infty \right\}$$

for $1 \leq p \leq \infty$.

In chapter 3, we constructed the space of test vectors and the reservoir in order to get a suitable distribution space for general $w$-integrable representations. For explicit representations, there might be other distribution spaces, described by another dual pair, that work just as fine as the reservoir. In this case, we can use the space $\mathcal{S}_0$ as a space of test functions and its dual $\mathcal{S}_0^\sim$ as a space of distributions.

First of all, the vector space $\mathcal{S}_0$ is invariant under the action of $\pi$, since translations and dilations do not affect vanishing moments or the Schwartz property. According to corollary 3.2.5, every $\psi \in \mathcal{S}_0 \setminus \{0\} \subset \mathcal{A}_{w_\rho}$ is cyclic in $\mathcal{H}^1_{w_\rho}$. Thus $\mathcal{S}_0$ is dense in $\mathcal{H}^1_{w_\rho}$.

The Schwartz space $\mathcal{S}$ is equipped with the topology that is induced by the family of norms

$$|f|_N = \max_{k,l \in \{0,\ldots,N\}} \sup_{x \in \mathbb{R}} |x|^k |f^{(l)}(x)|, \quad N \in \mathbb{N}_0.$$

By restricting this topology to $\mathcal{S}_0$, the space becomes a topological vector space. A calculation similar to proposition 4.1.4 shows that $\mathcal{S}_0$ is continuously embedded into $\mathcal{H}^1_{w_\rho}$.

Therefore, we have the continuous and dense embedding $\mathcal{S}_0 \hookrightarrow \mathcal{H}^1_{w_\rho}$, which implies by duality that $\mathcal{R}_{w_\rho}$ is continuously and densely embedded into $\mathcal{S}_0^\sim$. The space $\mathcal{S}_0^\sim$ can be understood as the quotient space $\mathcal{S}^\sim / \mathcal{P}$, where $\mathcal{S}^\sim$ is the space of tempered distributions and $\mathcal{P}$ denotes its subspace of all polynomials, since by definition the polynomials are exactly those distributions that act trivially on $\mathcal{S}_0$.

The wavelet coorbit spaces are therefore continuously embedded into the distribution space $\mathcal{S}_0^\sim = \mathcal{S}^\sim / \mathcal{P}$, and can even be constructed as subspaces of those distributions in the first place.

We now want to discuss the connection between the wavelet transform and the homogeneous Besov spaces. Latter are defined as follows.

**Definition 4.1.13.** [30, Sec. 5.1.3] Let $(\varphi_j)_{j \in \mathbb{Z}} \subset \mathcal{S}$ be a family of Schwartz functions with the following properties:

i. For all $j \in \mathbb{Z}$ we have $\mathrm{supp}(\varphi_j) \subset \{x \in \mathbb{R} \mid 2^{j-1} \leq |x| \leq 2^{j+1}\}$.

ii. For $x \in \mathbb{R}^*$ we have $\sum_{j \in \mathbb{Z}} \varphi_j(x) = 1$.

iii. For all $n \in \mathbb{N}_0$ there is a $C_n > 0$ such that $2^{jn} |\varphi_j^{(n)}| \leq C_n$ holds for all $j \in \mathbb{Z}$.

Furthermore, let $p, q \in [1, \infty]$ and $\sigma \in \mathbb{R}$.

Then the *homogeneous Besov space* $\dot{B}^\sigma_{p,q}$ is defined as the space of distributions $f \in \mathcal{S}_0^\sim$ for which the norm

$$\|f\|_{\dot{B}^\sigma_{p,q}} = \left( \sum_{j \in \mathbb{Z}} 2^{j\sigma q} \left\| \mathcal{F}^{-1}[\varphi_j \mathcal{F} f] \right\|_{L^p}^q \right)^{1/q}$$

is finite (with the usual adaptations for $q = \infty$).

The homogeneous Besov spaces are Banach spaces which describe the smoothness of distributions $f$ and which were used even before the wavelet transform became well-known. In the definition of the Besov norm, $f$ represents an entire equivalence class in $\mathcal{S}_0^{\sim} = \mathcal{S}^{\sim}/\mathcal{P}$.

For any distribution $f \in \mathcal{S}'$, each product $\varphi_j \mathcal{F} f$, $j \in \mathbb{Z}$ has compact support, thus $\mathcal{F}^{-1}[\varphi_j \mathcal{F} f]$ is indeed an (analytical) function of which the $L^p$-norm is well-defined. To prove that the overall norm $\|\cdot\|_{\dot{B}^{\sigma}_{p,q}}$ is well-defined, it has to be shown that $\|f + p\|_{\dot{B}^{\sigma}_{p,q}} = \|f\|_{\dot{B}^{\sigma}_{p,q}}$ is true for all polynomials $p$. This can be done by using some of the common computation rules of the Fourier transform. For more information about homogeneous Besov spaces, we refer to [30, Ch. 5].

The homogeneous Besov spaces can now be identified with the wavelet coorbit spaces. In order to make the relation more clear, we consider the norm

$$\|F\|_{L^{p,q}_s(\mathcal{A}\!f\!f)} = \left[\int_{\mathbb{R}} \left(\int_{\mathbb{R}} |F(b,a)|^p |a|^{-sp}\, db\right)^{q/p} \frac{da}{a^2}\right]^{1/q}$$

for functions $F : \mathcal{A}\!f\!f \to \mathbb{C}$, where the $L^p$-norm is taken with respect to $b$ and the $L^q$-norm is taken with respect to $a$. The function $F$ is additionally multiplied with the weight $|a|^{-s}$. For $p = \infty$ or $q = \infty$, the usual adjustments have to be made. The spaces $L^{p,q}_s(\mathcal{A}\!f\!f)$ that are defined through this norm are a generalization of the weighted $L^p$-spaces.

We now evaluate the $L^{p,q}_s(\mathcal{A}\!f\!f)$-norm of some wavelet transform $W_\psi f$. For that, we write $W_\psi f$ as in the proof of proposition 4.1.4 as

$$W_\psi(b,a) = |a|^{1/2} \int_{\mathbb{R}} \hat{f}(\omega)\overline{\hat{\psi}(a\omega)}e^{2\pi i b \omega}\, d\omega = |a|^{1/2} \mathcal{F}^{-1}\left[\overline{\hat{\psi}(a\cdot)}\mathcal{F}f\right](b).$$

This yields

$$\|W_\psi f\|_{L^{p,q}_s(\mathcal{A}\!f\!f)} = \left[\int_{\mathbb{R}^*} \left(\int_{\mathbb{R}} |W_\psi f(b,a)|^p |a|^{-sp}\, db\right)^{q/p} \frac{da}{a^2}\right]^{1/q}$$

$$= \left[\int_{\mathbb{R}^*} \left(\int_{\mathbb{R}} |a|^{p/2} \left|\mathcal{F}^{-1}\left[\overline{\hat{\psi}(a\cdot)}\mathcal{F}f\right](b)\right|^p |a|^{-sp}\, db\right)^{q/p} \frac{da}{a^2}\right]^{1/q}$$

$$= \left[\int_{\mathbb{R}^*} |a|^{q/2 - sq - 1} \left\|\mathcal{F}^{-1}\left[\overline{\hat{\psi}(a\cdot)}\mathcal{F}f\right]\right\|_{L^p}^q \frac{da}{|a|}\right]^{1/q}.$$

By setting $\sigma = s - 1/2 + 1/q$, this can be written as

$$\left[\int_{\mathbb{R}^*} \left\|\mathcal{F}^{-1}\left[\overline{\hat{\psi}(a\cdot)}\mathcal{F}f\right]\right\|_{L^p}^q |a|^{-\sigma q} \frac{da}{|a|}\right]^{1/q}.$$

The norm of $f$ in the space $\dot{B}^{\sigma}_{p,q}$ is essentially the discretization of the integral over $\frac{da}{|a|}$ by summing over the points $a = a_j = 2^{-j}$ for $j \in \mathbb{Z}$. The function $\overline{\hat{\psi}(2^{-j}\cdot)}$ corresponds to the function $\varphi_j$. That this discretization really leads to an equivalent norm (for suitable $\psi$) was proven by Triebel [29, Cor. 10] (independent from the wavelet theory).

For $p = q$, the above norm amounts to the $L^p_{m_s}(\mathcal{A}\!f\!f)$-norm. The wavelet coorbit space $\mathcal{C}o^p_{m_s}$ is a subspace of $\mathcal{S}_0^{\sim}$, just as the homogeneous Besov space $\dot{B}^{\sigma}_{p,q}$. Since the norms of both spaces are equivalent, we get the following theorem.

**Theorem 4.1.14.** *[cf. 14, Sec. 7.2] The coorbit space $\mathcal{C}o_{m_s}^p$ coincides with the homogeneous Besov space $\dot{B}_{p,p}^{s-1/2+1/p}$ and their norms are equivalent.*

It is also possible to define the more general coorbit spaces

$$\mathcal{C}o_{m_s}^{p,q} = \{f \in \mathcal{S}_0^\sim \mid W_\psi f \in L_s^{p,q}(\mathcal{A}\!f\!f)\}$$

by extending the necessary theory to the mixed norm spaces $L_s^{p,q}(\mathcal{A}\!f\!f)$. Then we have $\mathcal{C}o_{m_s}^{p,q} \cong \dot{B}_{p,q}^{s-1/2+1/q}$.

**Remark 4.1.15.** In this section, we described the wavelet transform in terms of the affine group $\mathcal{A}\!f\!f = \mathbb{R} \rtimes \mathbb{R}^*$ and its wavelet representation. It is possible to generalize this example through the quasi-regular representations from example 2.4.2. This leads to multidimensional wavelet coorbit theory, where integrability conditions are associated to a more abstract concept of vanishing moments [18; 19]. One prominent example of this multidimensional theory is the shearlet transform [5]. ◁

## 4.2 The Short-Time Fourier Transform

Our next example is the short-time Fourier transform

$$V_g f(x, \omega) = \int_\mathbb{R} f(t) \overline{g(t-x)} e^{-2\pi i t \omega} \, dt.$$

We want to describe this transform as the voice transform of some unitary representation of a locally compact group to apply our theory from the two preceding chapters. However, we have already seen in example 2.5.3 that there exists no unitary representation which is able to achieve that. We will nonetheless find a voice transform that is closely related to the short-time Fourier transform, circumventing the problem with only minor technical difficulties.

### The Heisenberg Group and its Schrödinger Representation

To find a suitable voice transform we will use the reduced Heisenberg group, which extends the time-frequency plane $\mathbb{R}^d \times \mathbb{R}^d$ by phase factors from $S^1$. There are different ways to describe this group; we are following the steps of [22, ch. 9].

**Definition 4.2.1.** For $d \in \mathbb{N}$ we define the *(d-dimensional) reduced Heisenberg group* as the topological space $\mathbb{H}_r^d = \mathbb{R}^d \times \mathbb{R}^d \times S^1$, equipped with the group structure

$$(x_1, \omega_1, \tau_1)(x_2, \omega_2, \tau_2) = \left(x_1 + x_2, \omega_1 + \omega_2, \tau_1 \tau_2 e^{\pi i(x_2\omega_1 - x_1\omega_2)}\right).$$

By $x_2\omega_1$ and $x_1\omega_2$ we mean the $d$-dimensional dot product.

The reduced Heisenberg group is obviously locally compact and Hausdorff. Since the group structure is continuous, it is a locally compact group. It is covered by the countable family of compact sets $[-n,n]^d \times [-n,n]^d \times S^1$, $n \in \mathbb{N}$, and is therefore $\sigma$-compact.

**Remark 4.2.2.** The *full Heisenberg group* is defined similarly as the space $\mathbb{H}^d = \mathbb{R}^d \times \mathbb{R}^d \times \mathbb{R}$ with the group structure

$$(x_1, \omega_1, t_1)(x_2, \omega_2, t_2) = \left(x_1 + x_2, \omega_1 + \omega_2, t_1 + t_2 + \frac{1}{2}(x_2\omega_1 - x_1\omega_2)\right).$$

The reduced Heisenberg group can then be written as the quotient

$$\mathbb{H}_r^d = \mathbb{H}^d \big/ \{0\} \times \{0\} \times \mathbb{Z}.$$

Since we need the compactness of the third factor $S^1$, we will only work with the reduced Heisenberg group. ◁

A Haar measure on the reduced Heisenberg group is induced by the usual Lebesgue measure. We write $dx$ and $d\omega$ for the Lebesgue measure on $\mathbb{R}^d$, and $d\tau$ for the rotational invariant normalized measure on $S^1$.

**Proposition 4.2.3.** *The product measure $d\mu(x, \omega, \tau) = dx\, d\omega\, d\tau$ defines a Haar measure on the reduced Heisenberg group $\mathbb{H}_r^d$, and $\mathbb{H}_r^d$ is unimodular.*

*Proof.* Suppose $F : \mathbb{H}_r^d \to \mathbb{C}$ is integrable with respect to the product measure $dx\, d\omega\, d\tau$. Then we have for fixed $(x_0, \omega_0, \tau_0) \in \mathbb{H}_r^d$

$$\int_{S^1} \int_{\mathbb{R}^d} \int_{\mathbb{R}^d} F((x_0, \omega_0, \tau_0)(x, \omega, \tau))\, dx\, d\omega\, d\tau$$
$$= \int_{S^1} \int_{\mathbb{R}^d} \int_{\mathbb{R}^d} F\left(x + x_0, \omega + \omega_0, \tau \tau_0 e^{\pi i(x\omega_0 - x_0\omega)}\right) dx\, d\omega\, d\tau.$$

The integral with respect to $x$ and $\omega$ is translation invariant, and the integral with respect to $\tau$ is rotation invariant. We can also change the order of integration as we like since $F$ is assumed to be integrable. Thus we get

$$\int_{S^1} \int_{\mathbb{R}^d} \int_{\mathbb{R}^d} F\left((x_0, \omega_0, \tau_0)(x, \omega, \tau)\right) dx\, d\omega\, d\tau$$
$$= \int_{S^1} \int_{\mathbb{R}^d} \int_{\mathbb{R}^d} F\left(x, \omega, \tau\tau_0 e^{\pi i(x\omega_0 - x_0\omega)}\right) dx\, d\omega\, d\tau$$
$$= \int_{\mathbb{R}^d} \int_{\mathbb{R}^d} \int_{S^1} F\left(x, \omega, \tau\tau_0 e^{\pi i(x\omega_0 - x_0\omega)}\right) d\tau\, dx\, d\omega$$
$$= \int_{\mathbb{R}^d} \int_{\mathbb{R}^d} \int_{S^1} F(x, \omega, \tau) d\tau\, dx\, d\omega$$
$$= \int_{S^1} \int_{\mathbb{R}^d} \int_{\mathbb{R}^d} F(x, \omega, \tau)\, dx\, d\omega\, d\tau.$$

The integration with respect to $d\mu(x, \omega, \tau) = dx\, d\omega\, d\tau$ is therefore invariant under the left action of $\mathbb{H}_r^d$, proving that $\mu$ is a Haar measure.

It can be shown the same way that $\mu$ is also a right Haar measure, implying that $\mathbb{H}_r^d$ is unimodular. □

The Schrödinger representation now acts on $L^2(\mathbb{R}^d)$ through time-frequency-shifts by $(x, \omega)$, extended by the multiplication with the phase factor $\tau e^{-\pi i \omega x}$.

**Definition 4.2.4.** The *Schrödinger representation* of $\mathbb{H}_r^d$ on $L^2(\mathbb{R})$ is defined as

$$\rho(x, \omega, \tau) : L^2(\mathbb{R}^d) \to L^2(\mathbb{R}^d), \quad \rho(x, \omega, \tau)f(t) = \tau e^{-\pi i \omega x} e^{2\pi i t \omega} f(t - x).$$

Using the modulation and translation operators $M_\omega$ and $T_x$, $\rho$ can be written as

$$\rho(x, \omega, \tau) = \tau e^{-\pi i \omega x} M_\omega T_x f(t).$$

We first show that $\rho$ indeed defines a unitary representation on $L^2(\mathbb{R}^d)$. The compatibility of $\rho$ with the group structure follows from a simple computation. Using the commutation relation $M_\omega T_x = e^{2\pi i x \omega} T_x M_\omega$, we have for $(x_1, \omega_1, \tau_1), (x_2, \omega_2, \tau_2) \in \mathbb{H}_r^d$

$$\rho((x_1,\omega_1,\tau_1)(x_2,\omega_2,\tau_2)) = \rho\left(x_1+x_2, \omega_1+\omega_2, \tau_1\tau_2 e^{\pi i(x_2\omega_1 - x_1\omega_2)}\right)$$
$$= \tau_1\tau_2 e^{\pi i(x_2\omega_1 - x_1\omega_2)} e^{-\pi i(x_1+x_2)(\omega_1+\omega_2)} M_{\omega_1+\omega_2} T_{x_1+x_2}$$
$$= \tau_1\tau_2 e^{-\pi i(2x_1\omega_2 + x_1\omega_1 + x_2\omega_2)} M_{\omega_1} M_{\omega_2} T_{x_1} T_{x_2}$$
$$= \tau_1\tau_2 e^{-\pi i(2x_1\omega_2 + x_1\omega_1 + x_2\omega_2)} e^{2\pi i x_1 \omega_2} M_{\omega_1} T_{x_1} M_{\omega_2} T_{x_2}$$
$$= \rho(x_1,\omega_1,\tau_1)\rho(x_2,\omega_2,\tau_2).$$

The short calculation

$$\|\rho(x,\omega,\tau)f\|_2^2 = \int_{\mathbb{R}^d} \left|\tau e^{\pi i x \omega} e^{2\pi i t \omega} f(t-x)\right|^2 dt$$
$$= \int_{\mathbb{R}^d} |f(t-x)|^2 dt = \|f\|_2^2$$

for $f \in L^2(\mathbb{R}^d)$ and $(x,\omega,\tau) \in \mathbb{H}_r^d$ shows that $\rho$ is unitary. The strong continuity of $\rho$ follows from the fact that $\rho$ is the concatenation of the continuous mappings $x \mapsto T_x f$, $\omega \mapsto M_\omega f$ and $(x,\omega,\tau) \mapsto \tau e^{-\pi i \omega x} f$.

Next we show that the Schrödinger representation is irreducible and square-integrable. As in the case of the wavelet transform, we do this by explicitly deriving the orthogonality relations of the associated voice transform. For $f, g \in L^2(\mathbb{R}^d)$, this transform is given by

$$\mathcal{V}_g f(x,\omega,\tau) = \int_{\mathbb{R}^d} f(t) \overline{\tau e^{-\pi i \omega x} M_\omega T_x f(t)}\, dt = \overline{\tau} e^{\pi i \omega x} V_g f(x,\omega).$$

Since $\mathcal{V}_g f$ and $V_g f$ only differ by a factor of absolute value 1, the orthogonality relations of the voice transform can be proven in the same way as for the short-time Fourier transform [22, Thm. 3.2.1].

**Proposition 4.2.5.** *For $g_1, g_2, f_1, f_2 \in L^2(\mathbb{R}^d)$ we have*

$$\langle \mathcal{V}_{g_1} f_1, \mathcal{V}_{g_2} f_2 \rangle_{L^2(\mathbb{H}_r^d)} = \overline{\langle g_1, g_2 \rangle}_{L^2(\mathbb{R}^d)} \langle f_1, f_2 \rangle_{L^2(\mathbb{R}^d)}. \tag{4.4}$$

*Proof.* We first assume that $\mathcal{V}_{g_1} f_1$ and $\mathcal{V}_{g_2} f_2$ are contained in $L^2(\mathbb{H}_r^d)$. Then their inner product is

$$\langle \mathcal{V}_{g_1} f_1, \mathcal{V}_{g_2} f_2 \rangle_{L^2(\mathbb{H}_r^d)} = \int_{\mathbb{R}^d} \int_{\mathbb{R}^d} \int_{S^1} \mathcal{V}_{g_1} f_1(x,\omega,\tau) \overline{\mathcal{V}_{g_2} f_2(x,\omega,\tau)}\, d\tau\, d\omega\, dx$$
$$= \int_{\mathbb{R}^d} \int_{\mathbb{R}^d} \int_{S^1} \left( \int_{\mathbb{R}^d} f_1(t) \overline{\tau} e^{\pi i x \omega} e^{-2\pi i t \omega} \overline{g_1(t-x)}\, dt \right)$$
$$\cdot \left( \int_{\mathbb{R}^d} \overline{f_2(t)} \tau e^{-\pi i x \omega} e^{2\pi i t \omega} g_2(t-x)\, dt \right) d\tau\, d\omega\, dx.$$

The factors $\overline{\tau} e^{\pi i x \omega}$ and $\tau e^{-\pi i x \omega}$ are independent of $t$ and cancel each other. The remaining integral with respect to $\tau$ evaluates to $|S^1| = 1$ since the integrand is no longer dependent on $\tau$. What remains is

$$\langle \mathcal{V}_{g_1} f_1, \mathcal{V}_{g_2} f_2 \rangle_{L^2(\mathbb{H}_r^d)}$$
$$= \int_{\mathbb{R}^d} \int_{\mathbb{R}^d} \left( \int_{\mathbb{R}^d} f_1(t) \overline{g_1(t-x)} e^{-2\pi i t \omega}\, dt \right) \overline{\left( \int_{\mathbb{R}^d} f_2(t) \overline{g_2(t-x)} e^{-2\pi i t \omega}\, dt \right)} dx\, d\omega.$$

The integrals with respect to $t$ can be understood as the Fourier transforms of $P_x = f_1 \overline{T_x g_1}$ reps. $Q_x = f_1 \overline{T_x g_2}$, evaluated at the point $\omega$, hence

$$\langle \mathcal{V}_{g_1} f_1, \mathcal{V}_{g_2} f_2 \rangle_{L^2(\mathbb{H}_r^d)} = \int_{\mathbb{R}^d} \int_{\mathbb{R}^d} \widehat{P_x}(\omega) \overline{\widehat{Q_x}(\omega)} \, dx \, d\omega = \int_{\mathbb{R}^d} \int_{\mathbb{R}^d} \widehat{P_x}(\omega) \overline{\widehat{Q_x}(\omega)} \, d\omega \, dx.$$

Using Plancherel's theorem it follows

$$\begin{aligned}\langle \mathcal{V}_{g_1} f_1, \mathcal{V}_{g_2} f_2 \rangle_{L^2(\mathbb{H}_r^d)} &= \int_{\mathbb{R}^d} \int_{\mathbb{R}^d} P_x(u) \overline{Q_x(u)} \, du \, dx \\ &= \int_{\mathbb{R}^d} \int_{\mathbb{R}^d} f_1(u) \overline{g_1(u-x)} \, \overline{f_2(u)} g_2(u-x) \, du \, dx \\ &= \overline{\langle g_1, g_2 \rangle}_{L^2(\mathbb{R}^d)} \langle f_1, f_2 \rangle_{L^2(\mathbb{R}^d)},\end{aligned}$$

proving equation (4.4).

We have shown in particular that $\|\mathcal{V}_g f\|_{L^2(\mathbb{H}_r^d)} = \|f\|_{L^2(\mathbb{R}^d)} \|g\|_{L^2(\mathbb{R}^d)}$, provided that $\mathcal{V}_g f \in L^2(\mathbb{H}_r^d)$. If we now assume $f_1, f_2, g_1, g_2 \in L^2(\mathbb{R}^d)$ to be arbitrary, we can follow the same steps of computation backwards to see that indeed $\mathcal{V}_{g_1} f_1, \mathcal{V}_{g_2} f_2 \in L^2(\mathbb{R}^d)$. This finishes the proof. □

The above proposition implies that all $g \in L^2(\mathbb{R}^d)$ are square-integrable with respect to the Schrödinger representation. Comparing our result to theorem 2.5.7, we see that the Duflo-Moore operator of the Schrödinger representation is the identity on $L^2(\mathbb{R}^d)$. The irreducibility of $\rho$ is now easy to show.

**Corollary 4.2.6.** *The Schrödinger representation is irreducible and square-integrable, with all functions in $L^2(\mathbb{R}^d)$ being square-integrable.*

*Proof.* According to corollary 2.5.9, the Schrödinger representation is irreducible if and only if the voice transform $\mathcal{V}_g$ is injective for all $g \in L^2(\mathbb{R}^d) \setminus \{0\}$. Since all those voice transforms are positive multiples of an isometry, this is obviously the case. □

**Coorbit Theory for the Short-Time Fourier Transform**

We now dedicate our attention to the coorbit theory of the Schrödinger representation. We begin by defining the relevant weight functions. We use the same weights as [20].

**Definition 4.2.7.** For $r, s \geq 0$ we define the weight function $v_{r,s}$ on $\mathbb{H}_r^d$ by

$$v_{r,s}(x, \omega, \tau) = v_{r,s}(x, \omega) = (1 + |x|)^r (1 + |\omega|)^s.$$

We also write for simplicity $v_r(x) = (1+|x|)^r$ and $v_s(\omega) = (1+|\omega|)^s$, so we have $v_{r,s}(x, \omega) = v_r(x) v_s(\omega)$ for all $x, \omega \in \mathbb{R}^d$.

Since the weight $v_{r,s}$ is independent of the phase factor $\tau$, we omit it most of the time. Moreover, we have

$$v_{r,s}((x_1, \omega_1, \tau_1)(x_2, \omega_2, \tau_2)) = v_{r,s}(x_1 + x_2, \omega_1 + \omega_2)$$

for all $(x_1, \omega_1, \tau_1), (x_2, \omega_2, \tau_2) \in \mathbb{H}_r^d$. The weight $v_{r,s}$ can therefore be understood as a weight on the additive group $\mathbb{R}^d \times \mathbb{R}^d$. As such, it defines the weighted $L^p$ space

$$L^p_{v_{r,s}}(\mathbb{R}^{2d}) = \left\{ F : \mathbb{R}^d \times \mathbb{R}^d \to \mathbb{C} \; \middle| \; \int_{\mathbb{R}^d} \int_{\mathbb{R}^d} |F(x, \omega)|^p v_{r,s}(x, \omega)^p \, dx \, d\omega < \infty \right\}.$$

Similarly, the weight function $v_r$ defines the space

$$L^p_{v_r}(\mathbb{R}^d) = \left\{ f : \mathbb{R}^d \to \mathbb{C} \,\middle|\, \int_{\mathbb{R}^d} |f(t)|^p v_r(t)^p \, dt < \infty \right\}.$$

The weights $v_{r,s}$, $r, s \geq 0$ are submultiplicative and moderate. The unimodularity of the reduced Heisenberg group even enables them to be their own control-weights.

**Lemma 4.2.8.** *Let $r, s \geq 0$.*

(i) *The weight $v_{r,s}$ is submultiplicative and moderate with respect to itself.*

(ii) *For $0 \leq r_0 \leq r$ and $0 \leq s_0 \leq s$, $v_{r,s}$ is a control-weight for $v_{r_0, s_0}$.*

*Proof.*

(i) For $x_1, x_2 \in \mathbb{R}^d$ we have

$$\begin{aligned} v_r(x_1 + x_2) &= (1 + |x_1 + x_2|)^r \leq (1 + |x_1| + |x_2| + |x_1 x_2|)^r \\ &= (1 + |x_1|)^r (1 + |x_2|)^r = v_r(x_1) v_r(x_2). \end{aligned}$$

Similarly $v_s(\omega_1 + \omega_2) \leq v_s(\omega_1) v_s(\omega_2)$ holds for $\omega_1, \omega_2 \in \mathbb{R}^d$. This implies the submultiplicativity and the moderateness of $v_{r,s}$.

(ii) With $m = v_{r_0, s_0}$ and $w = v_{r,s}$, inequality (2.6) holds if

$$v_{r,s}(x, \omega) \geq \max\{v_{r_0, s_0}(x, \omega), v_{r_0, s_0}(-x, -\omega)\}$$

is satisfied, as $\mathbb{H}^d_r$ is unimodular and the inverse of $(x, \omega, \tau)$ is $(-x, -\omega, \tau')$ for some $t' \in S^1$. Since $v_{r,s}$ is symmetric and $v_{r,s} \geq v_{r_0, s_0}$, this inequality indeed holds for all $x, \omega \in \mathbb{R}^d$.

$\square$

The sets and spaces used in coorbit theory are mainly defined through integrability conditions of the voice transform $\mathcal{V}_g f$ of certain functions $f, g \in L^2(\mathbb{R}^d)$. The fact that $|\mathcal{V}_g f| = |V_g f|$, combined with the compactness of the third factor $S^1$ in the reduced Heisenberg group, allows us to replace the voice transform in these integrability conditions by the short-time Fourier transform.

**Lemma 4.2.9.** *Suppose $f, g \in L^2(\mathbb{R}^d)$ and $1 \leq p \leq \infty$. Then $\mathcal{V}_g f \in L^p_{v_{r,s}}(\mathbb{H}^d_r)$ if and only if $V_g f \in L^p_{v_{r,s}}(\mathbb{R}^{2d})$ and the respective norms coincide.*

*Proof.* For $1 \leq p < \infty$ we have

$$\begin{aligned} \|\mathcal{V}_g f\|_{L^p_{v_{r,s}}(\mathbb{H}^d_r)} &= \int_{\mathbb{R}^d} \int_{\mathbb{R}^d} \int_{S^1} |\mathcal{V}_g f(x, \omega, \tau)|^p v_{r,s}(x, \omega, \tau)^p \, d\tau \, dx \, d\omega \\ &= \int_{\mathbb{R}^d} \int_{\mathbb{R}^d} \int_{S^1} |\overline{\tau} e^{\pi i x \omega}|^p |V_g f(x, \omega)|^p v_{r,s}(x, \omega)^p \, d\tau \, dx \, d\omega \\ &= \int_{\mathbb{R}^d} \int_{\mathbb{R}^d} |V_g f(x, \omega)|^p v_{r,s}(x, \omega)^p \, dx \, d\omega \\ &= \|V_g f\|_{L^p_{v_{r,s}}(\mathbb{R}^d \times \mathbb{R}^d)}, \end{aligned}$$

where we used the normalization $|S^1| = 1$. The case $p = \infty$ works similarly. $\square$

In the following, we will always express the integrability conditions in terms of the short-time Fourier transform.

The set of analyzing vectors of the Schrödinger representation is now given by

$$\mathcal{A}_{v_{r,s}} = \{g \in L^2(\mathbb{R}^d) \mid V_g g \in L^1_{v_{r,s}}(\mathbb{R}^d \times \mathbb{R}^d)\}.$$

Again, we want to prove a sufficient condition for a function $g$ to be an analyzing vector. For that, we use a simplified version of [20, Thm. 1.1]. It says essentially that every function $g$ with good enough time and frequency localization is contained in $\mathcal{A}_{v_{r,s}}$.

**Theorem 4.2.10.** *Let $r, s \geq 0$, $a > 2r + d$ and $\beta > 2s + d$. Suppose that $g \in L^2(\mathbb{R}^d)$ satisfies*

$$\|g\|_{L^1_{v_\alpha}(\mathbb{R}^d)} = \int_{\mathbb{R}^d} |g(x)|(1 + |x|)^\alpha \, dx < \infty$$

*as well as*

$$\|\widehat{g}\|_{L^1_{v_\beta}(\mathbb{R}^d)} = \int_{\mathbb{R}^d} |\widehat{g}(\omega)|(1 + |\omega|)^\beta \, d\omega < \infty.$$

*Then $g \in \mathcal{A}_{v_{r,s}}$.*

**Remark 4.2.11.** The original result of [20] is more general. It states that, if the function $f \in L^2(\mathbb{R}^d)$ satisfies

$$\|f\|_{L^p_{v_r}} < \infty \quad \text{and} \quad \|\widehat{f}\|_{L^q_{v_s}} < \infty$$

($p$ and $q$ independent), it is contained in a mixed norm coorbit space $\mathcal{C}o^{p',q'}_{v_{r',s'}}$ of the short-time Fourier transform, where the parameters $p, q, r, s$ and $p', q', r', s'$ are related by an inequality that has to hold. If $p' = q'$, these mixed norm coorbit spaces coincide with our coorbit spaces for weighted $L^p$-spaces. $\triangleleft$

To make the proof more accessible, we first derive some required inequalities. The first is similar to Hölder's inequality and helps us to work with the double integral $\int_{\mathbb{R}^d} \int_{\mathbb{R}^d} dx \, d\omega$.

**Lemma 4.2.12.** *[cf. 20, Lemma 2.1] For measurable functions $F, H : \mathbb{R}^d \times \mathbb{R}^d \to \mathbb{C}$ we have the inequality*

$$\int_{\mathbb{R}^d} \int_{\mathbb{R}^d} |F(x,\omega)||H(x,\omega)| \, dx \, d\omega$$

$$\leq \left( \int_{\mathbb{R}^d} \operatorname*{ess\,sup}_{\omega \in \mathbb{R}^d} |F(x,\omega)| \, dx \right) \left( \int_{\mathbb{R}^d} \operatorname*{ess\,sup}_{x \in \mathbb{R}^d} |H(x,\omega)| \, d\omega \right). \quad (4.5)$$

*Proof.* We apply Hölder's inequality regarding the exponents 1 and $\infty$ first to the inner, then to the outer integral. This gives us

$$\int_{\mathbb{R}^d} \int_{\mathbb{R}^d} |F(x,\omega)||H(x,\omega)| \, dx \, d\omega \leq \int_{\mathbb{R}^d} \left( \int_{\mathbb{R}^d} |F(x,\omega)| \, dx \right) \left( \operatorname*{ess\,sup}_{x \in \mathbb{R}^d} |H(x,\omega)| \right) d\omega$$

$$\leq \left( \operatorname*{ess\,sup}_{\omega \in \mathbb{R}^d} \int_{\mathbb{R}^d} |F(x,\omega)| \, dx \right) \left( \int_{\mathbb{R}^d} \operatorname*{ess\,sup}_{x \in \mathbb{R}^d} |H(x,\omega)| \, d\omega \right)$$

$$\leq \left( \int_{\mathbb{R}^d} \operatorname*{ess\,sup}_{\omega \in \mathbb{R}^d} |F(x,\omega)| \, dx \right) \left( \int_{\mathbb{R}^d} \operatorname*{ess\,sup}_{x \in \mathbb{R}^d} |H(x,\omega)| \, d\omega \right)$$

as claimed. □

Next we prove two estimates regarding the short-time Fourier transform.

**Lemma 4.2.13.** *Let $f, g \in L^2(\mathbb{R}^d)$ and $\alpha, \beta \geq 0$. Then we have*

$$\int_{\mathbb{R}^d} \sup_{\omega \in \mathbb{R}^d} |V_g f(x,\omega)| v_\alpha(x) \, dx \leq \|f\|_{L^1_{v_\alpha}(\mathbb{R}^d)} \|g\|_{L^1_{v_\alpha}(\mathbb{R}^d)} \tag{4.6}$$

*as well as*

$$\int_{\mathbb{R}^d} \sup_{x \in \mathbb{R}^d} |V_g f(x,\omega)| v_\beta(\omega) \, d\omega \leq \|\widehat{f}\|_{L^1_{v_\beta}(\mathbb{R}^d)} \|\widehat{g}\|_{L^1_{v_\beta}(\mathbb{R}^d)}. \tag{4.7}$$

*Proof.* For the first inequality, we substitute in the definition of the short-time Fourier transform and get

$$\int_{\mathbb{R}^d} \sup_{\omega \in \mathbb{R}^d} |V_g f(x,\omega)| v_\alpha(x) \, dx = \int_{\mathbb{R}^d} \sup_{\omega \in \mathbb{R}^d} \left| \int_{\mathbb{R}^d} f(t) \overline{g(t-x)} e^{-2\pi i t \omega} \, dt \right| v_\alpha(x) \, dx$$

$$\leq \int_{\mathbb{R}^d} \int_{\mathbb{R}^d} |f(t)| \, |g(t-x)| \, dt \, v_\alpha(x) \, dx.$$

Since we have $v_\alpha(x) = v_\alpha(t + (x-t)) \leq v_\alpha(t) v_\alpha(x-t)$ and $v_\alpha(x-t) = v_\alpha(t-x)$, it follows

$$\int_{\mathbb{R}^d} \sup_{\omega \in \mathbb{R}^d} |V_g f(x,\omega)| v_\alpha(x) \, dx \leq \int_{\mathbb{R}^d} \int_{\mathbb{R}^d} |f(t)| v_\alpha(t) \, |g(t-x)| v_\alpha(t-x) \, dt \, dx$$

$$= \|f\|_{L^1_{v_\alpha}(\mathbb{R}^d)} \|g\|_{L^1_{v_\alpha}(\mathbb{R}^d)},$$

proving the first estimate.

To show the second inequality, we rewrite the short-time Fourier transform as an integral over $\widehat{f}$ and $\widehat{g}$. By applying Plancherel's theorem we obtain

$$V_g f(x,\omega) = \int_{\mathbb{R}^d} f(t) \overline{g(t-x)} e^{-2\pi i x \omega} \, dt = \int_{\mathbb{R}^d} \widehat{f}(\xi) \overline{\widehat{g}(\xi-\omega)} e^{-2\pi i (\xi-\omega) x} \, d\xi.$$

Now the same computation as above yields

$$\int_{\mathbb{R}^d} \sup_{x \in \mathbb{R}^d} |V_g f(x,\omega)| v_\beta(\omega) \, d\omega \leq \int_{\mathbb{R}^d} \int_{\mathbb{R}^d} |\widehat{f}(\xi)| \, |\widehat{g}(\xi-\omega)| v_\beta(\omega) \, d\xi \, d\omega$$

$$\leq \|\widehat{f}\|_{L^1_{v_\beta}(\mathbb{R}^d)} \|\widehat{g}\|_{L^1_{v_\beta}(\mathbb{R}^d)},$$

so the second inequality follows. □

For the last auxiliary lemma we divide the time frequency plane $\mathbb{R}^d \times \mathbb{R}^d$ into the two disjoint domains

$$A_\sigma = \{(x,\omega) \in \mathbb{R}^d \times \mathbb{R}^d \mid |\omega| < |x|^\sigma\}$$

and

$$B_\sigma = \{(x,\omega) \in \mathbb{R}^d \times \mathbb{R}^d \mid |\omega| \geq |x|^\sigma\}$$

with respect to some $\sigma \in (0, \infty)$.

**Lemma 4.2.14.** *[cf. 20, Lemma 2.3] Let $\sigma > 0$.*

(i) For $r < 0$, $s > 0$ and $\sigma < -r/(s+d)$, the integral

$$\int_{\mathbb{R}^d} \sup_{x \in \mathbb{R}^d} v_r(x) v_s(\omega) \chi_{A_\sigma}(x, \omega) \, d\omega$$

is finite.

(ii) For $r > 0$, $s < 0$ and $\sigma > -(r+d)/s$, the integral

$$\int_{\mathbb{R}^d} \sup_{\omega \in \mathbb{R}^d} v_r(x) v_s(\omega) \chi_{B_\sigma}(x, \omega) \, dx$$

is finite.

*Proof.*

(i) First we can write

$$\int_{\mathbb{R}^d} \sup_{x \in \mathbb{R}^d} v_r(x) v_s(\omega) \chi_{A_\sigma}(x, \omega) \, d\omega = \int_{\mathbb{R}^d} \sup_{|x| > |\omega|^{1/\sigma}} v_r(x) v_s(\omega) \, d\omega$$

$$= \int_{\mathbb{R}^d} \sup_{|x| > |\omega|^{1/\sigma}} (1 + |x|)^r (1 + |\omega|)^s \, d\omega.$$

Since $r < 0$, the function $(1 + |x|)^r$ decreases with growing $|x|$, thus the supremum of that term is exactly $(1 + |\omega|^{1/\sigma})^r$. Therefore we have

$$\int_{\mathbb{R}^d} \sup_{x \in \mathbb{R}^d} v_r(x) v_s(\omega) \chi_{A_\sigma}(x, \omega) \, d\omega = \int_{\mathbb{R}^d} (1 + |\omega|^{1/\sigma})^r (1 + |\omega|)^s \, d\omega.$$

For $|\omega| \to \infty$, the factors in the integral are asymptotically bounded by $C|\omega|^{r/\sigma}$ resp. $C|\omega|^s$ for some constant $C > 0$. The integral over $\mathbb{R}^d$ now converges if the sum of the exponents is smaller then $-d$, that is when $r/\sigma + s < -d$. This is equivalent to the assumed inequality $\sigma < -r/(s+d)$, so the integral is indeed finite.

(ii) Similar to the first integral, we have

$$\int_{\mathbb{R}^d} \sup_{\omega \in \mathbb{R}^d} v_r(x) v_s(\omega) \chi_{B_\sigma}(x, \omega) \, dx = \int_{\mathbb{R}^d} \sup_{|\omega| \geq |x|^\sigma} (1 + |x|)^r (1 + |\omega|)^s \, dx$$

$$= \int_{\mathbb{R}^d} (1 + |x|)^r (1 + |x|^\sigma)^s \, dx.$$

This integral converges if the exponent sum $r + \sigma s$ is smaller than $-d$, which again is equivalent to $\sigma > -(r+d)/s$.

□

We are now able to prove theorem 4.2.10 using the above lemmas.

*Proof of theorem 4.2.10.* We have to show that the integral

$$\int_{\mathbb{R}^d} \int_{\mathbb{R}^d} |V_g g(x, \omega)| v_r(x) v_s(\omega) \, dx \, d\omega$$

is finite if $g$ satisfies the required integrability conditions. In order to do that, we divide the space $\mathbb{R}^d \times \mathbb{R}^d$ into the two domains $A_\sigma$ and $B_\sigma$ with respect to some $\sigma \in (0, \infty)$ that is specified later on. We estimate the integral over these domains separately.

For the integral over $A_\sigma$, we can apply lemma 4.2.12 as well as the multiplicativity $v_r = v_{r-\alpha} v_\alpha$ to obtain

$$\int_{A_\sigma} |V_g g(x,\omega)| v_r(x) v_s(\omega) \, d(x,\omega)$$
$$= \int_{\mathbb{R}^d} \int_{\mathbb{R}^d} \left( |V_g g(x,\omega)| v_\alpha(x) \right) \left( v_{r-\alpha}(x) v_s(\omega) \chi_{A_\sigma}(x,\omega) \right) dx \, d\omega$$
$$\leq \left( \int_{\mathbb{R}^d} \sup_{\omega \in \mathbb{R}^d} |V_g g(x,\omega)| v_\alpha(x) \, dx \right) \left( \int_{\mathbb{R}^d} \sup_{x \in \mathbb{R}^d} v_{r-\alpha}(x) v_s(\omega) \chi_{A_\sigma}(x,\omega) \, d\omega \right).$$

The first factor is at most $\|g\|^2_{L^1_{v_\alpha}(\mathbb{R}^d)}$, as we have proven in lemma 4.2.13. If we assume that $\sigma < -(r-\alpha)(s+d)$, lemma 4.2.14 states that the second factor is finite, too.

We similarly proceed with the integral over $B_\sigma$. Here an application of lemma 4.2.12 yields

$$\int_{B_\sigma} |V_g g(x,\omega)| v_r(x) v_s(\omega) \, dx \, d\omega$$
$$= \int_{\mathbb{R}^d} \int_{\mathbb{R}^d} \left( |V_g g(x,\omega)| v_\beta(\omega) \right) \left( v_r(x) v_{s-\beta}(\omega) \chi_{B_\sigma}(x,\omega) \right) dx \, d\omega$$
$$\leq \left( \int_{\mathbb{R}^d} \sup_{x \in \mathbb{R}^d} |V_g g(x,\omega)| v_\beta(\omega) \, d\omega \right) \left( \int_{\mathbb{R}^d} \sup_{\omega \in \mathbb{R}^d} v_r(x) v_{s-\beta}(\omega) \chi_{B_\sigma}(x,\omega) \, dx \right).$$

The first factor is at most $\|\hat{g}\|^2_{L^1_{v_\beta}(\mathbb{R}^d)}$ by lemma 4.2.13, the second finite by lemma 4.2.14 if the inequality $\sigma > -(r+d)/(s-\beta)$ holds.

Combining the two preceding estimates, it follows that $g$ is $v_{r,s}$-integrable if we can find some $\sigma$ that satisfies the two inequalities

$$-\frac{r+d}{s-\beta} < \sigma < -\frac{r-\alpha}{s+d}.$$

The existence of such $\sigma$ is equivalent to

$$-\frac{r+d}{s-\beta} < -\frac{r-\alpha}{s+d}$$
$$\iff (r+d)(s+d) < (\alpha - r)(\beta - s).$$

Since by requirement $\alpha > 2r+d$ and $\beta > 2s+d$, this condition is indeed satisfied, finishing the proof. □

The localization in time and frequency that is required in our sufficient condition 4.2.10 is satisfied, for instance, by every Schwartz function for any $\alpha, \beta > 0$. Thus, we immediately get the following corollary.

**Corollary 4.2.15.** *The Schwartz functions $\mathcal{S}$ are contained in $\mathcal{A}_{v_{r,s}}$ for all $r, s \geq 0$. In particular, the Schrödinger representation is $v_{r,s}$-integrable for every $r, s \geq 0$.*

We therefore can carry out the constructions from coorbit theory. For fixed $r, s \geq 0$, the space of test vectors is given by

$$M^1_{r,s} = \mathcal{H}^1_{v_{r,s}} = \{f \in L^2(\mathbb{R}^d) \mid V_g f \in L^1_{v_{r,s}}(\mathbb{R}^d \times \mathbb{R}^d) \text{ for some } g \in \mathcal{A}_{v_{r,s}} \setminus \{0\}\}.$$

In the unweighted case $r = s = 0$, this space is called *Feichtinger's algebra* and is noted $S_0$.

We have seen in proposition 3.2.6 that the set of analyzing vectors and test vectors coincide whenever the submultiplicative weight $w$ satisfies the symmetry relation $w = \Delta^{-1} w^\vee$. The reduced Heisenberg group is unimodular and the weight $v_{r,s} = v_{r,s}^\vee$ is symmetric, so it follows $\mathcal{A}_{v_{r,s}} = M^1_{r,s}$. Functions from $M^1_{r,s}$ therefore are not only test vectors to work with distributions from the reservoir, but also suitable window functions for the voice transform resp. the short-time Fourier transform.

In the following construction, the coorbit spaces are defined as subspaces of the reservoir $\mathcal{R}_{v_{r,s}} = (M^1_{r,s})^\sim$. Once again we want to replace the voice transform by the short-time Fourier transform, meaning that we have to extend the short-time Fourier transform to $\mathcal{R}_{v_{r,s}}$. This is done exactly as for the voice transform in section 3.3 by setting

$$V_g f(x, \omega) = \langle f, M_\omega T_x g \rangle_{\mathcal{R}_{v_{r,s}} \times M^1_{r,s}}, \quad (x, \omega) \in \mathbb{R}^d \times \mathbb{R}^d,$$

but we need to know beforehand that $M^1_{r,s}$ is invariant under time-frequency shifts $M_\omega T_x$. However, this is easy to see as the space $M^1_{r,s}$ is invariant under the Schrödinger representation $\rho$, and time frequency-shifts can be written as

$$M_\omega T_x = \rho(x, \omega, e^{\pi i x \omega}).$$

Now that the extended short-time Fourier transform is established, it can be proven similar to lemma 4.2.9 that

$$\|V_g f\|_{L^p_{v_{r,s}}(\mathbb{R}^d \times \mathbb{R}^d)} = \|\mathcal{V}_g f\|_{L^p_{v_{r,s}}(\mathbb{H}^d_r)}$$

for all $g \in M^1_{r,s}$ and $f \in \mathcal{R}_{v_{r,s}}$.

According to lemma 4.2.8, $v_{r,s}$ is its own control-weight. The coorbit spaces with respect to the Schrödinger representation can therefore be defined as

$$M^p_{r,s} = \mathcal{C}o^p_{v_{r,s}} = \{f \in \mathcal{R}_{v_{r,s}} \mid V_g f \in L^p_{v_{r,s}}(\mathbb{R}^d \times \mathbb{R}^d)\}$$

with $1 \leq p \leq \infty$ and some fixed $g \in M^1_{r,s}$. These coorbit spaces are called *modulation spaces*. They are of major importance in time-frequency analysis since their norm describes the time-frequency behaviour of (generalized) functions.

Like the coorbit spaces of the wavelet transform, the modulation spaces can be constructed as subspaces of another suitable distribution space than the reservoir, namely the space of tempered distributions. The space of Schwartz functions $\mathcal{S}$ is invariant under time frequency shifts resp. under the action of $\rho$. Since it is contained in $M^1_{r,s}$, it is already dense in that space according to corollary 3.2.5. It can also be shown that the embedding $\mathcal{S} \hookrightarrow M^1_{r,s}$ is continuous with a computation similar to 4.2.10 (alternatively, use [22, Thm. 11.2.5] and the obvious embedding $\mathcal{S}(\mathbb{R}^{2d}) \hookrightarrow L^1_{v_{r,s}}(\mathbb{R}^{2d})$). By duality, it follows that the reservoir is continuously embedded into the the space of tempered distributions $\mathcal{S}^\sim$, meaning that the modulation spaces themselves are spaces of tempered distributions.

***Remark 4.2.16.*** The coorbit theory of this section, regarding the Schrödinger representation and its voice transform, can be understood as the coorbit theory of the short-time Fourier transform. After checking that the Heisenberg group and Schrödinger representation are well-defined and satisfy all necessary assumptions, we were always able to replace the voice transform by the the short-time Fourier transform.

Christensen [2] generalized this concept to *projective representations*. These are unitary representations that are only well-defined up to a multiplicative constant of absolute value 1. That is, a projective representation $\pi$ of $\mathcal{G}$ satisfies $\pi(xy) = c_{x,y}\pi(x)\pi(y)$ with some (continuous) factor $c_{x,y} \in S^1$ for all $x, y \in \mathcal{G}$. Similar to what we have done with the reduced Heisenberg group, it is possible to define a (locally compact) group structure on $\widetilde{\mathcal{G}} = \mathcal{G} \times S^1$ such that $\pi$ induces a proper representation $\widetilde{\pi}$ of $\widetilde{\mathcal{G}}$, taking care of the factor $c$. Then, the coorbit theory of $\widetilde{\pi}$ can be understood as the coorbit theory of $\pi$.

The short-time Fourier transform arises in this sense from the projective representation $\pi(x, \omega) = M_\omega T_x$ of $\mathcal{G} = \mathbb{R}^d \times \mathbb{R}^d$ on $L^2(\mathbb{R}^d)$. ◁

CHAPTER 5

Banach Frames for Coorbit Spaces

The goal of this chapter is to describe the coorbit spaces and their elements in terms of sequence spaces. That means we want to characterize the membership of some $f$ to the coorbit space $\mathcal{C}o_m^p$ by the membership of an associated sequence $(c_i)_{i \in I}$ to a sequence space $\ell_{\widetilde{m}}^p(I)$, such that the coorbit norm and the sequence norm are equivalent. For (separable) Hilbert spaces, frames achieve this by describing an equivalent norm through coefficient sequences in $\ell^2$. For coorbit theory, we want to generalize this concept to Banach spaces. This leads to two different approaches, corresponding to the two dual concepts of *analysis* and *synthesis*.

The synthesis approach leads to *atomic decompositions* [15, Thm. 6.1] of coorbit spaces. There we write every $f \in \mathcal{C}o_m^p$ as a series $\sum_{i \in I} c_i \pi(x_i)\psi$ with a suitable discrete family $(x_i)_{i \in I} \subset \mathcal{G}$ and a suitable analyzing vector $\psi \in \mathcal{A}_w$. The coefficients $(c_i)_{i \in I}$ are not unique, but there is a bounded linear operator $D$ that maps each $f \in \mathcal{C}o_m^p$ to one of such coefficient sequences. In short, this approach utilizes this *coefficient operator* $D : \mathcal{C}o_m^p \to \ell_{\widetilde{m}}^p(I)$ to establish a norm equivalence $\|f\|_{\mathcal{C}o_m^p} \asymp \|Df\|_{\ell_{\widetilde{m}}^p}$, such that every $f$ can be written with the bounded *synthesis operator*

$$S : \ell_{\widetilde{m}}^p(I) \to \mathcal{C}o_m^p, \quad (c_i)_{i \in I} \mapsto \sum_{i \in I} c_i \, \pi(x_i)\psi$$

as $f = S(Df)$.

The analysis approach leads to *Banach frames* in the sense of Gröchenig [21]. We consider for each $f \in \mathcal{C}o_m^p$ the sequence $(\langle f, \pi(x_i)\psi \rangle_{\mathcal{R}_w \times \mathcal{H}_w^1})_{i \in I}$, again for suitable points $(x_i)_{i \in I} \subset \mathcal{G}$ and $\psi \in \mathcal{A}_w$. These coefficients are contained in some $\ell_{\widetilde{m}}^p(I)$ and determine $f$ uniquely. In short, this approach utilizes the *analysis operator*

$$C : \mathcal{C}o_m^p \to \ell_{\widetilde{m}}^p(I), \quad f \mapsto (\langle f, \pi(x_i)\psi \rangle_{\mathcal{R}_w \times \mathcal{H}_w^1})_{i \in I}$$

to establish a norm equivalence $\|f\|_{\mathcal{C}o_m^p} \asymp \|Cf\|_{\ell_{\widetilde{m}}^p}$, such that every $f$ can be recovered through some bounded *reconstruction operator* $R : \ell_{\widetilde{m}}^p(I) \to \mathcal{C}o_m^p$ as $f = R(Cf)$.

The *synthesis* and *analysis* operators are described through the family $(\pi(x_i)\psi)_{i \in I}$, while the *coefficient* and *reconstruction operators* are only ought to exist and be bounded. That is why both approaches are not compatible in general. Using the synthesis approach, the coefficients $Df$ are not necessarily of the form $(\langle f, \pi(x_i)\phi \rangle_{\mathcal{R}_w \times \mathcal{H}_w^1})_{i \in I}$. Using the analysis approach, the reconstruction operator is not necessarily given by a series $\sum_{i \in I} c_i \pi(x_i) \phi$. However, under stronger assumption on $\psi$ and the family $(x_i)_{i \in I}$, it is possible for $(\pi(x_i)\psi)_{i \in I}$ to be an atomic decomposition as well as a Banach frame, such that there exists a dual family $(e_i)_{i \in I} \subset \mathcal{H}_w^1$ with

$$f = \sum_{i \in I} \langle f, \pi(x_i)\psi \rangle_{\mathcal{R}_w \times \mathcal{H}_w^1} e_i = \sum_{i \in I} \langle f, e_i \rangle_{\mathcal{R}_w \times \mathcal{H}_w^1} \pi(x_i)\psi$$

for all $f \in \mathcal{C}o_m^p$ [21, Thm. U].

We will focus our attention on Banach frames and the analysis approach, using [4; 21] as our main sources. For results regarding atomic decompositions of coorbit spaces, we refer the reader to [15; 21].

We first need to give a precise definition of the term *Banach frame*.

**Definition 5.0.1.** Let $B$ be a Banach space. We call a countable family $(g_i)_{i \in I} \in B'$ in the dual space of $B$ a *Banach frame* if there exists a Banach space of sequences $B_d \subset \mathbb{C}^I$ such that the following is true:

  i. *Analysis:* The analysis operator

  $$C : B \to B_d, \quad f \mapsto (\langle f, g_i \rangle_{B \times B'})_{i \in I}$$

  is well-defined and establishes an equivalent norm on $B$, i.e. there exist some constants $0 < a \leq b < \infty$ such that

  $$a \|f\|_B \leq \left\| (\langle f, g_i \rangle_{B \times B'})_{i \in I} \right\|_{B_d} \leq b \|f\|_B$$

  for all $f \in B$. The constants $a$ and $b$ are called *frame bounds*.

  ii. *Reconstruction:* There is a bounded reconstruction operator $R : B_d \to B$ that satisfies $RC = \mathrm{id}_B$. That is, we have

  $$R\left( (\langle f, g_i \rangle_{B \times B'})_{i \in I} \right) = f$$

  for all $f \in B$.

The reconstruction operator $R$ already establishes the lower frame bound $a = \|R\|_{B_d \to B}^{-1}$, so the above definition essentially breaks down to the boundedness of the analysis operator and the existence of a bounded reconstruction operator.

**Remark 5.0.2.** A usual frame on a Hilbert space is also a Banach frame in the above sense. In that case, the sequence space is $\ell^2(I)$, thus the first condition is exactly the frame condition for Hilbert spaces.

The second condition already follows in this context from the first one, since any frame $(g_i)_{i \in I}$ in a Hilbert space gives rise to a dual frame $(\widetilde{g}_i)_{i \in I}$, such that $(c_i)_{i \in I} \mapsto \sum_{i \in I} c_i \widetilde{g}_i$ is a bounded reconstruction operator. For general Banach spaces, this is not the case, so the second condition is not redundant in general. ◁

The aim of this chapter is to construct Banach frames of the form $(\pi(x_i)\psi)_{i \in I}$ for $\mathcal{C}o_m^p$ with respect to the sequence space $\ell_{\widetilde{m}}^p(I)$, where we use the discrete weight function $\widetilde{m}(i) = m(x_i)$ on $I$. The coorbit spaces $\mathcal{C}o_m^p$ and $\mathcal{C}o_{1/m}^q$ with $p^{-1} + q^{-1} = 1$ are antidual to each other by the pairing

$$\langle f, g \rangle_{\mathcal{C}o_m^p \times \mathcal{C}o_{1/m}^q} = \int_{\mathcal{G}} \mathcal{V}_\psi f(y) \overline{\mathcal{V}_\psi g(y)} \, dy,$$

as we have seen in proposition 3.4.7. Since $(\pi(x_i)\psi)_{i\in I} \subset \mathcal{A}_w \subset \mathcal{C}o^q_{1/m}$, the frame vectors are indeed contained in $(\mathcal{C}o^p_m)^\sim$. Definition 5.0.1 is therefore applicable.

A major observation is now that the coefficient $\langle f, \pi(x_i)\psi\rangle_{\mathcal{R}_w \times \mathcal{H}^1_w}$ is just the voice transform $\mathcal{V}_\psi f$ evaluated at the point $x_i$. The property of $(\pi(x_i)\psi)_{i\in I}$ being a Banach frame can therefore be understood as the property of $\mathcal{V}_\psi f$ being 'stable' under sampling in the points $(x_i)_{i\in I}$.

We will now proceed as follows. In section 5.1, we take a closer look at the discrete families $(x_i)_{i\in I}$ and how they need to be distributed in $\mathcal{G}$. We will also relate the discrete spaces $\ell^p_{\widetilde{m}}(I)$ to the continuous spaces $L^p_m(\mathcal{G})$ using these families of points. In section 5.2, we investigate the oscillation of a function, which describes 'smoothness' of functions in a way we can utilize. After that, we are able in section 5.3 to show that certain function spaces with reproducing kernels are stable under sampling, provided the kernel is smooth enough. In section 5.4, we can finally use this stability under sampling in connection with the correspondence principle to construct Banach frames for coorbit spaces.

## 5.1 Discretization of the Underlying Group

We begin with discretizing the locally compact group $\mathcal{G}$ and its associated weighted $L^p$-spaces. In order to do that, we need families of points in $\mathcal{G}$ that are evenly distributed in a geometric sense.

**Definition 5.1.1.** *[31, Def. 2.3.8]* Let $X = (x_i)_{i\in I}$ be a countable family of points in $\mathcal{G}$ and $U \subset \mathcal{G}$ relatively compact.

a) We call $X$ *$U$-dense* if the family $(x_i U)_{i\in I}$ covers $\mathcal{G}$.

b) We call $X$ *relatively separated* if for every compact set $K \subset \mathcal{G}$ the supremum

$$C_K = \sup_{i \in I} |\{j \in I \mid x_i K \cap x_j K \neq \emptyset\}|$$

is finite, where $|M|$ denotes the cardinality of the set $M$.

c) We call $X$ *$U$-well-spread* if it is $U$-dense and relatively separated.

We first show that for every relatively compact neighbourhood of the neutral element $e \in \mathcal{G}$ there exists a $U$-well-spread family, so that the developed theory is not empty.

**Lemma 5.1.2.** *[12, Lemma 1]* *For every relatively compact $e$-neighbourhood $U \subset \mathcal{G}$, there exists a $U$-well-spread family $(x_i)_{i\in I}$.*

*Proof.* We first construct an open $e$-neighbourhood $V \subset \mathcal{G}$ that satisfies $V^2 \subset U$ and $V = V^{-1}$. The mapping $\mathcal{G} \times \mathcal{G} \to \mathcal{G}$, $(x,y) \mapsto xy$ is continuous, thus the preimage of $U$ under this mapping is open (in the product topology). This preimage contains the tuple $(e,e)$, and with that a rectangle $W_1 \times W_2 \subset \mathcal{G} \times \mathcal{G}$ of open $e$-neighbourhoods $W_1, W_2 \subset \mathcal{G}$ as well. The open set $V_0 = W_1 \cap W_2$ now satisfies $V_0^2 \subset U$. $V_0$ and $V_0^{-1}$ are both open $e$-neighbourhoods, so $V = V_0 \cap V_0^{-1}$ has the desired properties.

We now consider the set of all at most countable families $(y_i)_{i\in I} \subset \mathcal{G}$ that satisfy

$$y_i V \cap y_j V = \emptyset \text{ for all } i,j \in I, i \neq j.$$

This set is partially ordered by inclusion of families. It is non-empty, since it contains all single-point families. The chains (i.e. totally ordered subsets) in this set are all bounded from above by the union of their countably many families.

We now apply Zorn's lemma to get a maximal element $X = (x_i)_{i \in I}$ in this set. We will show that this maximal element is the desired $U$-well spread family.

We begin by proving the $U$-density. Assume that $x \in \mathcal{G}$ is not contained in any of the sets $x_i U$, $i \in I$. If there were $j \in I$ such that $xV \cap x_j V \neq \emptyset$, there would exist $v_1, v_2 \in V$ such that $x = x_j v_1 v_2^{-1} \in x_j U$, contradicting the assumption. So $xV \cap x_i V = \emptyset$ for all $i \in I$. But this contradicts the maximality of $X$, so $X$ is indeed $U$-dense.

Next we show the relative separateness of $X$. Let $K \subset \mathcal{G}$ be compact and $i \in I$ fixed for now. Assume that for some $j \in I$ the intersection $x_i K \cap x_j K$ is non-empty. Then there exist $k_1, k_2 \in K$ such that $x_j = x_i k_1 k_2^{-1}$. This implies the inclusion $x_j V \subset x_i K K^{-1} V$.

The set $x_i K K^{-1} V$ is relatively compact (because this is the case for $V \subset U$), thus it has a finite Haar measure. It can therefore only cover finitely many of the disjoint sets $x_j V$, $j \in I$, namely at most $\mu(KK^{-1}V)/\mu(V)$. It follows that

$$|\{j \in I \mid x_i K \cap x_j K \neq \emptyset\}| \leq |\{j \in I \mid x_j \subset x_i KK^{-1}V\}| \leq \mu(KK^{-1}V)/\mu(V)$$

for all $i \in I$, so $X$ is relatively separated. $\square$

We want to show in this section that for $U$-well spread families $(x_i)_{i \in I}$, moderate weights $m$ and $1 \leq p \leq \infty$, the two expressions

$$\left(\sum_{i \in I} |c_i|^p m(x_i)^p\right)^{1/p} \quad \text{and} \quad \left\|\sum_{i \in I} c_i \chi_{x_i U}\right\|_{L_m^p}$$

for $(c_i)_{i \in I} \subset \mathbb{C}$ define equivalent norms on the sequence space $\ell_{\widetilde{m}}^p(I)$ described by them. For that, we need the fact that the series on the right side consists of at most $N \in \mathbb{N}$ summands at each point. This follows as the special case $K = \{e\}$ from the next lemma.

**Lemma 5.1.3.** *[31, Lemma 2.3.10] Let $(x_i)_{i \in I} \subset \mathcal{G}$ be relatively separated and $U \subset \mathcal{G}$ relatively compact. Then for each compact set $K \subset \mathcal{G}$ the supremum*

$$\sup_{x \in \mathcal{G}} |\{j \in I \mid xK \cap x_j U \neq \emptyset\}| \qquad (5.1)$$

*is finite. In particular, the family $(x_i U)_{i \in I}$ is locally finite.*

*Proof.* Let $x \in \mathcal{G}$. Assume that for some $i \in I$ the intersection $xK \cap x_i U$ is non-empty. Then there are $k \in K$ and $u \in U$ such that $xk = x_i u$ and therefore $x = x_i u k^{-1}$. Thus we have $x \in x_i U K^{-1} \subset x_i \overline{U} K^{-1}$, where $\overline{U}$ is compact.

If $xK \cap x_j U$ is empty for all $j \in I$, $x$ does not contribute to the supremum in (5.1). We may therefore assume that $xK \cap x_j U \neq \emptyset$ for some $j \in I$. It then follows

$$\{i \in I \mid xK \cap x_i U \neq \emptyset\} \subset \{i \in I \mid x \in x_i \overline{U} K^{-1}\}$$
$$\subset \{i \in I \mid x_i \overline{U} K^{-1} \cap x_j \overline{U} K^{-1} \neq \emptyset\}.$$

Since $(x_i)_{i \in I}$ is relatively separated, the last set contains at most $C_{\overline{U}K}$ elements, where $C_{\overline{U}K}$ is the constant from definition 5.1.1(b). This constant is independent from $x$, so the supremum over all $x \in \mathcal{G}$ is finite.

To prove the local finiteness of $(x_i U)_{i \in I}$, we take any compact neighbourhood $K$ of $x$. This $K$ intersects only finitely many of the sets $(x_i U)$, $i \in I$, as we have just proven, so $(x_i U)_{i \in I}$ is locally finite. □

Next we show that the moderate weight $m$ is compatible with the discretization, that is, $m(y)$ is close to $m(x_i)$ provided $y$ is close to $x_i$.

**Lemma 5.1.4.** *[31, Cor. 2.2.23] Let $m$ be a moderate weight function on $\mathcal{G}$ and $U \subset \mathcal{G}$ a relatively compact neighbourhood of $e$. Then there are constants $0 < a \leq b < \infty$ such that for every $U$-well-spread family $(x_i)_{i \in I}$ we have*

$$am(x_i) \leq m(y) \leq bm(x_i)$$

*for all $i \in I$ and $y \in x_i U$.*

*Proof.* Let $i \in I$. Then we have for every $u \in U$

$$m(x_i) = m(x_i u u^{-1}) \leq m(x_i u) \beta_0(u^{-1}).$$

Thus we have the estimate

$$m(x_i) \leq m(y) \sup_{u \in U} \beta_0(u^{-1})$$

for all $y \in x_i U$. Since $\beta_0$ is locally bounded (lemma 2.2.4) and $U^{-1}$ is relatively compact, this proves the first inequality with $a^{-1} = \sup_{u \in U} \beta_0(u^{-1})$. The other inequality follows similarly. □

We are now able to show the claimed equivalence of norms.

**Proposition 5.1.5.** *[4, Lemma 4.10] Let $U$ be a relatively compact neighbourhood of $e$ and $(x_i)_{i \in I}$ a relatively separated family in $\mathcal{G}$. Furthermore, let $m$ be a moderate weight on $\mathcal{G}$ and $\widetilde{m}(i) = m(x_i)$, $i \in I$ its discretization on $I$. Then*

$$\left\| \sum_{i \in I} |c_i| \chi_{x_i U} \right\|_{L_m^p}, \quad (c_i)_{i \in I} \in \ell_{\widetilde{m}}^p(I)$$

*defines an equivalent norm on $\ell_{\widetilde{m}}^p(I)$. In particular, a sequence $(c_i)_{i \in I}$ of complex numbers is contained in $\ell_{\widetilde{m}}^p(I)$ if and only if the function $\sum_{i \in I} |c_i| \chi_{x_i U}$ is contained in $L_m^p(\mathcal{G})$.*

*Proof.* Let $(c_i)_{i \in I} \subset \mathbb{C}$. If we apply lemma 5.1.3 to $K = \{e\}$, it follows that the series

$$\sum_{i \in I} |c_i| \chi_{x_i U} \tag{5.2}$$

consists pointwise of at most $N$ summands for some constant $N \in \mathbb{N}$. This means that the series gives a pointwise well-defined function. We now have to show that the inequalities

$$a \|(c_i)_{i \in I}\|_{\ell_{\widetilde{m}}^p} \leq \left\| \sum_{i \in I} |c_i| \chi_{x_i U} \right\|_{L_m^p} \leq b \|(c_i)_{i \in I}\|_{\ell_{\widetilde{m}}^p}$$

hold for suitable constants $0 < a \leq b < \infty$.

We first assume $p \in [1, \infty)$. Since all norms on $\mathbb{C}^N$ are equivalent, we have the inequalities

$$a_1 \left(\sum_{k=1}^{N} |v_k|^p\right)^{1/p} \leq \sum_{k=1}^{N} |v_k| \leq b_1 \left(\sum_{k=1}^{N} |v_k|^p\right)^{1/p}$$

for some constants $0 < a_1 \leq b_1 < \infty$ and all finite sequences $v_1, \ldots, v_N \in \mathbb{C}$. This estimates hold pointwise for the series (5.2), so it follows

$$a_1^p \sum_{i \in I} |c_i|^p \chi_{x_i U} \leq \left(\sum_{i \in I} |c_i| \chi_{x_i U}\right)^p \leq b_1^p \sum_{i \in I} |c_i|^p \chi_{x_i U}. \tag{5.3}$$

We also need the estimate

$$a_2 m(x_i) \leq m(y) \leq b_2 m(x_i) \tag{5.4}$$

for suitable $0 < a_2 \leq b_2 < \infty$ and $y \in x_i U$ from lemma 5.1.4.

We abbreviate $E = \mu(U)$, so we also have $E = \mu(x_i U)$ for all $i \in I$. By using inequalities (5.3) and (5.4) we get

$$a_1^p a_2^p \|(c_i)_{i \in I}\|_{\ell_{\widetilde{m}}^p}^p = \sum_{i \in I} a_1^p a_2^p |c_i|^p m(x_i)^p$$

$$= \frac{1}{E} \sum_{i \in I} \int_{\mathcal{G}} a_1^p a_2^p |c_i|^p \chi_{x_i U}(y) m(x_i)^p \, dy$$

$$\leq \frac{1}{E} \int_{\mathcal{G}} \sum_{i \in I} a_1^p |c_i|^p \chi_{x_i U}(y) m(y)^p \, dy$$

$$\leq \frac{1}{E} \int_{\mathcal{G}} \left(\sum_{i \in I} |c_i| \chi_{x_i U}(y)\right)^p m(y)^p \, dy$$

$$= \frac{1}{E} \left\|\sum_{i \in I} |c_i| \chi_{x_i U}\right\|_{L_m^p}^p,$$

which proves the lower bound. The upper bound follows similarly.

For $p = \infty$ we have

$$a_2 \|(c_i)_{i \in I}\|_{\ell_{\widetilde{m}}^\infty} = \sup_{i \in I} |c_i| a_2 m(x_i)$$

$$\leq \sup_{y \in \mathcal{G}} \sum_{i \in I} |c_i| \chi_{x_i U}(y) a_2 m(x_i)$$

$$\leq \sup_{y \in \mathcal{G}} \sum_{i \in I} |c_i| \chi_{x_i U}(y) m(y)$$

$$= \left\|\sum_{i \in I} |c_i| \chi_{x_i U}\right\|_{L_m^\infty}$$

$$\leq N b_2 \sup_{i \in I} |c_i| m(x_i)$$

$$= N b_2 \|(c_i)_{i \in I}\|_{\ell_{\widetilde{m}}^\infty},$$

finishing the proof. □

## 5.2 The Oscillation of a Function

We have already noticed that the Banach frames of the coorbit spaces are of the form $(\mathcal{V}_\psi f(x_i))_{i \in I}$, meaning the functions $\mathcal{V}_\psi f$ have to be stable under sampling in the points $(x_i)_{i \in I}$. This is ensured by using a reproducing kernel $\mathcal{V}_\psi \psi$ with a certain 'smoothness'. The oscillation is a suitable way to describe this kind of smoothness.

**Definition 5.2.1.** Let $G : \mathcal{G} \to \mathbb{C}$ be a measurable function and $U \subset \mathcal{G}$ a compact neighbourhood of $e$. Then we call

$$\operatorname{osc}_U(G) : \mathcal{G} \to [0, \infty], \quad \operatorname{osc}_U(G)(x) = \operatorname*{ess\,sup}_{u \in U} |G(ux) - G(x)|$$

the *U-oscillation of $G$*, or simply the *oscillation of $G$*.

The following estimate will be useful to us.

**Lemma 5.2.2.** *[21, Lemma 4.6] Let $G : \mathcal{G} \to \mathbb{C}$ be measurable, $U$ a compact neighbourhood of the neutral element $e$ and $x, y, z \in \mathcal{G}$ such that $y \in xU$. Then*

$$|G(y^{-1}z) - G(x^{-1}z)| \leq \operatorname{osc}_U(G)(y^{-1}z).$$

*Proof.* Writing $y = xu$ for some $u \in U$ we get

$$|G(y^{-1}z) - G(x^{-1}z)| = |G(y^{-1}z) - G(uy^{-1}z)| \leq \operatorname{osc}_U(G)(y^{-1}z).$$

□

For the remainder of this section we need a submultiplicative weight function $w$ on $\mathcal{G}$. The $L^1_w(\mathcal{G})$-norm of the $U$-oscillation of the reproducing kernel will be important in the following theory. We first show that the finiteness of this norm is independent from the chosen $e$-neighbourhood $U$.

**Lemma 5.2.3.** *Let $G \in L^1_w(\mathcal{G})$ and $U \subset \mathcal{G}$ a compact neighbourhood of $e$ such that $\operatorname{osc}_U(G) \in L^1_w(\mathcal{G})$. Then $\operatorname{osc}_K(G) \in L^1_w(\mathcal{G})$ for all compact $e$-neighbourhoods $K$.*

*Proof.* Let $K$ be an arbitrary compact $e$-neighbourhood. $K$ is covered by the open sets $K \subset \bigcup_{x \in K} U^o x$, where $U^o$ denotes the topological interior of $U$. Thus there are finitely many points $x_1, \ldots, x_n \in K$ such that

$$K \subset \bigcup_{k=1}^n U^o x_k \subset \bigcup_{k=1}^n U x_k.$$

Now we have

$$\left\| \operatorname{osc}_K(G) \right\|_{L^1_w} = \int_\mathcal{G} \sup_{v \in K} |G(vy) - G(y)| w(y) \, dy$$

$$\leq \int_\mathcal{G} \left( \sup_{k=1,\ldots,n} \sup_{u \in U} |G(ux_k y) - G(y)| \right) w(y) \, dy$$

$$\leq \int_\mathcal{G} \left( \sup_{k=1,\ldots,n} \sup_{u \in U} |G(ux_k y) - G(x_k y)| + |G(x_k y)| + |G(y)| \right) w(y) \, dy$$

$$\leq \int_{\mathcal{G}} \left( \sum_{k=1}^{n} \left[ \sup_{u \in U} |G(ux_k y) - G(x_k y)| + |G(x_k y)| \right] + |G(y)| \right) w(y)\, dy$$

$$\leq \int_{\mathcal{G}} \left( \sum_{k=1}^{n} \operatorname{osc}_U(G)(x_k y) + |G(x_k y)| \right) w(y)\, dy + \|G\|_{L^1_w}.$$

By pulling the sum out of the integral and substituting $y \mapsto x_k^{-1} y$ we get

$$\left\| \operatorname{osc}_K(G) \right\|_{L^1_w} \leq \sum_{k=1}^{n} \int_{\mathcal{G}} \left( \operatorname{osc}_U(G)(x_k y) + |G(x_k y)| \right) w(y)\, dy + \|G\|_{L^1_w}$$

$$= \sum_{k=1}^{n} \int_{\mathcal{G}} \left( \operatorname{osc}_U(G)(y) + |G(y)| \right) w(x_k^{-1} y)\, dy + \|G\|_{L^1_w}$$

$$\leq \sum_{k=1}^{n} w(x_k^{-1}) \left( \left\| \operatorname{osc}_U(G) \right\|_{L^1_w} + \|G\|_{L^1_w} \right) + \|G\|_{L^1_w}.$$

The last term is finite, so $\operatorname{osc}_K(G) \in L^1_w(\mathcal{G})$. □

The following definition now makes sense.

**Definition 5.2.4.** Let $G \in L^1_w(\mathcal{G})$. We say $G$ has a *w-integrable oscillation*, if $\operatorname{osc}_U(G)$ is $w$-integrable for any (and thus for every) compact $e$-neighbourhood $U \subset \mathcal{G}$. We denote the set of functions with $w$-integrable oscillation by $W(L^1_w)$.

**Remark 5.2.5.** We use the notation $W(L^1_w)$ because the set coincides with the *Wiener amalgam space* with global component $L^1_w(\mathcal{G})$ [21, Lemma 4.6(i)]. The Wiener amalgam spaces are a class of function spaces that combine local and global behaviour of functions into a single norm.

It is possible to replace the oscillation methods in the following theory by methods utilizing Wiener amalgam spaces and their convolution relations [15; 31]. We opted to use the oscillation as it is less technical and the estimates are more straightforward. ◁

If the function $G \in W(L^1_w)$ is continuous, the $L^1_w(\mathcal{G})$-norm of the $U$-oscillation becomes arbitrarily small when choosing $U$ small enough.

**Lemma 5.2.6.** *[21, Lemma 4.6] Let $G \in L^1_w(\mathcal{G})$ be continuous with $w$-integrable oscillation. Then*
$$\lim_{U \to \{e\}} \left\| \operatorname{osc}_U(G) \right\|_{L^1_w} = 0,$$
*where $U$ runs through all compact neighbourhoods of $e$. That is, for every $\varepsilon > 0$ there is a compact neighbourhood $Q$ of $e$ such that $\left\| \operatorname{osc}_U(G) \right\|_{L^1_w} < \varepsilon$ for all compact $e$-neighbourhoods $U \subset Q$.*

*Proof.* Let $\varepsilon > 0$ and $V$ a compact $e$-neighbourhood. Since the norm $\left\| \operatorname{osc}_V(G) \right\|_{L^1_w}$ is finite and $\mathcal{G}$ is $\sigma$-compact, there exists some compact $K \subset \mathcal{G}$ such that

$$\int_{\mathcal{G} \setminus K} \operatorname{osc}_V(G)(x)\, w(x)\, dx \leq \frac{\varepsilon}{2}.$$

Now, the restriction of $G$ to $K$ is uniformly continuous. That means that for every $c > 0$ there is some $e$-neighbourhood $Q$ such that $|G(x) - G(y)| < c$ for all $x, y \in K$ with

$xy^{-1} \in Q$ [26, Ch. XII, Prop. 1.1]. We can therefore estimate the $Q$-oscillation of $y \in K$ by

$$\operatorname{osc}_Q(G)(y) = \sup_{u \in Q} |G(uy) - G(y)| < c. \tag{5.5}$$

Technically speaking, this inequality is only valid for $y \in \bigcap_{u \in Q} u^{-1}K$, as only then the group elements $uy$ are all contained in the compact set $K$. But if we take some compact $e$ neighbourhood $Q_0$, choose all constants regarding the uniform continuity of $G$ with respect to the larger set $Q_0 K$ (which is still compact), and assume all used neighbourhoods $Q$ to be contained in $Q_0$, the above estimate is still achievable.

According to [31, Thm. 2.2.22] submultiplicative weights are locally bounded, so the supremum $a = \sup_{x \in K} w(x)$ is finite. Thus we can choose $Q \subset V$ small enough such that (5.5) holds for $c = \varepsilon/(2\mu(K)a)$. Then we have for all compact $e$-neighbourhoods $U \subset Q$

$$\int_{\mathcal{G}} \operatorname{osc}_U(G)(y) w(y)\, dy = \int_{\mathcal{G} \setminus K} \operatorname{osc}_U(G)(y) w(y)\, dy + \int_K \operatorname{osc}_U(G)(y) w(y)\, dy$$
$$\leq \int_{\mathcal{G} \setminus K} \operatorname{osc}_V(G)(y) w(y)\, dy + \int_K \operatorname{osc}_Q(G)(y) w(y)\, dy$$
$$< \frac{\varepsilon}{2} + \mu(K) \frac{\varepsilon}{2\mu(K)a} a$$
$$= \varepsilon,$$

where we used the monotonicity $\operatorname{osc}_U \leq \operatorname{osc}_Q \leq \operatorname{osc}_K$. This finishes the proof. □

## 5.3 Banach Frames for Reproducing Kernel Spaces

In this section, we assume $G$ to be a function from $W(L^1_w) \cap L^q_{1/m}(\mathcal{G})$ that fulfils the equalities $G = G^\triangledown = G * G$. As usual, $m$ is supposed to be a moderate weight with the $p$-control-weight $w$, $1 \leq p \leq \infty$ and $p^{-1} + q^{-1} = 1$.

Then we define the *reproducing kernel space* associated to $G$ by

$$\mathcal{M}^p_m = \{F \in L^p_m(\mathcal{G}) \mid F = F * G\}. \tag{5.6}$$

It is a closed subspace of $L^p_m(\mathcal{G})$ since the convolution operator $F \mapsto F * G$ is bounded on $L^p_m(\mathcal{G})$ by Young inequality (2.8), thus it is a Banach space. Moreover, the convolution operator a bounded projection from $L^p_m(\mathcal{G})$ onto $\mathcal{M}^p_m$ since we have $F * G = F * (G * G) = (F * G) * G$ for all $F \in L^p_m(\mathcal{G})$.

The technical requirement $G \in L^q_{1/m}(\mathcal{G})$ is important for multiple reasons. To see that, we write the convolution $F * G$ for $F \in \mathcal{M}^p_m$ and $x \in \mathcal{G}$ as

$$(F * G)(x) = \int_{\mathcal{G}} F(y) G(y^{-1}x)\, dy = \int_{\mathcal{G}} F(y) \overline{G(x^{-1}y)}\, dy$$
$$= \langle F, L_x G \rangle_{L^p_m \times L^q_{1/m}} = \langle L_{x^{-1}} F, G \rangle_{L^p_m \times L^q_{1/m}}.$$

First of all, the integrals can really be understood as an antiduality between $L^p_m(\mathcal{G})$ and $L^q_{1/m}(\mathcal{G})$ as both spaces are translation invariant. In particular is the convolution integral absolutely convergent at every point. Since either $p$ or $q$ is finite, and the translation operator is strongly continuous for finite exponents by lemma 2.2.15, every function $F = F * G$ in $\mathcal{M}^p_m$ is continuous.

Now point evaluations in $\mathcal{M}_m^p$ are well-defined and can be expressed through the convolution with $G$. This yields the inequality

$$|F(x)| = |(F * G)(x)| \leq \|F\|_{L_m^p} \|L_x G\|_{L_{1/m}^q} \leq w(x) \|F\|_{L_m^p} \|G\|_{L_{1/m}^q},$$

which means that the evaluation functionals $E_x : \mathcal{M}_m^p \to \mathbb{C}$, $F \mapsto F(x)$ are bounded by $w(x)\|G\|_{L_{1/m}^q}$. This can be expressed in short by the continuous embedding $\mathcal{M}_m^p \hookrightarrow L_{1/w}^\infty(\mathcal{G})$. The space $\mathcal{M}_m^p$ can be seen as a generalization of the reproducing kernel Hilbert spaces from section 2.6.

The definition of $\mathcal{M}_m^p$ coincides with the definition used in (3.7) for the kernel $G = \mathcal{V}_\psi \psi$. There, the correspondence principle 3.4.4 states that $\mathcal{M}_m^p$ is isometrically isomorphic to the coorbit space $\mathcal{C}o_m^p$. Note that the proof of the correspondence principle also relies on the fact $\mathcal{V}_\psi \psi \in L_{1/m}^q(\mathcal{G})$.

We want to show in this section that for certain $G$, the space $\mathcal{M}_m^p$ possesses a Banach frame of the form $(L_{x_i}G)_{i\in I} \subset L_{1/m}^q(\mathcal{G})$ with respect to a suitable well-spread family $(x_i)_{i\in I}$. In this case, the analysis operator is just the sampling operator

$$C : \mathcal{M}_m^p \to \ell_{\tilde{m}}^p(I), \quad CF = \left(\langle F, L_{x_i} G\rangle_{L_m^p \times L_{1/m}^q}\right)_{i\in I} = (F(x_i))_{i\in I}.$$

Thus we have to show that the frame inequalities

$$a\|F\|_{L_m^p} \leq \|F(x_i)_{i\in I}\|_{\ell_{\tilde{m}}^p} \leq b\|F\|_{L_m^p}$$

hold, and that there exists a bounded reconstruction operator which recovers $F$ from its samples in $(x_i)_{i\in I}$. For both, the function $F$ has to possess a certain 'smoothness', which is why $G$ has to have a $w$-integrable oscillation.

We begin by showing that an upper frame bound $b$ exists. This only requires that the family $(x_i)_{i\in I}$ is relatively separated.

**Proposition 5.3.1** (Upper frame bound, [4, Lemma 4.14]). *Let $(x_i)_{i\in I} \subset \mathcal{G}$ be a relatively separated family and $U$ a compact neighbourhood of $e$. Then there exists a constant $0 < b < \infty$ such that*

$$\|F(x_i)_{i\in I}\|_{\ell_{\tilde{m}}^p} \leq b\|F\|_{L_m^p}$$

*holds for all $F \in \mathcal{M}_m^p$.*

*Proof.* Let $F \in \mathcal{M}_m^p$. According to proposition 5.1.5, the inequality

$$\|F(x_i)_{i\in I}\|_{\ell_{\tilde{m}}^p} \leq C \left\|\sum_{i\in I} |F(x_i)| \chi_{x_i U}\right\|_{L_m^p}$$

holds with some constant $C > 0$. Thus it suffices to estimate the $L_m^p$-norm of the series $\sum_{i\in I} |F(x_i)| \chi_{x_i U}$ by $C'\|F\|_{L_m^p}$. In order to do that, we first consider the series pointwise.

For $x \in \mathcal{G}$, the convolution relation $F = F * G$ implies

$$\sum_{i\in I} |F(x_i)| \chi_{x_i U}(x) = \sum_{i\in I} \left|\int_\mathcal{G} F(y) G(y^{-1} x_i)\, dy\right| \chi_{x_i U}(x)$$

$$\leq \int_\mathcal{G} |F(y)| \left(\sum_{i\in I} |G(y^{-1} x_i)| \chi_{x_i U}(x)\right) dy. \tag{5.7}$$

We take a closer look at the series in the brackets. Using the triangle inequality, it follows that

$$\sum_{i \in I} |G(y^{-1}x_i)|\chi_{x_iU}(x) \leq \sum_{i \in I} |G(y^{-1}x_i) - G(y^{-1}x)|\chi_{x_iU}(x) + \sum_{i \in I} |G(y^{-1}x)|\chi_{x_iU}(x)$$

$$= \sum_{i \in I} |G(x_i^{-1}y) - G(x^{-1}y)|\chi_{x_iU}(x) + \sum_{i \in I} |G(x^{-1}y)|\chi_{x_iU}(x),$$

where we used $G = G^\triangledown$ in the last step. Now each series has at most $N$ non-vanishing summands by lemma 5.1.3, where $N \in \mathbb{N}$ is independent of $x$. If a summand does not vanish in $x$, then $x \in x_iU$, thus we obtain with lemma 5.2.2

$$\sum_{i \in I} |G(y^{-1}x_i)|\chi_{x_iU}(x) \leq N \left( \mathrm{osc}_U(G)(x^{-1}y) + |G(x^{-1}y)| \right).$$

By applying this estimate to (5.7), we get

$$\sum_{i \in I} |F(x_i)|\chi_{x_iU}(x) \leq N \int_{\mathcal{G}} |F(y)| \left( \mathrm{osc}_U(G)(x^{-1}y) + |G(x^{-1}y)| \right) dy$$

$$= N(|F| * \mathrm{osc}_U(G)^\vee)(x) + N(|F| * |G^\vee|)(x).$$

Young inequality (2.8) now implies

$$\left\| F(x_i)_{i \in I} \right\|_{\ell^p_{\tilde{m}}} \leq C \left\| \sum_{i \in I} |F(x_i)|\chi_{x_iU} \right\|_{L^p_m}$$

$$\leq CN\|F\|_{L^p_m} \left( \left\| \mathrm{osc}_U(G) \right\|_{L^1_w} + \|G\|_{L^1_w} \right),$$

completing the proof. $\square$

The existence of a lower frame bound (and a reconstruction operator) is substantially more difficult to show. For that, we need a partition of unity that is compatible with some well-spread family.

**Definition 5.3.2.** [31, Def. 2.3.9] Let $U$ be a relatively compact neighbourhood of $e$. Then a family of measurable functions $\varphi_i : \mathcal{G} \to \mathbb{R}$, $i \in I$ is called *bounded uniform partition of unity of size $U$*, in short *$U$-BUPU*, if the following is true:

  i. It is $0 \leq \varphi_i \leq 1$ for all $i \in I$.

  ii. It is $\sum_{i \in I} \varphi_i \equiv 1$.

  iii. There is a relatively separated family $X = (x_i)_{i \in I}$ such that $\mathrm{supp}\,\varphi_i \subset x_iU$ for all $i \in I$.

The family $X$ is called the *localizing family of* $(\varphi_i)_{i \in I}$

A $U$-BUPU is a partition of unity subordinate to the cover $(x_iU)_{i \in I}$. In particular, $(x_iU)_{i \in I}$ is indeed a cover of $\mathcal{G}$, as otherwise the series $\sum_{i \in I} \varphi_i$ could not be non-zero everywhere. The localizing family $X$ is therefore implicitly assumed to be $U$-dense, and thus also $U$-well-spread.

Note that the series $\sum_{i \in I} \varphi_i$ of a $U$-BUPU consists pointwise of at most $N$ non-vanishing terms, where $N \in \mathbb{N}$ is independent from the point (by lemma 5.1.3).

**Remark 5.3.3.** If $U$ is any relatively compact neighbourhood of $e$, there exists a $U$-well-spread family $(x_i)_{i \in I}$ according to lemma 5.1.2. Then the functions

$$\varphi_i = \chi_{x_i U} \left( \sum_{j \in I} \chi_{x_j U} \right)^{-1}, \quad i \in I$$

define a $U$-BUPU with localizing family $(x_i)_{i \in I}$. Thus, we can always find a $U$-BUPU to a given $e$-neighbourhood $U$.

Continuous partitions of unity do also exist. If the group $\mathcal{G}$ is second-countable and $U$ is open, this fact already follows from general topology since the family $(x_i U)_{i \in I}$ is a locally finite cover of $\mathcal{G}$ (again for $U$-well-spread $(x_i)_{i \in I}$, [26, Ch. IX, Thm. 5.3]). This is the case for all Lie groups, and therefore for the affine group $\mathcal{A}\!f\!f$ and the reduced Heisenberg group $\mathbb{H}_r^d$.

It is also possible to construct continuous $U$-BUPUs in general, see for instance [12, Thm. 2]. ◁

The next statement extends the lower estimate of proposition 5.1.5 to $U$-BUPUs.

**Lemma 5.3.4.** *Let $U$ be a relatively compact neighbourhood of $e \in \mathcal{G}$ and $(\varphi_i)_{i \in I}$ be a $U$-BUPU with a localizing family $(x_i)_{i \in I}$. Then there is a constant $C > 0$ such that the estimate*

$$\left\| \sum_{i \in I} c_i \varphi_i \right\|_{L_m^p} \leq C \left\| (c_i)_{i \in I} \right\|_{\ell_{\widetilde{m}}^p}$$

*holds for all sequences $(c_i)_{i \in I} \subset \ell_m^p \widetilde{m}(I)$. In other words, the operator*

$$\ell_{\widetilde{m}}^p(I) \to L_m^p(\mathcal{G}), \quad (c_i)_{i \in I} \mapsto \sum_{i \in I} c_i \varphi_i$$

*is bounded.*

*Proof.* It is $\operatorname{supp} \varphi_i \subset x_i U$ for $i \in I$. This implies the pointwise inequality

$$\left| \sum_{i \in I} c_i \varphi_i \right| \leq \sum_{i \in I} |c_i| \varphi_i \leq \sum_{i \in I} |c_i| \chi_{x_i U}.$$

The statement now follows from proposition 5.1.5. □

We now use a $U$-BUPU $(\varphi_i)_{i \in I}$ to approximate functions $F \in \mathcal{M}_m^p$ by the series $\sum_{i \in I} F(x_i) \varphi_i$. This approximation is possible due to the $w$-integrable oscillation of the kernel $G$, which leads to small pointwise 'oscillations' of the function $F = F * G$.

**Lemma 5.3.5.** *[21, Thm. 4.11] Let $U$ be a compact neighbourhood of $e$ and $(\varphi_i)_{i \in I}$ a $U$-BUPU with localizing family $(x_i)_{i \in I}$. Then we have for all $F \in \mathcal{M}_m^p$, $i \in I$ and $y \in x_i U$ the estimates*

$$|F(y) - F(x_i)| \leq \left( |F| * \operatorname{osc}_U(G)^{\vee} \right)(y) \tag{5.8}$$

*as well as*

$$\left\| F - \sum_{j \in I} F(x_j) \varphi_j \right\|_{L_m^p} \leq \|F\|_{L_m^p} \left\| \operatorname{osc}_U(G) \right\|_{L_w^1}. \tag{5.9}$$

*Proof.* We first consider the difference $|F(y) - F(x_i)|$. We write $y = x_i u$ for some $u \in U$ and obtain
$$F(x_i) = F(yu^{-1}) = (F * G)(yu^{-1}) = (F * R_{u^{-1}}G)(y).$$
Therefore we have
$$|F(y) - F(x_i)| = |(F * (G - R_{u^{-1}}G))(y)| \leq (|F| * |G - R_{u^{-1}}G|)(y). \tag{5.10}$$

We can estimate the last term further by using the oscillation of $G$. For that, we evaluate the difference of $G$ and $R_{u^{-1}}G$ in $x \in \mathcal{G}$ and get
$$|G(x) - R_{u^{-1}}G(x)| = |G(x) - G(xu^{-1})| = |G(x^{-1}) - G(ux^{-1})| \leq \operatorname{osc}_U(G)(x^{-1}).$$

Thus it is $|G - R_{u^{-1}}G| \leq \operatorname{osc}_U(G)^\vee$. By inserting this in (5.10) it follows
$$|F(y) - F(x_i)| \leq (|F| * \operatorname{osc}_U(G)^\vee)(y),$$
so the first inequality is proven.

In order to show (5.9), we estimate for arbitrary $y \in \mathcal{G}$
$$\left|F(y) - \sum_{i \in I} F(x_i)\varphi_i(y)\right| = \left|\sum_{i \in I}(F(y) - F(x_i))\varphi_i(y)\right|$$
$$\leq \sum_{i \in I} |F(y) - F(x_i)|\varphi_i(y).$$

The summand with index $i$ can only be non-zero when $y \in \operatorname{supp} \varphi_i \subset x_i U$. Hence, we can apply inequality (5.8) to each summand individually to obtain
$$\left|F(y) - \sum_{i \in I} F(x_i)\varphi_i(y)\right| \leq \sum_{i \in I}(|F| * \operatorname{osc}_U(G)^\vee)(y)\varphi_i(y)$$
$$= (|F| * \operatorname{osc}_U(G)^\vee)(y).$$

Now an application of Young inequality (2.8) finishes the proof. $\square$

We are now finally able to show the existence of a lower frame bound and a reconstruction operator. We have to make additional assumptions about the reproducing kernel $G$, though.

**Proposition 5.3.6** (lower frame bound, [cf. 21, Thm. 5.3]). *Let $U$ be a compact neighbourhood of $e$ and $(x_i)_{i \in I}$ a $U$-well-spread family in $\mathcal{G}$ such that*
$$\|G\|_{L^1_w}\|\operatorname{osc}_U(G)\|_{L^1_w} < 1.$$

*Then there exists a lower frame bound $a > 0$ such that the inequality*
$$a\|F\|_{L^p_m} \leq \|F(x_i)_{i \in I}\|_{\ell^p_{\tilde{m}}}$$

*holds for all $F \in \mathcal{M}^p_m(\mathcal{G})$. Moreover, there is a bounded reconstruction operator $R : \ell^p_{\tilde{m}}(I) \to \mathcal{M}^p_m$ which satisfies $R(F(x_i)_{i \in I}) = F$.*

*Proof.* Suppose $(\varphi_i)_{i \in I}$ is a $U$-BUPU with localizing family $(x_i)_{i \in I}$; such a BUPU exists according to remark 5.3.3. We define the operator

$$T : \mathcal{M}_m^p \to \mathcal{M}_m^p, \quad TF = \sum_{i \in I} F(x_i) \varphi_i * G.$$

This operator is well-defined and bounded, as we have for $F \in \mathcal{M}_m^p$

$$\left\| \sum_{i \in I} F(x_i) \varphi_i * G \right\|_{L_m^p} \leq \left\| \sum_{i \in I} F(x_i) \varphi_i \right\|_{L_m^p} \|G\|_{L_w^1}$$
$$\leq C \left\| F(x_i)_{i \in I} \right\|_{\ell_{\widetilde{m}}^p} \|G\|_{L_w^1}$$
$$\leq bC \|F\|_{L_m^p} \|G\|_{L_w^1},$$

where we used Young inequality (2.8), lemma 5.3.4 and the upper frame bound from proposition 5.3.1.

We now consider the difference $\mathrm{id} - T : \mathcal{M}_m^p \to \mathcal{M}_m^p$, where $\mathrm{id}$ denotes the identity operator on $\mathcal{M}_m^p$. By applying this difference to some $F \in \mathcal{M}_m^p$, we obtain

$$\left\| (\mathrm{id} - T) F \right\|_{L_m^p} = \left\| F - \sum_{i \in I} F(x_i) \varphi_i * G \right\|_{L_m^p}$$
$$= \left\| \left( F - \sum_{i \in I} F(x_i) \varphi_i \right) * G \right\|_{L_m^p}$$
$$\leq \left\| F - \sum_{i \in I} F(x_i) \varphi_i \right\|_{L_m^p} \|G\|_{L_w^1}$$
$$\leq \|F\|_{L_m^p} \left\| \mathrm{osc}_U(G) \right\|_{L_w^1} \|G\|_{L_w^1},$$

where we used Young inequality (2.8) as well as the estimate (5.9). The operator norm of this difference can therefore be estimated by

$$\|\mathrm{id} - T\|_{\mathcal{M}_m^p \to \mathcal{M}_m^p} \leq \left\| \mathrm{osc}_U(G) \right\|_{L_w^1} \|G\|_{L_w^1} < 1,$$

implying that $T$ has a bounded inverse [26, Ch. IV, Thm. 2.1].

Now we can write $F \in \mathcal{M}_m^p$ as

$$F = T^{-1}(TF) = T^{-1} \left( \sum_{i \in I} F(x_i) \varphi_i * G \right).$$

Thus, the operator

$$R : \ell_{\widetilde{m}}^p(I) \to \mathcal{M}_m^p, \quad (c_i)_{i \in I} \mapsto T^{-1} \left( \sum_{i \in I} c_i \varphi_i * G \right)$$

is able to recover $F$ from its samples. The fact that $R$ is well-defined and bounded can be seen similarly as for $T$. The boundedness of $R$ implies in particular

$$\|F\|_{L_m^p} = \left\| R(F(x_i)_{i \in I}) \right\|_{L_m^p} \leq \|R\|_{\ell_{\widetilde{m}}^p \to \mathcal{M}_m^p} \left\| F(x_i)_{i \in I} \right\|_{\ell_{\widetilde{m}}^p}$$

for $F \in \mathcal{M}_m^p$, so the existence of a lower frame bound follows immediately, too. $\square$

Combining the two propositions 5.3.1 and 5.3.6, we finally get the frame statement.

**Theorem 5.3.7.** *Let $G \in W(L_w^1)$ and $U$ be a compact neighbourhood of $e$ such that*

$$\left\|\operatorname{osc}_U(G)\right\|_{L_w^1} \|G\|_{L_w^1} < 1. \tag{5.11}$$

*Then for every $U$-well-spread family $(x_i)_{i \in I}$ there are frame bounds $0 < a \leq b < \infty$ such that*

$$a\|F\|_{L_m^p} \leq \left(\sum_{i \in I} |F(x_i)|^p m(x_i)^p\right)^{1/p} \leq b\|F\|_{L_m^p}$$

*for all $F \in \mathcal{M}_m^p$. In addition, there is a bounded reconstruction operator $R : \ell_{\widetilde{m}}^p(I) \to \mathcal{M}_m^p$ with $R(F(x_i)_{i \in I}) = F$.*

*This means that $(L_{x_i} G)_{i \in I}$ is a Banach frame of $\mathcal{M}_m^p$ with respect to the sequence space $\ell_{\widetilde{m}}^p(I)$.*

To obtain a Banach frame through theorem 5.3.7, the inequality (5.11) must hold. If we assume $G$ to be continuous, we can apply lemma 5.2.6 to see that there exists a small enough compact $e$-neighbourhood $U$ such that that inequality is indeed satisfied. According to lemma 5.1.2, it is then possible to find a $U$-well-spread family $(x_i)_{i \in I}$. Thus, *if the reproducing kernel is continuous, the space $\mathcal{M}_m^p$ admits a Banach frame of the form $(L_{x_i} G)_{i \in I}$.*

## 5.4 Banach Frames for Coorbit Spaces

Now that we have established Banach frames of the form $(L_{x_i} G)_{i \in I}$ for the reproducing kernel spaces $\mathcal{M}_m^p$, we can use the correspondence principle pull them back to the coorbit spaces $\mathcal{C}o_m^p$.

Let $m$ be a moderate weight on $\mathcal{G}$ with the $p$-control-weight $w$ and let $\pi$ be a $w$-integrable unitary representation of $\mathcal{G}$ on some Hilbert space $\mathcal{H}$. According to the correspondence principle 3.4.4, the voice transform $\mathcal{V}_\psi$ of an admissible analyzing vector $\psi \in \mathcal{A}_w \setminus \{0\}$ defines an isometric isomorphism

$$\mathcal{V}_\psi : \mathcal{C}o_m^p \to \mathcal{M}_m^p = \{F \in L_m^p(\mathcal{G}) \mid F = F * \mathcal{V}_\psi \psi\}.$$

To apply the results of the preceding section 5.3, we have to make sure that the reproducing kernel $\mathcal{V}_\psi \psi$ satisfies all required conditions. This leads to the following definition.

**Definition 5.4.1.** We define the set of *w-atoms* by

$$\mathcal{B}_w = \{\psi \in \mathcal{H} \mid \mathcal{V}_\psi \psi \in W(L_w^1)\}$$
$$= \{\psi \in \mathcal{A}_w \mid \operatorname{osc}_U(\mathcal{V}_\psi \psi) \in L_w^1(\mathcal{G}) \text{ for some compact } e\text{-neighbourhood } U\}.$$

We have shown in lemma 3.4.3 that $\mathcal{V}_\psi \psi \in L_{1/m}^q(\mathcal{G})$ for all analyzing vectors $\psi$. Thus, for an admissible $w$-atom $\psi$, the kernel $\mathcal{V}_\psi \psi$ satisfies all technical prerequisites of section 5.3.

**Theorem 5.4.2.** *[21, Thm. 5.3] Let $\psi \in \mathcal{B}_w$ be an admissible $w$-atom, and let $U$ be a compact neighbourhood of the neutral element $e \in \mathcal{G}$ such that*

$$\big\|\mathrm{osc}_U(\mathcal{V}_\psi \psi)\big\|_{L^1_w} \big\|\mathcal{V}_\psi \psi\big\|_{L^1_w} < 1.$$

*Then, for any $U$-well-spread family $(x_i)_{i \in I} \subset \mathcal{G}$, the family $(\pi(x_i)\psi)_{i \in I}$ is a Banach frame of $\mathcal{C}o^p_m$ with respect to the sequence space $\ell^p_{\widetilde{m}}(I)$.*

*Proof.* Let $(x_i)_{i \in I}$ be a $U$-well-spread family. We can right away apply theorem 5.3.7 to obtain the Banach frame $(L_{x_i}\mathcal{V}_\psi \psi)_{i \in I}$ of the space

$$\mathcal{M}^p_m = \{F \in L^p_m(\mathcal{G}) \mid F = F * \mathcal{V}_\psi \psi\}$$

with respect to the sequence space $\ell^p_{\widetilde{m}}(I)$. Thus we have the frame inequalities

$$a\|F\|_{L^p_m} \leq \big\|F(x_i)_{i \in I}\big\|_{\ell^p_{\widetilde{m}}} \leq b\|F\|_{L^p_m}$$

for all $F \in \mathcal{M}^p_m$ and some constants $0 < a \leq b < \infty$. There also is a bounded reconstruction operator $R : \ell^p_{\widetilde{m}}(I) \to \mathcal{M}^p_m$.

Since the voice transform $\mathcal{V}_\psi : \mathcal{C}o^p_m \to \mathcal{M}^p_m$ is an isometric isomorphism, the frame inequalities imply

$$a\|\mathcal{V}_\psi f\|_{L^p_m} \leq \big\|\mathcal{V}_\psi f(x_i)_{i \in I}\big\|_{\ell^p_{\widetilde{m}}} \leq b\|\mathcal{V}_\psi f\|_{L^p_m}$$

and further

$$a\|f\|_{\mathcal{C}o^p_m} \leq \bigg\|\Big(\langle f, \pi(x_i)\psi\rangle_{\mathcal{R}_w \times \mathcal{H}^1_w}\Big)_{i \in I}\bigg\|_{\ell^p_{\widetilde{m}}} \leq b\|f\|_{\mathcal{C}o^p_m}$$

for all $f \in \mathcal{C}o^p_m$. This proves the frame inequalities for $(\pi(x_i)\psi)_{i \in I}$.

A bounded reconstruction operator $\widetilde{R} : \ell^p_{\widetilde{m}}(I) \to \mathcal{C}o^p_m$ is given by $\widetilde{R} = \mathcal{V}_\psi^{-1} R$. □

The reproducing kernel $\mathcal{V}_\psi \psi$ is a continuous function. Thus, provided it has a $w$-integrable oscillation, the norm $\big\|\mathrm{osc}_U(\mathcal{V}_\psi \psi)\big\|_{L^1_w}$ becomes arbitrarily small if $U$ is chosen small enough (lemma 5.2.6). If $\psi$ is an admissible $w$-atom, the requirements of theorem 5.4.2 can therefore be fulfilled by choosing such a $U$. This means that every admissible $w$-atom $\psi$ induces a Banach frame of the form $(\pi(x_i)\psi)_{i \in I}$ of $\mathcal{C}o^p_m$ for a suitable family $(x_i)_{i \in I}$. In the next lemma, we will see that $w$-atoms do exist for any $w$-integrable representation, thus *every coorbit space admits a Banach frame.*

**Lemma 5.4.3.** *[14, Lemma 6.1] We have $\mathcal{B}_w \neq \{0\}$.*

*Proof.* Let $F, G : \mathcal{G} \to \mathbb{C}$ be measurable functions, $U$ some compact $e$-neighbourhood and $x \in \mathcal{G}$. Then we can estimate the $U$-oscillation of $F * G$ by

$$\begin{aligned}
\mathrm{osc}_U(G * F)(x) &= \sup_{u \in U} \bigg| \int_\mathcal{G} G(y) F(y^{-1}x)\, dy - \int_\mathcal{G} G(y) F(y^{-1}ux)\, dy \bigg| \\
&= \sup_{u \in U} \bigg| \int_\mathcal{G} G(y) F(y^{-1}x)\, dy - \int_\mathcal{G} G(uy) F(y^{-1}x)\, dy \bigg| \\
&= \sup_{u \in U} \bigg| \int_\mathcal{G} (G(y) - G(uy)) F(y^{-1}x)\, dy \bigg| \\
&\leq \int_\mathcal{G} \sup_{u \in U} |G(y) - G(uy)|\, |F(y^{-1}x)|\, dy \\
&= (\mathrm{osc}_U(G) * |F|)(x).
\end{aligned}$$

With Young inequality (2.5) it follows

$$\left\|\mathrm{osc}_U(G*F)\right\|_{L^1_w} \leq \left\|\mathrm{osc}_U(G)\right\|_{L^1_w}\|F\|_{L^1_w} \tag{5.12}$$

for all $G \in W(L^1_w)$ and $F \in L^1_w(\mathcal{G})$, hence $G*F \in W(L^1_w)$.

In the proof of lemma 3.4.3, we have constructed an analyzing vector $\phi \in \mathcal{A}_w\backslash\{0\}$ that satisfies

$$\mathcal{V}_\phi\phi = \varphi * \mathcal{V}_\psi\psi * \varphi^\nabla$$

for some $\varphi \in C_c(\mathcal{G})$ and some $\psi \in \mathcal{A}_w\backslash\{0\}$. The function $\varphi$ has a $w$-integrable oscillation, since the compactness of its support implies that the integral

$$\int_\mathcal{G} \mathrm{osc}_U(\varphi)(y)\,w(y)\,dy = \int_\mathcal{G} \sup_{u \in U}|\varphi(y) - \varphi(uy)|w(y)\,dy$$

is finite. Using the fact that $\mathcal{V}_\psi\psi * \varphi^\nabla \in L^1_w(\mathcal{G})$ together with inequality (5.12), we obtain

$$\mathcal{V}_\phi\phi = \varphi * (\mathcal{V}_\psi\psi * \varphi^\nabla) \in W(L^1_w),$$

and therefore $\phi \in \mathcal{B}_w\backslash\{0\}$. $\square$

The essential prerequisite to obtain Banach frames through theorem 5.4.2 is the $w$-integrability of the oscillation of $\mathcal{V}_\psi\psi$, combined with the explicit estimate

$$\left\|\mathrm{osc}_U(\mathcal{V}_\psi\psi)\right\|_{L^1_w}\left\|\mathcal{V}_\psi\psi\right\|_{L^1_w} < 1.$$

However, such an inequality is not always easy to establish for specific examples, as it is necessary to find reasonably small upper bounds for the $L^1_w$-norm of $U$-oscillations. Additionally, the family $(x_i)_{i \in I}$ does depend on $\psi$ and $U$ and might need to lie quite densely.

Fortunately, there is another way to obtain Banach frames for coorbit spaces, namely by extending tight frames from the Hilbert space $\mathcal{H}$. A *tight frame for $\mathcal{H}$* is a family $(g_i)_{i \in I} \subset \mathcal{H}$ such that

$$c\|f\|_\mathcal{H} = \left(\sum_{i \in I}|\langle f, g_i\rangle_\mathcal{H}|^2\right)^{1/2} \tag{5.13}$$

holds for all $f \in \mathcal{H}$ and some constant $c > 0$. In particular, such a tight frame has the reconstruction property

$$f = c^{-2}\sum_{i \in I}\langle f, g_i\rangle_\mathcal{H} g_i$$

for all $f \in \mathcal{H}$, where the series converges unconditionally [3, Cor. 5.1.7]. This reconstruction property can be extended to the coorbit space $\mathcal{C}o^p_m$, which only requires the family $(x_i)_{i \in I}$ to be relatively separated.

**Theorem 5.4.4.** *[21, 6.6(c)] Let $\psi \in \mathcal{B}_w$ be an admissible $w$-atom and $(x_i)_{i \in I} \subset \mathcal{G}$ a relatively separated family. If $(\pi(x_i)\psi)_{i \in I}$ is a tight frame for $\mathcal{H}$, then it also is a Banach frame for $\mathcal{C}o^p_m$ with respect to $\ell^p_{\widetilde{m}}(I)$.*

*Proof.* The existence of an upper frame bound $b < \infty$ follows from proposition 5.3.1, applied to $G = \mathcal{V}_\psi\psi$ and $F = \mathcal{V}_\psi f$, $f \in \mathcal{C}o^p_m$. It remains to show that there exists a bounded reconstruction operator $R : \ell^p_{\widetilde{m}}(I) \to \mathcal{C}o^p_m$, as a lower frame bound is then given by $a = \|R\|^{-1}_{\ell^p_{\widetilde{m}} \to \mathcal{C}o^p_m}$.

We define the reconstruction operator weakly by

$$R : \ell^p_{\widetilde{m}}(I) \to \mathcal{C}o^p_m, \quad \langle R(c_i)_{i\in I}, g\rangle_{\mathcal{C}o^p_m \times \mathcal{C}o^q_{1/m}} = c^{-2} \sum_{i\in I} c_i \overline{\mathcal{V}_\psi g(x_i)}, \quad g \in \mathcal{C}o^q_{1/m},$$

where we use the antiduality of $\mathcal{C}o^p_m$ and $\mathcal{C}o^q_{1/m}$ from proposition 3.4.7, and where $c$ is the frame constant of the tight frame (5.13). According to proposition 5.3.1 we have

$$\sum_{i\in I} |\mathcal{V}_\psi g(x_i)|^q / m(x_i)^q \leq C \|\mathcal{V}_\psi g\|_{L^q_{1/m}} = C\|g\|_{\mathcal{C}o^q_{1/m}}$$

for some constant $C > 0$ and all $g \in \mathcal{C}o^q_{1/m}$, so by duality of $\ell^p_{\widetilde{m}}(I)$ and $\ell^q_{1/\widetilde{m}}(I)$ it follows that $R$ is well-defined and bounded.

It remains to show that $R$ actually is a reconstruction operator. For that, we need the fact that $(\pi(x_i)\psi)_{i\in I}$ is a tight frame, i.e. that

$$c^2 \|f\|^2_{\mathcal{H}} = \sum_{i\in I} |\langle f, \pi(x_i)\psi\rangle_{\mathcal{H}}|^2 = \sum_{i\in I} |\mathcal{V}_\psi f(x_i)|^2$$

is true for all $f \in \mathcal{H}$. Polarizing this identity yields

$$c^2 \langle f, g\rangle_{\mathcal{H}} = \sum_{i\in I} \mathcal{V}_\psi f(x_i) \overline{\mathcal{V}_\psi g(x_i)}$$

for all $f, g \in \mathcal{H}$. It follows that for $f \in \mathcal{C}o^p_m$ and $g \in \mathcal{C}o^q_{1/m}$ the equation

$$c^2 \langle f, g\rangle_{\mathcal{C}o^p_m \times \mathcal{C}o^q_{1/m}} = \sum_{i\in I} \mathcal{V}_\psi f(x_i) \overline{\mathcal{V}_\psi g(x_i)}$$

holds. To see that, we first use proposition 5.3.1 again to verify that both sides are bounded (anti)functionals in $f$ resp. $g$. Since the equation is true on the dense subsets $\mathcal{C}o^p_m \cap \mathcal{H} \supset \mathcal{H}^1_w$ and $\mathcal{C}o^q_{1/m} \cap \mathcal{H} \supset \mathcal{H}^1_w$ (see lemma 3.4.8 for the density of $\mathcal{H}^1_w$ in $\mathcal{C}o^p_m$ and $\mathcal{C}o^q_{1/m}$), it carries over to all such $f$ and $g$.

Now we get

$$\left\langle R\left[(\langle f, \pi(x_i)\psi\rangle_{\mathcal{R}_w \times \mathcal{H}^1_w})_{i\in I}\right], g\right\rangle_{\mathcal{C}o^p_m \times \mathcal{C}o^q_{1/m}} = c^{-2} \sum_{i\in I} \mathcal{V}_\psi f(x_i) \overline{\mathcal{V}_\psi g(x_i)}$$
$$= \langle f, g\rangle_{\mathcal{C}o^p_m \times \mathcal{C}o^q_{1/m}},$$

which implies $R(\langle f, \pi(x_i)\psi\rangle_{\mathcal{R}_w \times \mathcal{H}^1_w})_{i\in I} = f$. Thus, $R$ is indeed a reconstruction operator. □

**Remark 5.4.5.** The prerequisites of the two theorems 5.4.2 and 5.4.4 only depend on the $w$-atom $\psi$ and the $p$-control-weight $w$. This means that the family $(\pi(x_i)\psi)_{i\in I}$ from these theorems is a Banach frame *for all coorbit spaces* $\mathcal{C}o^p_m$ for which $w$ is a $p$-control-weight of $m$. The frame bounds do depend on $p$ and $m$, though. ◁

CHAPTER 6

Examples of Banach Frame Constructions

In chapter 4, we have examined the wavelet transform and the short-time Fourier transform, and constructed their respective coorbit spaces. We now aim to find examples of Banach frames for these spaces in the sense of chapter 5. This includes sufficient conditions for analyzing vectors to be atoms, as well as the definition of well-spread families of points that are typically used in the respective contexts.

## 6.1 Banach Frames and the Wavelet Transform

We described the wavelet transform in section 4.1 in terms of the wavelet representation $\pi$ of the affine group $\mathcal{A}\!f\!f$. We also identified its coorbit spaces with the homogeneous Besov spaces $\mathcal{C}o^p_{m_s} \cong \dot{B}^{s-1/2+1/p}_{p,p}$, where we used the moderate weights $m_s(b,a) = |a|^{-s}$ and the $p$-control-weights $w_\rho(b,a) = |a|^\rho + |a|^{-\rho}$.

In order to construct Banach frames for these coorbit spaces, we need to find families of points $(b_i, a_i)_{i \in I} \subset \mathcal{A}\!f\!f$ that lie well-spread with respect to a given compact neighbourhood of the neutral element $(0,1)$. This can be achieved through regular affine lattices, like for instance in [5, sec. 3.2.3.1; 13]. We need such lattices to cover both positive and negative scale factors $a > 0$ and $a < 0$ of the affine group, which is why we include the factor $\varepsilon = \pm 1$.

**Definition 6.1.1.** Let $\alpha > 1$ and $\beta > 0$. Then we define the *affine lattice* as the family

$$\Lambda(\beta, \alpha) = \left(\varepsilon\alpha^j \beta k, \varepsilon\alpha^j\right)_{j,k \in \mathbb{Z}, \varepsilon = \pm 1} \subset \mathcal{A}\!f\!f.$$

For fixed $\alpha > 1$ and $\beta > 0$, the lattice $\Lambda(\beta, \alpha)$ is well-spread with respect to the rectangle

$$A_{\beta,\alpha} = [-\beta/2, \beta/2] \times [\alpha^{-1/2}, \alpha^{1/2}], \tag{6.1}$$

which is a compact neighbourhood of $(0,1)$. In general, the point $\lambda \in \Lambda(\beta, \alpha)$ has the compact neighbourhood $\lambda A_{\beta,\alpha}$, see figure 6.1. Note that all these rectangles have the same Haar measure in $\mathcal{A}\!f\!f$.

**Proposition 6.1.2.** *Let $\alpha > 1$ and $\beta > 0$. Then $\Lambda(\beta, \alpha)$ is $A_{\beta,\alpha}$-well-spread.*

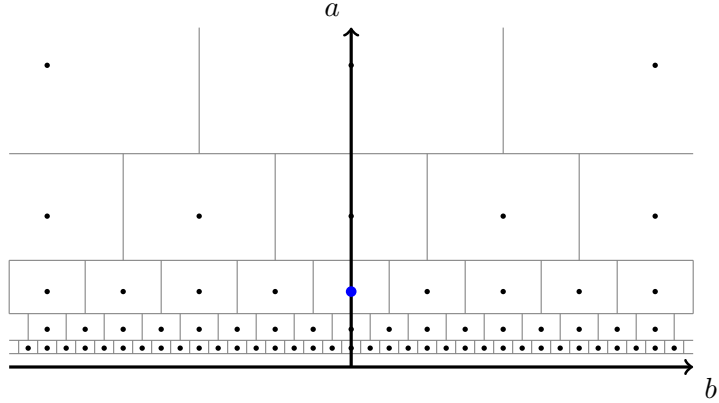

Figure 6.1: The subfamily $\varepsilon = 1$ of the affine lattice $\Lambda(1,2)$. The blue dot is the neutral element $(0,1)$. The rectangles are of the form $\lambda A_{1,2}$, where $\lambda \in \Lambda(1,2)$ is the center point of the respective rectangle.

*Proof.* We need to show that the family $\Lambda(\beta, \alpha)$ is $A_{\beta,\alpha}$-dense, that is

$$\bigcup_{\varepsilon = \pm 1} \bigcup_{j,k \in \mathbb{Z}} (\varepsilon \alpha^j \beta k, \varepsilon \alpha^j) A_{\beta,\alpha} = \mathcal{A}\!f\!f,$$

and that it is relatively separated, that is

$$\sup_{j,k \in \mathbb{Z}, \varepsilon = \pm 1} |\{(m, l, \eta) \in \mathbb{Z}^2 \times \{\pm 1\} \mid (\varepsilon \alpha^j \beta k, \varepsilon \alpha^j) K \cap (\eta \alpha^l \beta m, \eta \alpha^l) K \neq \emptyset\}| < \infty$$

for all compact $K \subset \mathcal{A}\!f\!f$.

For the $A_{\beta,\alpha}$-density we write the lattice points as products

$$(\varepsilon \alpha^j \beta k, \varepsilon \alpha^j) = (0, \varepsilon \alpha^j)(\beta k, 1).$$

Now we have

$$\bigcup_{k \in \mathbb{Z}} (\beta k, 1) A_{\beta,\alpha} = \bigcup_{k \in \mathbb{Z}} [\beta(k - 1/2), \beta(k + 1/2)] \times [\alpha^{-1/2}, \alpha^{1/2}] = \mathbb{R} \times [\alpha^{-1/2}, \alpha^{1/2}]$$

as well as

$$\bigcup_{\varepsilon = \pm 1} \bigcup_{j \in \mathbb{Z}} (0, \varepsilon \alpha^j) \left( \mathbb{R} \times [\alpha^{-1/2}, \alpha^{1/2}] \right) = \bigcup_{\varepsilon = \pm 1} \bigcup_{j \in \mathbb{Z}} \varepsilon \alpha^j \mathbb{R} \times \varepsilon \alpha^j [\alpha^{-1/2}, \alpha^{1/2}] = \mathbb{R} \times \mathbb{R}^*.$$

This already proves the $A_{\beta,\alpha}$-density.

To show that $\Lambda(\beta, \alpha)$ is relatively separated, let $K \subset \mathcal{A}\!f\!f$ be compact. Then $KK^{-1}$ is compact as well, so we have $KK^{-1} \subset [-\beta M, \beta M] \times \{\alpha^{-N} \leq |a| \leq \alpha^N\}$ for some $M, N \in \mathbb{N}$.

Let $j, k, l, m \in \mathbb{Z}$ and $\varepsilon, \eta \in \{\pm 1\}$. Then

$$(\varepsilon \alpha^j \beta k, \varepsilon \alpha^j) K \cap (\eta \alpha^l \beta m, \eta \alpha^l) K \neq \emptyset$$

implies that

$$(\eta \alpha^l \beta m, \eta \alpha^l)^{-1} (\varepsilon \alpha^j \beta k, \varepsilon \alpha^j) \in KK^{-1} \subset [-\beta M, \beta M] \times \{\alpha^{-N} \leq |a| \leq \alpha^N\}.$$

The left-hand side amounts to

$$(\eta\alpha^l\beta m, \eta\alpha^l)^{-1}(\varepsilon\alpha^j\beta k, \varepsilon\alpha^j) = (-\beta m, \eta\alpha^{-l})(\varepsilon\alpha^j\beta k, \varepsilon\alpha^j)$$
$$= (\eta\varepsilon\alpha^{j-l}\beta k - \beta m, \eta\varepsilon\alpha^{j-l}).$$

Now fix $k, j \in \mathbb{Z}$ and $\varepsilon = \pm 1$. The scale factor $\eta\varepsilon\alpha^{j-l}$ can only be contained in $\{\alpha^{-N} \leq |a| \leq \alpha^N\}$ if $l \in \{j - N, \ldots, j + N\}$. Similarly, if the number $\eta\varepsilon\alpha^{j-l}\beta k - \beta m$ is contained in $[-\beta M, \beta M]$, then

$$m \in [\eta\varepsilon\alpha^{j-l}k - M, \eta\varepsilon\alpha^{j-l}k + M] = \eta\varepsilon\alpha^{j-l}k + [-M, M].$$

Thus there are at most $(2N + 1)(2M + 1)$ pairs of values $(l, m) \in \mathbb{Z}^2$ and the two values $\eta = \pm 1$ for which

$$(\eta\alpha^l\beta m, \eta\alpha^l)^{-1}(\varepsilon\alpha^j\beta k, \varepsilon\alpha^j) \in KK^{-1}$$

can hold. This is independent of $k, j$ and $\varepsilon$, so it follows

$$\sup_{j,k\in\mathbb{Z},\varepsilon=\pm 1} |\{(m, l, \eta) \in \mathbb{Z}^2 \times \{\pm 1\} \mid (\varepsilon\alpha^j\beta k, \varepsilon\alpha^j)K \cap (\eta\alpha^l\beta m, \eta\alpha^l)K \neq \emptyset\}|$$
$$\leq 2(2N + 1)(2M + 1).$$

This proves that $\Lambda(\beta, \alpha)$ is relatively separated. □

The set of rectangles $\{A_{\beta,\alpha} \mid \alpha > 1, \beta > 0\}$ is a neighbourhood basis of the neutral element $(0, 1)$, that is, every neighbourhood $U$ of $(0, 1)$ contains a rectangle of the form $A_{\beta,\alpha}$. The lattices $\Lambda(\beta, \alpha)$ therefore become arbitrarily 'dense' for $\alpha \searrow 1$ and $\beta \searrow 0$.

**Remark 6.1.3.** It might look like a good idea to remove the negative scale factors and replace $\mathcal{Aff}$ by its subgroup $\mathcal{Aff}_+ = \mathbb{R} \rtimes (0, \infty)$, as in most practical applications only positive scale factors are used. However, the restriction of the wavelet representation $\pi$ to $\mathcal{Aff}_+$ is no longer irreducible, since positive scale factors are unable to relate positive and negative frequencies. For instance, any function $g \in L^2(\mathbb{R})$ with purely positive frequencies $\operatorname{supp} \hat{g} \subset (0, \infty)$ would generate a subspace

$$\widetilde{\mathcal{E}}_g = \overline{\operatorname{span}\{\pi(b, a)g \mid b \in \mathbb{R}, a > 0\}} \subset L^2(\mathbb{R})$$

of functions that again contain purely positive frequencies, meaning such $g$ cannot be cyclic with respect to $\mathcal{Aff}_+$. ◁

We now want to find a sufficient condition for a wavelet $\psi$ to belong to the set of $w_\rho$-atoms

$$\mathcal{B}_{w_\rho} = \left\{\psi \in \mathcal{A}_{w_\rho} \mid \operatorname{osc}_U(W_\psi\psi) \in L^1_{w_\rho}(\mathcal{Aff})\right\},$$

where $U$ is an arbitrary compact $(1, 0)$-neighbourhood (like a rectangle $A_{\beta,\alpha}$). Thus, we need a way to estimate the $L^1_{w_\rho}$-norm of the oscillation $\operatorname{osc}_U(W_\psi\psi)$. We proceed similarly as in section 4.1, where we derived a sufficient condition for the $L^1_{w_\rho}$-integrability of $W_\psi\psi$ itself.

First we prove an auxiliary lemma that helps to estimate the integral over local suprema of a function. It mirrors lemma 4.1.7.

**Lemma 6.1.4.** *Let $\psi \in L^1(\mathbb{R})$ differentiable such that $\psi'$ is integrable and its first moment converges absolutely. Let $U \subset \mathcal{A}\!f\!f$ be a compact neighbourhood of $(0,1)$. Then the integral*

$$\int_{\mathbb{R}} \sup_{(\delta,\tau)\in U} |\psi(\tau t + \delta)|\, dt$$

*is finite.*

*Proof.* We first note that the first moment of $\psi'$ vanishes. This is because the integral

$$\int_{-\infty}^{\infty} \psi'(x)\, dx = \lim_{y\to\infty} \psi(y) - \psi(-y)$$

converges absolutely and amounts to zero as $\psi$ is continuous and integrable. It follows in particular that

$$|\psi(x)| = \left|\int_{-\infty}^{x} \psi'(t)\, dt\right| = \left|-\int_{x}^{\infty} \psi'(t)\, dt\right| \leq \int_{|t|\geq |x|} |\psi'(t)|\, dt. \tag{6.2}$$

The set $U$ is compact in $\mathcal{A}\!f\!f$ and therefore contained in a double-rectangle $[-\beta,\beta] \times \{\alpha^{-1} \leq |a| \leq \alpha\}$ for some $\alpha > 1$ and $\beta > 0$. Thus, for $(\delta,\tau) \in U$ and $t \in \mathbb{R}$, the inequality $|\tau t + \delta| \geq |\tau t| - \beta \geq \alpha^{-1}|t| - \beta$ holds. With (6.2) it follows

$$\int_{\mathbb{R}} \sup_{(\delta,\tau)\in U} |\psi(\tau t + \delta)|\, dt$$

$$\leq \int_{\mathbb{R}} \sup_{(\delta,\tau)\in U} \int_{|x|\geq |\tau t+\delta|} |\psi'(x)|\, dx\, dt$$

$$\leq \int_{\mathbb{R}} \int_{|x|\geq \alpha^{-1}|t|-\beta} |\psi'(x)|\, dx\, dt$$

$$= \int_{-\alpha\beta}^{\alpha\beta} \int_{|x|\geq \alpha^{-1}|t|-\beta} |\psi'(x)|\, dx\, dt + \int_{|t|\geq \alpha\beta} \int_{|x|\geq \alpha^{-1}|t|-\beta} |\psi'(x)|\, dx\, dt,$$

where we used the last step to divide the integral over $\mathbb{R}$ into the two parts $|t| < \alpha\beta$ and $|t| \geq \alpha\beta$.

For $|t| < \alpha\beta$ it is $\alpha^{-1}|t| - \beta < 0$, hence the inner integral of the first summand runs over $\mathbb{R}$ entirely, independent of $t$. In the second summand, we change the order of integration. It follows

$$\int_{\mathbb{R}} \sup_{(\delta,\tau)\in U} |\psi(\tau t + \delta)|\, dt \leq \int_{-\alpha\beta}^{\alpha\beta} \int_{\mathbb{R}} |\psi'(x)|\, dx\, dt + \int_{\mathbb{R}} \int_{\alpha\beta \leq |t| \leq \alpha(|x|+\beta)} |\psi'(x)|\, dt\, dx$$

$$= 2\alpha\beta \|\psi'\|_{L^1(\mathbb{R})} + 2\alpha \int_{\mathbb{R}} |x||\psi'(x)|\, dx$$

$$< \infty,$$

completing the proof. $\square$

We can now formulate and proof our sufficient condition for $\psi \in \mathcal{B}_{w_\rho}$. This is very similar to the sufficient condition in proposition 4.1.10.

**Proposition 6.1.5.** *Suppose $\psi \in L^2(\mathbb{R})$ has $L \in \mathbb{N}$ vanishing moments and $L+1$ absolutely convergent moments. Suppose further that $\psi$ is $L$-times differentiable such that all derivatives $\psi^{(k)}$, $k = 0, \ldots, L$ are integrable, and that the $L$-th derivative satisfies the integrability condition*

$$\int_{\mathbb{R}} \sup_{(\delta,\tau) \in U} |\psi^{(L)}(\tau t + \delta)| \, dt < \infty \tag{6.3}$$

*for some compact neighbourhood $U \subset \mathcal{A}\!f\!f$ of $(0, 1)$.*

*Then $\mathrm{osc}_U(W_\psi \psi) \in L^1_{w_\rho}(\mathcal{A}\!f\!f)$ for $0 \leq \rho < L - 1/2$, and in particular $\psi \in \mathcal{B}_{w_\rho}$ for such values of $\rho$.*

*Proof.* The compact set $U \in \mathcal{A}\!f\!f$ is contained in a double-rectangle $[-\beta, \beta] \times \{\alpha^{-1} \leq |a| \leq \alpha\}$ for some $\alpha > 1$ and $\beta > 0$. This allows us to express the integrability condition (6.3) in terms of an oscillation condition

$$C_L = \int_{\mathbb{R}} \sup_{(\delta,\tau) \in U} \left| |\tau|^{1/2} \tau^L \psi^{(L)}(\tau t + \delta) - \psi^{(L)}(t) \right| dt \tag{6.4}$$

$$\leq \alpha^{L+1/2} \int_{\mathbb{R}} \sup_{(\delta,\tau) \in U} |\psi^{(L)}(\tau t + \delta)| \, dt + \left\| \psi^{(L)} \right\|_{L^1(\mathbb{R})}$$

$$< \infty,$$

which will be needed at a later point.

Similar to the proof of 4.1.10, the moment requirements of $\psi$ allow us to construct the family of antiderivatives

$$\psi^{(-k)}(x) = \int_{-\infty}^{x} \psi^{(-k+1)}(t) \, dt, \quad k = 1, \ldots, L.$$

All these antiderivatives are well-defined and integrable according to corollary 4.1.8. By induction over lemma 4.1.7, the function $\psi^{(-L+1)}$ has two absolutely convergent moments. We can therefore apply lemma 6.1.4 to $\psi^{(-L)}$ to obtain another oscillation condition

$$C_{-L} = \int_{\mathbb{R}} \sup_{(\delta,\tau) \in U} \left| |\tau|^{1/2} \tau^{-L} \psi^{(-L)}(\tau t + \delta) - \psi^{(-L)}(t) \right| dt \tag{6.5}$$

$$\leq \alpha^{L-1/2} \int_{\mathbb{R}} \sup_{(\delta,\tau) \in U} |\psi^{(-L)}(\tau t + \delta)| \, dt + \left\| \psi^{(-L)} \right\|_{L^1(\mathbb{R})}$$

$$< \infty.$$

We now consider the $U$-oscillation of $W_\psi \psi$ at the point $(b, a) \in \mathcal{A}\!f\!f$. It is defined as

$$\mathrm{osc}_U(W_\psi \psi)(b, a)$$
$$= \sup_{(\delta,\tau) \in U} |W_\psi \psi((\delta, \tau)(b, a)) - W_\psi \psi(b, a)|$$
$$= \sup_{(\delta,\tau) \in U} \left| |a\tau|^{-1/2} \int_{\mathbb{R}} \psi(t) \overline{\psi\left(\frac{t - \delta - \tau b}{\tau a}\right)} dt - |a|^{-1/2} \int_{\mathbb{R}} \psi(t) \overline{\psi\left(\frac{t - b}{a}\right)} dt \right|$$
$$= \sup_{(\delta,\tau) \in U} |a|^{-1/2} \left| \int_{\mathbb{R}} |\tau|^{1/2} \psi(\tau t + \delta) \overline{\psi\left(\frac{t - b}{a}\right)} dt - \int_{\mathbb{R}} \psi(t) \overline{\psi\left(\frac{t - b}{a}\right)} dt \right|$$
$$= \sup_{(\delta,\tau) \in U} |a|^{-1/2} \left| \int_{\mathbb{R}} \left( |\tau|^{1/2} \psi(\tau t + \delta) - \psi(t) \right) \overline{\psi\left(\frac{t - b}{a}\right)} dt \right|.$$

We apply the partial integration formula from lemma 4.1.9 $L$-times to the last integral. This gives us the two representations

$$\operatorname{osc}_U(W_\psi\psi)(b,a)$$

$$= \sup_{(\delta,\tau)\in U} |a|^{-1/2+L} \left| \int_\mathbb{R} \left(|\tau|^{1/2}\tau^L\psi^{(L)}(\tau t+\delta) - \psi^{(L)}(t)\right) \overline{\psi^{(-L)}\left(\frac{t-b}{a}\right)} dt \right| \quad (6.6)$$

$$= \sup_{(\delta,\tau)\in U} |a|^{-1/2-L} \left| \int_\mathbb{R} \left(|\tau|^{1/2}\tau^{-L}\psi^{(-L)}(\tau t+\delta) - \psi^{(-L)}(t)\right) \overline{\psi^{(L)}\left(\frac{t-b}{a}\right)} dt \right|. \quad (6.7)$$

Next we integrate the oscillation with respect to $db$. The first representation (6.6) leads to the estimate

$$\int_\mathbb{R} \operatorname{osc}_U(W_\psi\psi)(b,a)\, db$$

$$\leq |a|^{-1/2+L} \int_\mathbb{R} \sup_{(\delta,\tau)\in U} \int_\mathbb{R} \left||\tau|^{1/2}\tau^L\psi^{(L)}(\tau t+\delta) - \psi^{(L)}(t)\right| \left|\psi^{(-L)}\left(\frac{t-b}{a}\right)\right| dt\, db$$

$$\leq |a|^{-1/2+L} \int_\mathbb{R} \int_\mathbb{R} \sup_{(\delta,\tau)\in U} \left||\tau|^{1/2}\tau^L\psi^{(L)}(\tau t+\delta) - \psi^{(L)}(t)\right| \left|\psi^{(-L)}\left(\frac{t-b}{a}\right)\right| dt\, db$$

$$= |a|^{-1/2+L} \int_\mathbb{R} \sup_{(\delta,\tau)\in U} \left||\tau|^{1/2}\tau^L\psi^{(L)}(\tau t+\delta) - \psi^{(L)}(t)\right| \int_\mathbb{R} \left|\psi^{(-L)}\left(\frac{t-b}{a}\right)\right| db\, dt$$

$$= |a|^{1/2+L} \left\|\psi^{(-L)}\right\|_{L^1(\mathbb{R})} \int_\mathbb{R} \sup_{(\delta,\tau)\in U} \left||\tau|^{1/2}\tau^L\psi^{(L)}(\tau t+\delta) - \psi^{(L)}(t)\right| dt$$

$$= |a|^{1/2+L} \left\|\psi^{(-L)}\right\|_{L^1(\mathbb{R})} C_L, \quad (6.8)$$

where we abbreviated the last integral by the (finite) constant $C_L$ from (6.4). Similarly, the second representation (6.7) yields the estimate

$$\int_\mathbb{R} \operatorname{osc}_U(W_\psi\psi)(b,a)\, db$$

$$\leq |a|^{-1/2-L} \int_\mathbb{R} \sup_{(\delta,\tau)\in U} \int_\mathbb{R} \left||\tau|^{1/2}\tau^{-L}\psi^{(-L)}(\tau t+\delta) - \psi^{(-L)}(t)\right| \left|\psi^{(L)}\left(\frac{t-b}{a}\right)\right| dt\, db$$

$$\leq |a|^{1/2-L} \left\|\psi^{(L)}\right\|_{L^1(\mathbb{R})} C_{-L} \quad (6.9)$$

with the constant $C_{-L}$ from (6.5).

It remains to integrate with respect to $(|a|^\rho + |a|^{-\rho})\frac{da}{a^2}$ to obtain the $L^1_{w_\rho}$-norm of $\operatorname{osc}_U(W_\psi\psi)$. For that we separate $\mathbb{R}^*$ into the two parts $|a|\leq 1$ and $|a|>1$, and use the estimates (6.8) and (6.9) on those parts respectively. This gives us the inequality

$$\left\|\mathrm{osc}_U(W_\psi\psi)\right\|_{L^1_{w_\rho}(\mathcal{A}\!f\!f)}$$

$$= \int_\mathbb{R} \int_\mathbb{R} \sup_{(\delta,\tau)\in U} |W_\psi\psi((\delta,\tau)(b,a)) - W_\psi\psi(b,a)| \, db \, (|a|^\rho + |a|^{-\rho}) \frac{da}{a^2}$$

$$\leq \int_{|a|\leq 1} |a|^{1/2+L} \left\|\psi^{(-L)}\right\|_{L^1(\mathbb{R})} C_L (|a|^\rho + |a|^{-\rho}) \frac{da}{a^2}$$

$$+ \int_{|a|>1} |a|^{1/2-L} \left\|\psi^{(L)}\right\|_{L^1(\mathbb{R})} C_{-L} (|a|^\rho + |a|^{-\rho}) \frac{da}{a^2}$$

$$\leq C \left( \int_{|a|\leq 1} |a|^{\rho+L-3/2} + |a|^{-\rho+L-3/2} \, da + \int_{|a|>1} |a|^{\rho-L-3/2} + |a|^{-\rho-L-3/2} \, da \right)$$

with some constant $C < \infty$. The two remaining integrals converge exactly for $\rho < L - 1/2$ resp. $\rho < L + 1/2$, thus $\mathrm{osc}_U(W_\psi\psi) \in L^1_{w_\rho}(\mathcal{A}\!f\!f)$ if $\rho \in [0, L + 1/2)$.

Since $\psi$ fulfils all requirements of proposition 4.1.10, $W_\psi\psi \in L^1_{w_\rho}(\mathcal{A}\!f\!f)$ is also true, so it follows $\psi \in \mathcal{B}_{w_\rho}$ as claimed. $\square$

Using the well-spread families $\Lambda(\beta,\alpha)$ for $\alpha > 1$, $\beta > 0$ and the just proven sufficient condition for $\psi \in \mathcal{B}_{w_\rho}$, we are now able to give statements about Banach frames for the wavelet coorbit spaces $\mathcal{C}o^p_{m_s}$. For that we first recapitulate what we need and fix some abbreviating notation.

If $s \in \mathbb{R}$, then $w_\rho$ is a simultaneous control-weight of $m_s$ for $\rho \geq |s| + 1$. In order to use sufficient condition 6.1.5, we need a wavelet $\psi$ that satisfies all the stated conditions for some $L > \rho + 1/2 \geq |s| + 3/2$. For fixed $\alpha > 1$, $\beta > 0$ we discretize the weight $m_s(b,a) = |a|^{-s}$ with respect to the family $\Lambda(\beta,\alpha)$ by evaluating it at the lattice points. This yields

$$\widetilde{m}^{\beta,\alpha}_s(j,k,\varepsilon) = m_s(\varepsilon\alpha^j\beta k, \varepsilon\alpha^j) = \alpha^{-js}.$$

The functions contained in the wavelet family $(\pi(\lambda)\psi)_{\lambda\in\Lambda(\beta,\alpha)}$ can be written as

$$\psi^{\beta,\alpha}_{j,k,\varepsilon}(x) = \pi(\varepsilon\alpha^j\beta k, \varepsilon\alpha^j)\psi(x) = \alpha^{-j/2}\psi(\varepsilon\alpha^{-j}x - \beta k), \quad j,k \in \mathbb{Z}, \, \varepsilon = \pm 1.$$

**Theorem 6.1.6.** *Let $L > |s| + 3/2$ a natural number and $1 \leq p \leq \infty$. Suppose $\psi \in L^2(\mathbb{R})$ satisfies all requirements of proposition 6.1.5.*

*Then there is $\alpha > 1$ and $\beta > 0$ such that the family $(\psi^{\beta,\alpha}_{j,k,\varepsilon})_{j,k\in\mathbb{Z},\varepsilon=\pm 1}$ is a Banach frame for $\mathcal{C}o^p_{m_s}$ with respect to the sequence space*

$$\ell^p_{\widetilde{m}^{\beta,\alpha}_s}(\mathbb{Z}^2 \times \{\pm 1\}) = \left\{ (c_{j,k,\varepsilon})_{j,k\in\mathbb{Z},\varepsilon=\pm 1} \, \middle| \, \left( \sum_{\varepsilon=\pm 1} \sum_{j,k\in\mathbb{Z}} |c_{j,k,\varepsilon}|^p \alpha^{-jsp} \right)^{1/p} < \infty \right\}.$$

*In particular, there are frame constants $0 < A \leq B < \infty$ such that the inequalities*

$$A^p \|f\|^p_{\mathcal{C}o^p_{m_s}} \leq \sum_{\varepsilon=\pm 1} \sum_{j,k\in\mathbb{Z}} \alpha^{-jps} \left| \langle f, \psi^{\beta,\alpha}_{j,k,\varepsilon} \rangle_{\mathcal{R}_{w_\rho} \times \mathcal{H}^1_{w_\rho}} \right|^p \leq B^p \|f\|^p_{\mathcal{C}o^p_{m_s}}$$

*hold for all $f \in \mathcal{C}o^p_{m_s}$ (with modifications for $p = \infty$).*

*Proof.* Let $\rho = |s| + 1$. According to proposition 6.1.5 we have $\psi \in \mathcal{B}_{w_\rho}$, so $W_\psi \psi$ has a $w_\rho$-integrable oscillation. Since $W_\psi \psi$ is continuous, lemma 5.2.6 implies that there is a small neighbourhood $Q \subset \mathcal{A}\!f\!f$ of $(0, 1)$ such that the inequality

$$\|W_\psi \psi\|_{L^1_{w_\rho}(\mathcal{A}\!f\!f)} \|\mathrm{osc}_Q(W_\psi \psi)\|_{L^1_{w_\rho}(\mathcal{A}\!f\!f)} < 1$$

holds. Now there exist $\alpha > 1$ and $\beta > 0$ for which $A_{\beta,\alpha} \subset Q$ (with $A_{\beta,\alpha}$ as in (6.1)), and the above estimate still holds if $Q$ is replaced by $A_{\beta,\alpha}$. The family $\Lambda(\beta, \alpha)$ is $A_{\beta,\alpha}$-well-spread. Now, the statement follows from theorem 5.4.2. □

The sufficient condition 6.1.5 is particularly satisfied by the Schwartz functions with infinitely many vanishing moments $\mathcal{S}_0$ for any $\rho \geq 0$. If $\psi \in \mathcal{S}_0$, lemma 6.1.4 can be applied to $\psi^{(L)}$ for any $L \in \mathbb{N}$ to see that the integrability condition (6.3) is always satisfied, hence $\psi \in \mathcal{B}_{w_\rho}$ for all $\rho \geq 0$. Combining this with theorem 6.1.6 we see that for any Schwartz function with infinitely many vanishing moments and any $s \in \mathbb{R}$ there is a lattice $\Lambda(\beta, \alpha)$ that is dense enough such that $(\pi(\lambda)\psi)_{\lambda \in \Lambda(\beta,\alpha)}$ is a Banach frame for the coorbit space $\mathcal{C}o^p_{m_s} \cong \dot{B}^{s-1/2+1/p}_{p,p}$.

We can also use theorem 5.4.4 to obtain Banach frames that are extensions of tight wavelet frames for $L^2(\mathbb{R})$. Such wavelet frames are often constructed with respect to the lattice $(2^j k, 2^j)_{j,k \in \mathbb{Z}} \subset \mathcal{A}\!f\!f$, which uses only positive scale factors $2^j > 0$. We use the common notation $\psi_{j,k} = 2^{j/2} \psi(2^j \cdot + k) = \pi(-2^{-j}k, 2^{-j})\psi$. Note that the sign of $j$ (and $k$) is inverted compared to our previous notation.

**Theorem 6.1.7.** *Suppose $\psi \in L^2(\mathbb{R})$ satisfies all conditions from proposition 6.1.5 for some $L \in \mathbb{N}$, $L \geq 2$. Suppose further that the family $(\psi_{j,k})_{j,k \in \mathbb{Z}}$ is a tight frame for $L^2(\mathbb{R})$, that is*

$$c\|f\|_{L^2(\mathbb{R})} = \left( \sum_{j,k \in \mathbb{Z}} |\langle f, \psi_{j,k} \rangle_{L^2(\mathbb{R})}|^2 \right)^{1/2}$$

*for all $f \in L^2(\mathbb{R})$ and a constant $c > 0$.*

*Then $(\psi_{j,k})_{j,k \in \mathbb{Z}}$ is a Banach frame for all coorbit spaces $\mathcal{C}o^p_{m_s}$, $|s| < L - 3/2$ and $p \in 1 \leq p \leq \infty$, with respect to the sequence space*

$$\left\{ (c_{j,k})_{j,k \in \mathbb{Z}} \; \middle| \; \|(c_{j,k})_{j,k \in \mathbb{Z}}\| = \left( \sum_{j,k \in \mathbb{Z}} |c_{j,k}|^p \, 2^{jsp} \right)^{1/p} < \infty \right\}.$$

Tight wavelet frames can be constructed through (generalized) multiresolution analyses, c.f. [32, Ch. 2] or [3, Ch. 17]. There it is possible to find wavelets with certain prescribed properties like vanishing moments, smoothness or compact support. All three properties are simultaneously achievable, as was proven by Daubechies [6].

If $\psi$ is $L$-times differentiable and has compact support, the $L$th derivative automatically satisfies the integrability condition (6.3). Thus, if such $\psi$ also is $L$-times differentiable and has $L$ vanishing moments as well as $L + 1$ absolutely convergent moments, it is contained in $\mathcal{B}_{w_\rho}$ for $0 \leq \rho < L - 3/2$. Multiresolution analyses can therefore be used to obtain Banach frames for wavelet coorbit spaces.

## *6.2 Banach Gabor Frames for Modulation Spaces*

We now aim to construct Banach frames for the modulation spaces, i.e. the coorbit spaces associated to the short-time Fourier transform. Such frames are called *Banach Gabor frames*.

We have defined the modulation spaces $M_{r,s}^p$ in section 4.2 as the coorbit spaces that belong to the Schrödinger representation $\rho$ of the reduced Heisenberg group $\mathbb{H}_r^d$. The considered weight functions are

$$v_{r,s}(x,\omega,\tau) = v_{r,s}(x,\omega) = (1+|x|)^r(1+|\omega|)^s$$

for $r,s \geq 0$. These weights are their own control-weights. We have also proven that the space of analyzing vectors coincides with the space of test vectors $M_{r,s}^1$.

In section 4.2, we were able to replace all integrability conditions regarding the voice transform

$$\mathcal{V}_g f(x,\omega,\tau) = \overline{\tau} e^{\pi i x \omega} V_g f(x,\omega), \quad (x,\omega,\tau) \in \mathbb{H}_r^d$$

by integrability conditions regarding the short-time Fourier transform. This allowed us to describe the modulation spaces solely in terms of the short-time Fourier transform. We want to use the same principle again to construct Banach frames that only use time-frequency shifts without any phase factor $\tau \in S^1$. We do this by showing that certain Banach frames of the form $(\rho(x_j,\omega_j,\tau_k)g)_{j\in I, k\in J}$ can be simplified to a Banach frame $(M_{\omega_j} T_{x_j} g)_{j\in I}$. A similar approach is used in [2, Thm. 6.1] for atomic decompositions associated to projective representations.

**Lemma 6.2.1.** *Let $(x_j,\omega_j,\tau_k)_{j\in I, k\in J} \subset \mathbb{H}_r^d$ be a family with countable index set $I$ and finite index set $J = \{1,\ldots,N\}$. Suppose $g \in M_{r,s}^1$ such that $(\rho(x_j,\omega_j,\tau_k)g)_{j\in I, k\in J}$ is a Banach frame of $M_{r,s}^p$ with respect to the sequence space*

$$\ell_{\widetilde{v}_{r,s}}^p(I \times J) = \left\{ (c_{j,k})_{j\in I, k\in J} \;\middle|\; \sum_{j=1}^N \sum_{j\in I} |c_{j,k}|^p \, v_{r,s}(x_j,\omega_j)^p < \infty \right\}.$$

*Then $(M_{\omega_j} T_{x_j} g)_{j\in I}$ defines a Banach frame of $M_{r,s}^p$ with respect to the sequence space*

$$\ell_{\widetilde{v}_{r,s}}^p(I) = \left\{ (c_j)_{j\in I} \;\middle|\; \sum_{j\in I} |c_j|^p \, v_{r,s}(x_j,\omega_j)^p < \infty \right\}.$$

*Proof.* Let $f \in M_{r,s}^p$. The the equation $\rho(x,\omega,\tau) = \tau e^{-\pi i x \omega} M_\omega T_x$ implies

$$\left\| (\langle f, \rho(x_j,\omega_j,\tau_k)g \rangle_{\mathcal{R}_{v_{r,s}} \times M_{r,s}^1})_{j\in I, k\in J} \right\|_{\ell_{\widetilde{v}_{r,s}}^p(I\times J)}^p$$

$$= \sum_{j=1}^N \sum_{k\in I} \left| \langle f, \tau_k e^{-\pi i x_j \omega_j} M_{\omega_j} T_{x_j} g \rangle_{\mathcal{R}_{v_{r,s}} \times M_{r,s}^1} \right|^p v_{r,s}(x_j,\omega_j)^p$$

$$= \sum_{k=1}^N \sum_{j\in I} \left| \langle f, M_{\omega_j} T_{x_j} g \rangle_{\mathcal{R}_{v_{r,s}} \times M_{r,s}^1} \right|^p v_{r,s}(x_j,\omega_j)^p$$

$$= N \left\| (\langle f, M_{\omega_j} T_{x_j} g \rangle_{\mathcal{R}_{v_{r,s}} \times M_{r,s}^1})_{j\in I} \right\|_{\ell_{\widetilde{v}_{r,s}}^p(I)}^p.$$

Now the family $(\rho(x_j, \omega_j, \tau_k)g)_{j \in I, k \in J}$ satisfies the frame inequalities

$$A\|f\|_{M^p_{r,s}} \leq \left\|(\langle f, \rho(x_j, \omega_j, \tau_k)g\rangle_{\mathcal{R}_{v_{r,s}} \times M^1_{r,s}})_{j \in I, k \in J}\right\|_{\ell^p_{\widetilde{v}_{r,s}}(I \times J)} \leq B\|f\|_{M^p_{r,s}}$$

for some constants $0 < A \leq B < \infty$. The above computation shows that

$$A/N\|f\|_{M^p_{r,s}} \leq \left\|(\langle f, M_{\omega_j} T_{x_j} g\rangle_{\mathcal{R}_{v_{r,s}} \times M^1_{r,s}})_{j \in I}\right\|_{\ell^p_{\widetilde{v}_{r,s}}(I)} \leq B/N\|f\|_{M^p_{r,s}}$$

holds, so $(M_{\omega_j} T_{x_j} g)_{j \in I}$ has the frame bounds $0 < A/N \leq B/N < \infty$.

Let $R : \ell^p_{\widetilde{v}_{r,s}}(I \times J) \to M^p_{r,s}$ be the reconstruction operator of the Banach frame $(\rho(x_j, \omega_j, \tau_k)g)_{j \in I, k \in J}$. Then we define $\widetilde{R} : \ell^p_{\widetilde{v}_{r,s}}(I) \to M^p_{r,s}$ by

$$\widetilde{R}\left((c_j)_{j \in I}\right) = R\left((c_j \overline{\tau}_k e^{\pi i x_j \omega_j})_{j \in I, k \in J}\right).$$

This operator is bounded since the map

$$\ell^p_{\widetilde{v}_{r,s}}(I) \to \ell^p_{\widetilde{v}_{r,s}}(I \times J), \quad (c_j)_{j \in I} \mapsto (c_j \overline{\tau}_k e^{\pi i x_j \omega_j})_{j \in I, k \in J}$$

is bounded. It also defines a reconstruction operator for $(M_{\omega_j} T_{x_j} g)_{j \in I}$ as we have

$$\widetilde{R}\left((\langle f, M_{\omega_j} T_{x_j} g\rangle_{\mathcal{R}_{v_{r,s}} \times M^1_{r,s}})_{j \in I}\right) = R\left((\langle f, \tau_k e^{-\pi i x_j \omega_j} M_{\omega_j} T_{x_j} g\rangle_{\mathcal{R}_{v_{r,s}} \times M^1_{r,s}})_{j \in I, k \in J}\right)$$
$$= R\left((\langle f, \rho(x_j, \omega_j, \tau_k)g\rangle_{\mathcal{R}_{v_{r,s}} \times M^1_{r,s}})_{j \in I, k \in J}\right)$$
$$= f.$$

This concludes the proof. $\square$

Next we want to find a sufficient condition for a function $g \in L^2(\mathbb{R}^d)$ to be a $v_{r,s}$-atom. For that, we use the following proposition, which gives a simple answer to that question for general *IN-groups*. A group $\mathcal{G}$ is called an IN-group if there is a compact neighbourhood $U$ of the neutral element which is invariant under inner group automorphisms, i.e. which satisfies $Ux = xU$ for all $x \in \mathcal{G}$.

**Proposition 6.2.2.** *[16, Lemma 7.2] Let $\mathcal{G}$ be an IN-group and $\pi$ a w-integrable representation of $\mathcal{G}$ on $\mathcal{H}$ for a submultiplicative weight $w$. Then $\mathcal{A}_w = \mathcal{B}_w$.*

*Proof.* Let $F \in L^1_w(\mathcal{G})$ and $G \in W(L^1_w)$. We have already proven the inequality (5.12)

$$\left\|\mathrm{osc}_U(G * F)\right\|_{L^1_w} \leq \left\|\mathrm{osc}_U(G)\right\|_{L^1_w} \|F\|_{L^1_w}.$$

Now we can use the invariant compact $e$-neighbourhood $U$ to also show the corresponding inequality for $F * G$.

Observe first that since $G(xU) = G(Ux)$ for all $x \in \mathcal{G}$, we have

$$\operatorname*{ess\,sup}_{u \in U} |G(ux) - G(x)| = \operatorname*{ess\,sup}_{u \in U} |G(xu) - G(x)|.$$

Thus it follows

$$\begin{aligned}
\operatorname{osc}_U(F*G)(x) &= \operatorname*{ess\,sup}_{u\in U} |(F*G)(ux) - (F*G)(x)| \\
&= \operatorname*{ess\,sup}_{u\in U} \left| \int_{\mathcal{G}} F(y)\left(G(y^{-1}ux) - G(y^{-1}x)\right) dy \right| \\
&\leq \int_{\mathcal{G}} |F(y)| \operatorname*{ess\,sup}_{u\in U} \left|G(y^{-1}ux) - G(y^{-1}x)\right| dy \\
&\leq \int_{\mathcal{G}} |F(y)| \operatorname*{ess\,sup}_{u\in U} \left|G(uy^{-1}x) - G(y^{-1}x)\right| dy \\
&= (|F| * \operatorname{osc}_U(G))(x),
\end{aligned}$$

so we can apply Young inequality (2.5) to obtain

$$\left\|\operatorname{osc}_U(F*G)\right\|_{L^1_w} \leq \|F\|_{L^1_w} \left\|\operatorname{osc}_U(G)\right\|_{L^1_w}. \tag{6.10}$$

Now let $\phi \in \mathcal{B}_w \setminus \{0\}$ be admissible and $\psi \in \mathcal{A}_w$ be arbitrary. Then we can use the convolution relations from corollary 2.6.4 to write the kernel $\mathcal{V}_\psi \psi$ as

$$\mathcal{V}_\psi \psi = \mathcal{V}_\phi \psi * \mathcal{V}_\phi \phi * \mathcal{V}_\psi \phi.$$

The convolution inequalities (5.12) and (6.10) therefore yield

$$\left\|\operatorname{osc}_U(\mathcal{V}_\psi \psi)\right\|_{L^1_w} \leq \left\|\mathcal{V}_\phi \psi\right\|_{L^1_w} \left\|\operatorname{osc}_U(\mathcal{V}_\phi \phi)\right\|_{L^1_w} \left\|\mathcal{V}_\psi \phi\right\|_{L^1_w} < \infty.$$

We see that $\mathcal{V}_\psi \psi$ has a $w$-integrable oscillation, so $\psi \in \mathcal{B}_w$. This shows that $\mathcal{B}_w = \mathcal{A}_w$. □

The reduced Heisenberg group is an IN-group, as for instance $[-1,1]^d \times [-1,1]^d \times S^1$ is an invariant $(0,0,1)$-neighbourhood. Thus all $v_{r,s}$-integrable vectors are already $v_{r,s}$-atoms, that is $\mathcal{B}_{v_{r,s}} = M^1_{r,s}$.

The fact that the space of test vectors coincides with the space of atoms is an important observation. The space of Schwartz functions $\mathcal{S}$ is contained in $M^1_{r,s}$ for all $r,s \geq 0$, so this space $\mathcal{B}_{v_{r,s}} = M^1_{r,s}$ of 'nice window functions' extends the often used Schwartz space, and makes it possible to replace it in a lot of contexts. This particularly underlines the importance of Feichtinger's algebra $S_0 = M^1_{0,0}$, which shares many properties with $\mathcal{S}$ but has the advantage of being a Banach space [22, sec. 12.1; 25].

Before discussing special well-spread families in $\mathbb{R}^d \times \mathbb{R}^d$, we write down the Banach frame statement in more generality. If we have a well-spread family $(x_j, \omega_j)_{j\in I} \subset \mathbb{R}^d \times \mathbb{R}^d$, we use the discretized weight

$$\widetilde{v}_{r,s}(j) = v_{r,s}(x_j, \omega_j) = (1+|x_j|)^r (1+|\omega_j|)^s, \quad j \in I.$$

The frame then consists of the functions

$$g_j(t) = M_{\omega_j} T_{x_j} g(t) = e^{2\pi i \omega_j t} g(t - x_j), \quad t \in \mathbb{R}^d.$$

**Proposition 6.2.3.** *Let $g \in M^1_{r,s}$ be admissible (i.e. $\|g\|_{L^2} = 1$) and $1 \leq p \leq \infty$. Then there is a compact $(0,0)$-neighbourhood $U \subset \mathbb{R}^d \times \mathbb{R}^d$ such that for all $U$-well-spread families $(x_j, \omega_j)_{j\in I}$ the functions $(g_j)_{j\in I}$ define a Banach frame for $M^p_{r,s}$ with respect to the sequence space $\ell^p_{\widetilde{v}_{r,s}}(I)$.*

*Proof.* First there is a compact neighbourhood $Q \subset \mathbb{H}_r^d$ of $(0,0,1)$ such that

$$\|\mathcal{V}_g g\|_{L^1_{v_{r,s}}(\mathbb{H}_r^d)} \|\mathrm{osc}_Q(\mathcal{V}_g g)\|_{L^1_{v_{r,s}}(\mathbb{H}_r^d)} < 1 \tag{6.11}$$

because $\mathcal{V}_g g$ has a $v_{r,s}$-integrable oscillation and is continuous, see lemma 5.2.6. We can assume that $Q$ has the form $U \times W$ with some compact $(0,0)$-neighbourhood $U \subset \mathbb{R}^d \times \mathbb{R}^d$ and some symmetric 1-neighbourhood

$$1 \in W = \{e^{2\pi i \theta} \mid \theta \in [-\varepsilon, \varepsilon]\} \subset S^1, \quad \varepsilon > 0,$$

as products of this form define a neighbourhood basis of $(0,0,1)$.

Now there are finitely many $\tau_1, \ldots, \tau_N \in S^1$ such that $S^1$ is covered by $\cup_{k=1}^N \tau_k W$. Let $(x_j, \omega_j)_{j \in I}$ be an arbitrary $U$-well-spread family in $\mathbb{R}^d \times \mathbb{R}^d$. Then $(x_j, \omega_j, \tau_k)_{j \in I, k \in J}$ with $J = \{1, \ldots, N\}$ is $U \times W$-well-spread in $\mathbb{H}_r^d$. Theorem 5.4.2 therefore implies that

$$(\rho(x_j, \omega_j, \tau_k) g)_{j \in I, k \in J}$$

is a Banach frame for $M_{r,s}^p$ with respect to the sequence space $\ell^p_{\widetilde{v}_{r,s}}(I \times J)$. By 6.2.1, this Banach frame reduces to the Gabor Banach frame $(M_{\omega_j} T_{x_j} g)_{j \in I}$ with respect to the sequence space $\ell^p_{\widetilde{v}_{r,s}}(I)$. $\square$

In time-frequency analysis, lattices are commonly used as well-spread families in $\mathbb{R}^d \times \mathbb{R}^d$. These are subgroups of $(\mathbb{R}^{2d}, +)$ that can be written as $\Lambda = A\mathbb{Z}^{2d}$, where $A$ is an invertible $2d \times 2d$ matrix. Such lattices are always relatively separated as we will see in the next proof. To achieve the required density, we scale the lattice by a factor $c > 0$, $c \searrow 0$.

**Proposition 6.2.4.** *Let $g \in M_{r,s}^1$ be admissible and $A \in \mathrm{GL}(\mathbb{R}^{2d})$. Then there exists a $c > 0$ such that the lattice $\Lambda = cA\mathbb{Z}^{2d}$ induces the Banach Gabor frame $(M_\omega T_x g)_{(x,\omega) \in \Lambda}$ for $M_{r,s}^p$ with respect to the sequence space $\ell^p_{\widetilde{v}_{r,s}}(I)$.*

*Proof.* We only need to show that to any compact $(0,0)$-neighbourhood $U \subset \mathbb{R}^d \times \mathbb{R}^d$ there exists a $c > 0$ such that $\Lambda = cA\mathbb{Z}^{2d}$ is $U$-well-spread. The statement then follows from proposition 6.2.3.

To show that $\Lambda$ is relatively separated, we assume without loss of generality that $c = 1$. Let $K \subset \mathbb{R}^{2d}$ be compact. We need to prove that there is a constant $C_K > 0$ such that for any $m \in \mathbb{Z}^{2d}$

$$|\{n \in \mathbb{Z}^{2d} \mid (Am + K) \cap (An + K) \neq \emptyset\}| \leq C_K$$

holds.

If $m$ is fixed and $(Am + K) \cap (Mn + K) \neq \emptyset$, it follows that $An \in Mm + K - K$ and therefore $n \in m + A^{-1}(K - K)$, where $K - K = \{x - y \mid x, y \in K\}$. Since $m + A^{-1}(K - K)$ is compact and bounded, the intersection $m + A^{-1}(K - K) \cap \mathbb{Z}^{2d}$ contains at most $C_K$ elements for some $C_K < \infty$. The cardinality of this intersection is invariant under translations from $\mathbb{Z}^{2d}$, so $C_K$ is indeed independent of $m$. Thus $\Lambda$ is relatively separated.

For the $U$-density we assume that $U \subset \mathbb{R}^{2d}$ is a compact 0-neighbourhood. Then $A^{-1}U$ is a 0-neighbourhood as well, so there is a constant $c > 0$ such that $[-c/2, c/2]^{2d} \subset A^{-1}U$. Thus we have

$$cA[-1/2, 1/2]^{2d} \subset U$$

and therefore

$$\mathbb{R}^{2d} = cA\mathbb{R}^{2d} = \bigcup_{m \in \mathbb{Z}^{2d}} cA\left(m + [-1/2, 1/2]^{2d}\right) \subset \bigcup_{m \in \mathbb{Z}^{2d}} cAm + U.$$

This shows that the lattice $\Lambda = cA\mathbb{Z}^{2d}$ is $U$-dense. $\square$

In the above proposition, the constant $c$ does not depend on the exponent $p$. Thus, for any admissible function $g \in M^1_{r,s}$ (that is $g \in M^1_{r,s}$ with $\|g\|_{L^2} = 1$), there is a lattice $\Lambda \subset \mathbb{R}^d \times \mathbb{R}^d$ that is dense enough so that $(M_\omega T_x g)_{(x,\omega) \in \Lambda}$ is a Banach Gabor frame for all modulation spaces $M^p_{r,s}$, $1 \leq p \leq \infty$.

Often lattices of the form $\Lambda = a\mathbb{Z}^d \times b\mathbb{Z}^d$ for some $a, b > 0$ are used. With the above result we see that for any admissible $g \in M^1_{r,s}$ and small enough constants $\alpha, \beta > 0$ (depending on $g$) the Gabor system

$$\mathcal{G}(g; \alpha, \beta) = \{M_{\alpha k} T_{\beta l} g \mid k, l \in \mathbb{Z}\}$$

is a Banach Gabor frame for $M^p_{r,s}$ for all $1 \leq p \leq \infty$.

---

Table of Notation

---

| | |
|---|---|
| $\mathcal{A}ff$ | affine group (p. 46) |
| $\mathcal{A}_w$ | analyzing ($w$-integrable) vectors (p. 30) |
| $\mathcal{B}_w$ | $w$-atoms (p. 82) |
| $\mathcal{C}o^p_m$ | coorbit spaces (p. 40) |
| $F^\vee, F^\triangledown$ | involutions: $F^\vee(x) = F(x^{-1})$ and $F^\triangledown(x) = \overline{F(x^{-1})}$ |
| $F * G$ | convolution (p. 9) |
| $\widehat{f}, \mathcal{F}f$ | Fourier transform: $\widehat{f}(\omega) = \int_{\mathbb{R}^d} f(t) e^{-2\pi i \omega t}\, dt$ |
| $\check{f}$ | inverse Fourier transform |
| $\mathcal{G}$ | locally compact group (p. 5) |
| $\mathcal{H}$ | Hilbert space |
| $\mathcal{H}^1_w$ | test vectors (p. 32) |
| $\mathbb{H}^d_r$ | reduced Heisenberg group (p. 57) |
| $L^p_m(\Omega), \ell^p_{\widetilde{m}}(I)$ | weighted $L^p$-space (p. 9) |
| $\mathcal{M}^p_m$ | reproducing kernel space (pp. 42, 76) |
| $\mu, d\mu$ | Haar measure (p. 6) |
| $\mathrm{osc}_U(G)$ | oscillation (p. 74) |
| $\mathcal{R}_w$ | reservoir (p. 36) |
| $T^*$ | adjoint operator on Hilbert spaces |
| $T'$ | adjoint operator between dual spaces |
| $T^\sim$ | adjoint operator between antidual spaces |
| $V'$ | dual space |
| $V^\sim$ | antidual space |
| $\mathcal{V}_g f$ | (extended) voice transform (pp. 21, 37) |
| $V_g f$ | short-time Fourier transform (pp. 22, 57) |
| $W_g f$ | wavelet transform (pp. 21, 46) |
| $W(L^1_w)$ | functions with $w$-integrable oscillation (p. 75) |